%% file: main.tex
\documentclass[12pt,final,CPage]{ufthesis} 

\include{packages}

\SetFullName{Marly Gotti}%
\SetThesisType{Dissertation}
\SetDegreeType {Doctor of Philosophy}
\SetGradMonth{December}%
\SetGradYear{2019}%
\SetDepartment{Mathematics}%
\SetChair{Peter Sin}%
\SetTitle{Atomicity and Factorization of Puiseux Monoids}

\include{usersetcommands}

\begin{document}



\bibliographystyle{bst/abbrv}


\maketitle 
\makecopyright


\dedication{\input{tex/dedication}} 


\include{tex/acknowledgements} 


\include{tex/TOC} 




\addtocontents{toc}{\protect\addvspace{10pt}\noindent{CHAPTER}\protect\hfill\par}{}
\phantomsection
\include{tex/abstract} 



\include{tex/ch0}

\include{tex/ch1}
\include{tex/ch2}
\include{tex/ch3}
\include{tex/ch4}
\include{tex/ch5}

\include{tex/ch6}








\bibliography{bib/references}

\include{bio}


\end{document}

%% file: packages.tex
\usepackage{graphicx}
\usepackage{amsmath}
\usepackage{amsthm}
\usepackage{algpseudocode}
\usepackage{tabularx}
\usepackage{url}
\usepackage[letterpaper,hmargin=1in,vmargin=1in]{geometry}
\usepackage{lscape}
\usepackage{longtable}
\usepackage{amsfonts}
\usepackage{amssymb}
\usepackage{subfigure}
\usepackage{rotating}
\usepackage{calc}
\usepackage{setspace}
\usepackage{ufenumerate}
\usepackage{latexsym}
\usepackage{epsf}
\usepackage{epsfig}
\usepackage{euscript}
\usepackage{stmaryrd}
\usepackage[table]{xcolor}
\usepackage{amscd}
\usepackage[format=hang,justification=raggedright,singlelinecheck=0,labelsep=period]{caption}
\usepackage[numbers,sort&compress]{natbib}
\usepackage{enumitem} 

\definecolor{bleudefrance}{RGB}{205,205,250}
\definecolor{blizzardblue}{RGB}{235,235,250}

\usepackage[hyperfootnotes=false]{hyperref}
\hypersetup{colorlinks=true,linkcolor=blue,anchorcolor=blue,citecolor=blue,filecolor=blue,urlcolor=blue,bookmarksnumbered=true,pdfview=FitB} %

\usepackage{lipsum}  

%% file: usersetcommands.tex
\hyphenchar\font=-1

\newlength{\widthOfItem}

\renewenvironment{itemize}{%
	\begin{Itemize}
		\setlength{\itemsep}{0.5\baselineskip}
		\setlength{\labelwidth}{2em}
		\setlength{\listparindent}{.32in}%
		\setlength{\leftmargin}{.32in}
		\setlength{\rightmargin}{0in}
		\settowidth{\widthOfItem}{\labelitemi}
		\setlength{\labelsep}{\leftmargin-\widthOfItem}
		
		\singlespacing}{%
	\end{Itemize}}

\newtheorem*{maintheorem*}{Main Theorem}
\newtheorem{theorem}{Theorem}[section]
\newtheorem{proposition}[theorem]{Proposition}
\newtheorem{lemma}[theorem]{Lemma}
\newtheorem{corollary}[theorem]{Corollary}
\theoremstyle{definition}
\newtheorem{definition}[theorem]{Definition}
\newtheorem{remark}[theorem]{Remark}
\newtheorem{example}[theorem]{Example}

\newcommand{\nn}{\mathbb{N}}
\newcommand{\pp}{\mathbb{P}}
\newcommand{\qq}{\mathbb{Q}}
\newcommand{\rr}{\mathbb{R}}
\newcommand{\zz}{\mathbb{Z}}

\newcommand\pval{\mathsf{v}_p}
\newcommand{\gp}{\textsf{gp}}

\newcommand{\supp}{\mathsf{Supp}}
\newcommand\lcm{\text{lcm}}



%% file: tex/dedication.tex
\vspace{-20 mm}
To my grandmother \\


%% file: tex/acknowledgements.tex
\acknowledge{%

I would like to thank professors Philip Boyland, Yunmei Chen, Kevin Keating, Jennifer Rothschild, and Peter Sin for kindly agreeing to be part of my PhD committee. In particular, I would like to thank my advisor, Peter Sin, for his counseling and the many hours spent discussing interesting mathematics.

I am grateful to my coauthors Scott Chapman, Felix Gotti, and Harold Polo for the fruitful discoveries and enjoyable time together. Thank you, Scott, for introducing me to research in PURE Math 2013 and for being my research mentor since then. Thank you, Felix, for introducing me to atomicity and Puiseux monoids: this dissertation would not exist without your constant and effective guidance.

In addition, I would like to thank Blanqui, Dayron, Juandingo, Harold, and Meli for the many encouraging words and for our friendship. Thank you, Dayron, Meli, and Harold for the needed distraction, and thank you, Blanqui and Juandingo, for taking care of me when I lacked the time.

I would like to thank my parents, Josefa Marlen Marcos Vazquez and Carlos Armando Corrales Leon, for the lifelong emotional support and for the constant love during all these years of distant graduate studies. My deepest gratitude goes to my lovely grandmother Azucena Margarita Vazquez de la Garza not only for teaching me the first mathematics I ever knew, but also for being a source of inspiration at many times of my life. Finally, I am immensely grateful to my husband, Felix Gotti, for encouraging me to be the best version of myself, for teaching me how to navigate life, and for giving me his unconditional love.

}

%% file: tex/TOC.tex
\pdfbookmark[0]{TABLE OF CONTENTS}{tableofcontents}
\tableofcontents %
\listoftables %
\listoffigures %


%% file: tex/abstract.tex
\begin{abstract}

A commutative and cancellative monoid (or an integral domain) is called atomic if each non-invertible element can be expressed as a product of irreducibles. Many algebraic properties of monoids/domains are determined by their atomic structure. For instance, a ring of algebraic integers has class group of size at most two if and only if it is half-factorial (i.e., the lengths of any two irreducible factorizations of an element are equal). Most atomic monoids are not unique factorization monoids (UFMs). In factorization theory one studies how far is an atomic monoid from being a UFM. During the last four decades, the factorization theory of many classes of atomic monoids/domains, including numerical/affine monoids, Krull monoids, Dedekind domains, and Noetherian domains, have been systematically investigated.

A Puiseux monoid is an additive submonoid of the nonnegative cone of rational numbers. Although Puiseux monoids are torsion-free rank-one monoids, their atomic structure is rich and highly complex. For this reason, they have been important objects to construct crucial examples in commutative algebra and factorization theory. In 1974 Anne Grams used a Puiseux monoid to construct the first example of an atomic domain not satisfying the ACCP, disproving Cohn's conjecture that every atomic domain satisfies the ACCP. Even recently, Jim Coykendall and Felix Gotti have used Puiseux monoids to construct the first atomic monoids with monoid algebras (over a field) that are not atomic, answering a question posed by Robert Gilmer back in the 1980s.

This dissertation is focused on the investigation of the atomic structure and factorization theory of Puiseux monoids. Here we established various sufficient conditions for a Puiseux monoid to be atomic (or satisfy the ACCP). We do the same for two of the most important atomic properties: the finite-factorization property and the bounded-factorization property. Then we compare these four atomic properties in the context of Puiseux monoids. This leads us to construct and study several classes of Puiseux monoids with distinct atomic structure. Our investigation provides sufficient evidence to believe that the class of Puiseux monoids is the simplest class with enough complexity to find monoids satisfying almost every fundamental atomic behavior.

\end{abstract}

%% file: tex/ch0.tex
\chapter{Introduction} \label{chap:intro}

Factorization theory studies the phenomenon of non-unique factorizations into irreducibles in commutative cancellative monoids and integral domains. Factorization theory originated from commutative algebra and algebraic number theory, most of its initial motivation was the study of factorization into primes in ring of integers, Dedekind domains, and Krull domains. During the last four decades factorization theory has become an autonomous field and has been actively investigated in connection with other areas, including number theory, combinatorics, and convex geometry. The primary goal of factorization theory is to measure how far an atomic monoid or an integral domain is from being factorial or half-factorial (i.e., the number of irreducible factors in any two factorizations of a given element are the same).

The origin of factorization theory lies in algebraic number theory, and one of the primary motivations was the fact that the ring of integers $\mathcal{O}_K$ of an algebraic number field $K$ usually fails to be a UFD. Consider for example the ring of integers $\mathbb{Z}[\sqrt{-5}]$. We can write $6$ in $\mathbb{Z}[\sqrt{-5}]$ as
\[
6 = 2 \cdot 3 = (1 - \sqrt{-5})(1 + \sqrt{-5}),
\]
and the elements $2, 3, 1 - \sqrt{-5}$ and $1 + \sqrt{-5}$ are non-associate irreducible elements in $\mathbb{Z}[\sqrt{-5}]$ (this example is surveyed in one of my recent papers, \cite{CGG19a}). Throughout history some theories have been developed to understand the phenomenon of non-unique factorizations, including C. F. Gauss's theory of binary quadratic forms for quadratic fields, L. Kronecker's divisor theory, and R. Dedekind's ideal theory. In the mid-twentieth century, L. Carlitz characterized the half-factorial rings of integers in terms of their class number~\cite{lC60}, and W. Narkiewicz began a systematic study of non-unique factorizations on ring of integers (see~\cite{wN79} and references therein). On the other hand, R. Gilmer studied factorization properties of more general integral domains~\cite{dA06}. During the last four decades and motivated by the work of Gilmer, many authors have influenced the development of factorization theory in general integral domains (see \cite{AF11} and references therein). Since a large number of factorization properties of integral domains do not depend on the domains's additive structure, to investigate such properties it often suffices to focus on the corresponding multiplicative monoids. As a result, several techniques to measure the non-uniqueness of factorizations have been systematically developed and abstracted to the context of commutative cancellative monoids. Many factorization invariants and arithmetic statistics have been introduced, including the set of lengths, the union of sets of lengths, the elasticity, the catenary degree, and the tame degree.

In this thesis, we present results on the factorization theory and atomic structure of a class of commutative and cancellative monoids called Puiseux monoids (i.e, additive submonoids of $\mathbb{Q}_{\ge 0}$). Although the atomicity of Puiseux monoids has earned attention only in the last few years (see~\cite{CGG19,fG19a} and references therein), since the 1970s Puiseux monoids have been crucial in the construction of numerous examples in commutative ring theory. Back in 1974, A.~Grams~\cite{aG74} used an atomic Puiseux monoid as the main ingredient to construct the first example of an atomic integral domain that does not satisfy the ACCP, and thus she refuted P. Cohn's assumption that every atomic integral domain satisfies the ACCP. In addition, in~\cite{AAZ90}, A.~Anderson et al. appealed to Puiseux monoids to construct various needed examples of integral domains satisfying certain prescribed properties. More recently, Puiseux monoids have played an important role in~\cite{CG19}, where J. Coykendall and F. Gotti partially answered a question on the atomicity of monoid rings posed by R.~Gilmer back in the 1980s (see~\cite[page 189]{rG84}).

Puiseux monoids have also been important in factorization theory. For instance, the class of Puiseux monoids comprises the first (and only) example known so far of primary atomic monoids with irrational elasticity (this class was found in ~\cite[Section~4]{GGT19} via~\cite[Theorem~3.2]{GO19}). A Puiseux monoid is a suitable additive structure containing simultaneously several copies of numerical monoids independently generated. This fact has been harnessed by A. Geroldinger and W. Schmid to achieve a nice realization theorem for the sets of lengths of numerical monoids~\cite[Theorem~3.3]{GS18}. In~\cite{fG18} Puiseux monoids were studied in connection with Krull monoids and transfer homomorphisms. In addition, Puiseux monoids have been recently studied in~\cite{BG20} in connection to factorizations of upper triangular matrices. Finally, some connections between Puiseux monoids and music theory have been recently highlighted by M. Bras-Amoros in the Monthly article~\cite{mBA19}. A brief survey on Puiseux monoids can be found in~\cite{CGG19b}.

This thesis is organized as follows. Chapter 1 provides the background information and sets up the notation needed to study the atomic structure of Puiseux monoids given in Chapter 2 and Chapter 3. In these two chapters, the atomic structure of some families of Puiseux monoids is fully characterized (Proposition \ref{prop:atomic classification of multiplicative cyclic Puiseux monoids}, Corollary \ref{cor:prime reciprocal PM are atomic}, Proposition \ref{prop:atomicity p-adic}). In addition, the chain of implications \ref{eq:monoid atomicity taxonomy} is shown not to be reversible by providing results and examples in the realm of Puiseux monoids (Corollary \ref{cor:atomic and not ACCP}, Theorem \ref{thm:a class of ACCP monoids}, Theorem \ref{thm:BF sufficient condition}).

Furthermore, Chapters 4, 5, and 6 explore the factorization invariants of Puiseux monoids and Puiseux algebras. Sets of lengths, the union of set of lengths, the elasticity, $k$-elasticities, among other factorization invariants are analyzed (Theorem \ref{thm:sets of lengths}, Proposition \ref{prop:union of sets of lengths: infinite case}, Theorem \ref{theo:sufficient conditions for finite elasticity in ppm}). In addition, the connection between molecules and atomicity in Puiseux monoids and Puiseux algebras is investigated (Theorem \ref{thm:F(M) equals A(M)}, Theorem \ref{thm:U-UFD Puiseux algebras}).

%% file: tex/ch1.tex
\chapter{Algebraic Background of Puiseux Monoids} \label{chap:algebraic background}

\section{Preliminaries}
\label{sec:prelim}

In this section we introduce most of the relevant concepts on commutative monoids and factorization theory required to follow our exposition. General references for background information can be found in~\cite{pG01} for commutative monoids and in~\cite{GH06} for atomic monoids and factorization theory.
\smallskip

\subsection{General Notation}
We let $\nn := \{1,2,\dots\}$ denote the set of positive integers and set $\nn_0 := \nn \cup \{0\}$. In addition, we let $\pp$ denote the set of all prime numbers. For $X \subseteq \rr$ and $r \in \rr$, we set
\[
	X_{\ge r} := \{x \in X : x \ge r\}
\]
and we use the notations $X_{> r}, X_{\le r}$, and $X_{< r}$ in a similar manner. If $q \in \qq_{> 0}$, then we call the unique $n,d \in \nn$ such that $q = n/d$ and $\gcd(n,d)=1$ the \emph{numerator} and \emph{denominator} of $q$ and denote them by $\mathsf{n}(q)$ and $\mathsf{d}(q)$, respectively. Finally, for $Q \subseteq \qq_{>0}$, we set
\[
	\mathsf{n}(Q) := \{\mathsf{n}(q) : q \in Q\} \quad \text{ and } \quad \mathsf{d}(Q) := \{\mathsf{d}(q) : q \in Q\}.
\]
\smallskip

\subsection{Commutative Monoids} Throughout this thesis, the term \emph{monoid} stands for a commutative and cancellative semigroup with identity. Unless we specify otherwise, monoids are written additively. Let $M$ be a monoid. We let $M^\bullet$ denote the set of nonzero elements of $M$ while we let $U(M)$ denote the set of invertible elements of~$M$. When $M^\bullet = \emptyset$ we say that $M$ is \emph{trivial} and when $U(M) = \{0\}$ we say that $M$ is \emph{reduced}. 

For $S \subseteq M$, we let $\langle S \rangle$ denote the smallest (under inclusion) submonoid of $M$ containing $S$, and we call it the submonoid of $M$ \emph{generated} by $S$. The monoid $M$ is \emph{finitely generated} if $M$ can be generated by a finite set. An element $a \in M \setminus U(M)$ is an \emph{atom} provided that the equality $a = x+y$ for $x,y \in M$ implies that either $x \in U(M)$ or $y \in U(M)$. The set of atoms of $M$ is denoted by $\mathcal{A}(M)$.
\smallskip

\begin{definition}
	A monoid $M$ is \emph{atomic} if each element of $M \setminus U(M)$ can be expressed as a sum of atoms, and $M$ is \emph{antimatter} if $\mathcal{A}(M) = \emptyset$.
\end{definition}
 Every finitely generated monoid is atomic~\cite[Proposition~2.7.8(4)]{GH06}, while it immediately follows that every abelian group is an antimatter monoid.

A subset $I$ of $M$ is an \emph{ideal} of $M$ if $I + M = I$ (or, equivalently, $I + M \subseteq I$). An ideal $I$ is \emph{principal} if $I = x + M$ for some $x \in M$, and $M$ satisfies the \emph{ascending chain condition on principal ideals} (or \emph{ACCP}) provided that every increasing sequence of principal ideals of $M$ eventually stabilizes. It is well known that every monoid satisfying the ACCP must be atomic \cite[Proposition~1.1.4]{GH06}.

An equivalence relation $\rho \subseteq M \times M$ is a \emph{congruence} if it is compatible with the operation of the monoid $M$, i.e., for all $x,y,z \in M$ with $(x,y) \in \rho$ it follows that $(z + x, z + y) \in \rho$. It can be readily verified that the set $M/\rho$ consisting of the equivalence classes of a congruence $\rho$ is a commutative semigroup with identity. For $x,y \in M$, we say that $x$ \emph{divides} $y$ \emph{in} $M$ and write $x \mid_M y$ provided that $x + x' = y$ for some $x' \in M$. Two elements $x,y \in M$ are \emph{associates} if $y = u + x$ for some $u \in U(M)$. Being associates defines a congruence on $M$ whose semigroup of classes is a reduced monoid, which we denote by $M_{\text{red}}$. Observe that $M$ is reduced if and only if $M_{\text{red}} = M$.


The \emph{Grothendieck group} $\gp(M)$ of $M$ is the abelian group (unique up to isomorphism) satisfying that any abelian group containing a homomorphic image of $M$ also contains a homomorphic image of $\gp(M)$. The \emph{rank} of a monoid $M$ is the rank of the $\zz$-module $\gp(M)$ or, equivalently, the dimension of the $\qq$-vector space $\qq \otimes_\zz \gp(M)$. The monoid $M$ is \emph{torsion-free} if for all $x,y \in M$ and $n \in \nn$, the equality $nx = ny$ implies that $x=y$. Clearly, $M$ is a torsion-free monoid if and only if $\gp(M)$ is a torsion-free group.

A \emph{numerical monoid} is a submonoid $N$ of $(\nn_0,+)$ satisfying that $|\nn_0 \setminus N| < \infty$. If $N \neq \nn_0$, then $\max (\nn_0 \! \setminus \! N)$ is the \emph{Frobenius number} of $N$. Numerical monoids are finitely generated and, therefore, atomic with finitely many atoms. The \emph{embedding dimension} of $N$ is the cardinality of $\mathcal{A}(N)$. For an introduction to numerical monoids, see \cite{GR09}, and for some of their many applications, see~\cite{AG16a}.
\smallskip

\subsection{Factorizations} A multiplicative monoid $F$ is called \emph{free} with basis~$P$ if every element $x \in F$ can be written uniquely in the form
\[
	x = \prod_{p \in P} p^{\pval(x)},
\]
where $\pval(x) \in \nn_0$ and $\pval(x) > 0$ only for finitely many elements $p \in P$. The monoid $F$ is determined by $P$ up to isomorphism. By the Fundamental Theorem of Arithmetic, the multiplicative monoid $\nn$ is free on the set of prime numbers. In this case, we can extend $\pval$ to $\qq_{\ge 0}$ as follows. For $r \in \qq_{> 0}$ let $\pval(r) := \pval(\mathsf{n}(r)) - \pval(\mathsf{d}(r))$ and set $\pval(0) = \infty$. The map $\pval$ is called the $p$-\emph{adic valuation} on $\qq$. It is not hard to verify that $\pval$ is \emph{semi-additive}, i.e.,
\begin{equation*}
\pval(r + s) \ge \min\{\pval(r), \pval(s) \} \text{ for all } \ r,s \in \mathbb{Q}_{\ge 0}. \label{eq:semiadditivity of valuations}
\end{equation*}
\smallskip

\begin{definition}
	Let $M$ be a reduced monoid. The \emph{factorization monoid} of $M$, denoted by $\mathsf{Z}(M)$, is the free commutative monoid on $\mathcal{A}(M)$. The elements of $\mathsf{Z}(M)$ are called \emph{factorizations}.
\end{definition}
 If $z = a_1 \cdots a_n \in \mathsf{Z}(M)$, where $a_1, \dots, a_n \in \mathcal{A}(M)$, then $|z| := n$ is the \emph{length} of $z$. The unique monoid homomorphism $\pi \colon \mathsf{Z}(M) \to M$ satisfying that $\pi(a) = a$ for all $a \in \mathcal{A}(M)$ is the \emph{factorization homomorphism} of~$M$. For each $x \in M$,
\[
	\mathsf{Z}(x) := \pi^{-1}(x) \subseteq \mathsf{Z}(M) \quad \text{and} \quad \mathsf{L}(x) := \{|z| : z \in \mathsf{Z}(x)\}
\]
are the \emph{set of factorizations} and the \emph{set of lengths} of $x$, respectively. Factorization invariants stemming from the sets of lengths have been studied for several classes of atomic monoids and domains; see, for instance, \cite{CGGR01,CGTV16,CGP14,CHM06}. In particular, the sets of lengths of numerical monoids have been studied in~\cite{ACHP07,CDHK10,GS18}. In~\cite{GS18} the sets of lengths of numerical monoids were studied using techniques involving Puiseux monoids. An overview of sets of lengths and the role they play in factorization theory can be found in the Monthly article~\cite{aG16}.

By restricting the size of the sets of factorizations/lengths, one obtains subclasses of atomic monoids that have been systematically studied by many authors. We say that a reduced atomic monoid $M$ is
\begin{enumerate}
	\item a \emph{UFM} (or a \emph{factorial monoid}) if $|\mathsf{Z}(x)| = 1$ for all $x \in M$,
	\vspace{2pt}
	
	\item an \emph{HFM} (or a \emph{half-factorial monoid}) if $|\mathsf{L}(x)| = 1$ for all $x \in M$,
	\vspace{2pt}
	
	\item an \emph{FFM} (or a \emph{finite-factorization monoid}) if $|\mathsf{Z}(x)| < \infty$ for all $x \in M$, and
	\vspace{2pt}
	
	\item a \emph{BFM} (or a \emph{bounded-factorization monoid}) if $|\mathsf{L}(x)| < \infty$ for all $x \in M$.
\end{enumerate}
\medskip

\section{Closures and Conductor}
\label{sec:closure nad conductor}

In this section we study some algebraic aspects of Puiseux monoids.

\begin{definition}
	A \emph{Puiseux monoid} is an additive submonoid of $\qq_{\ge 0}$.
\end{definition}

Note that Puiseux monoids are natural generalizations of numerical monoids. As numerical monoids, it is clear that Puiseux monoids are reduced. However, as we shall see later, Puiseux monoids are not, in general, finitely generated or atomic.
\smallskip

\subsection{The Grothendieck Group} Recall that a monoid $M$ is torsion-free if for all $x,y \in M$ and $n \in \nn$, the equality $nx = ny$ implies that $x=y$. Each Puiseux monoid $M$ is obviously torsion-free and, therefore, $\gp(M)$ is a torsion-free group. Moreover, for a Puiseux monoid $M$, one can take the Grothendieck group $\gp(M)$ to be an additive subgroup of $\qq$, specifically,
\begin{equation} \label{eq:difference group of a PM}
	\gp(M) = \{ x-y : x,y \in M \}.
\end{equation}

Puiseux monoids can be characterized as follows.

\begin{proposition} \label{prop:characterization property of Puiseux monoids}
	For a nontrivial monoid $M$ the following statements are equivalent.
	\begin{enumerate}
		\item $M$ is a rank-$1$ torsion-free monoid that is not a group.
		
		\item $M$ is isomorphic to a Puiseux monoid.
	\end{enumerate}
\end{proposition}

\begin{proof}
	To argue (1) $\Rightarrow$ (2), first note that $\gp(M)$ is a rank-$1$ torsion-free abelian group. Therefore it follows from \cite[Section 85]{lF73} that $\gp(M)$ is isomorphic to a subgroup of $(\qq,+)$, and one can assume that $M$ is a submonoid of $(\qq,+)$. Since $M$ is not a group, \cite[Theorem 2.9]{rG84} ensures that either $M \subseteq \qq_{\le 0}$ or $M \subseteq \qq_{\ge 0}$. So $M$ is isomorphic to a Puiseux monoid. To verify (2) $\Rightarrow$ (1), let us assume that $M \subseteq \gp(M) \subseteq \qq$. As $\gp(M)$ is a subgroup of $(\qq,+)$, it is a rank-$1$ torsion-free abelian group. This implies that~$M$ is a rank-$1$ torsion-free monoid. Since $M$ is nontrivial and reduced, it cannot be a group, which completes our proof.
\end{proof}

Puiseux monoids are abundant, as the next proposition illustrates.

\begin{proposition} \label{prop:there are uncountably many PMs}
	There are uncountably many non-isomorphic Puiseux monoids.
\end{proposition}

\begin{proof}
	Consider the assignment $G \mapsto M_G := G \cap \qq_{\ge 0}$ sending each subgroup $G$ of $(\qq,+)$ to a Puiseux monoid. Clearly, $\gp(M_G) \cong G$. In addition, for all subgroups $G$ and $G'$ of $(\qq,+)$, each monoid isomorphism between $M_G$ and $M_{G'}$ naturally extends to a group isomorphism between $G$ and $G'$. Hence our assignment sends non-isomorphic groups to non-isomorphic monoids. It follows from \cite[Corollary 85.2]{lF73} that there are uncountably many non-isomorphic rank-$1$ torsion-free abelian groups. As a result, there are uncountably many non-isomorphic Puiseux monoids.
\end{proof}
\smallskip

\subsection{Closures and Conductor}
Given a monoid $M$ with Grothendieck group $\gp(M)$, the sets
\begin{itemize}
	\item $M' := \big\{ x \in \gp(M) : \text{ there exists } N \in \nn \text{ such that } nx \in M \ \text{for all} \ n \ge N \big\}$,
	
	\item $\widetilde{M} := \big\{ x \in \gp(M) : \text{ there exists } n \in \mathbb{N}  \text{ such that } nx \in M \big\}$, and
	
	\item $\widehat{M} := \big\{ x \in \gp(M) : \text{ there exists } c \in M \text{ such that } c + nx \in M \text{ for all } n \in \nn \big\}$
\end{itemize}
are called the \emph{seminormal closure}, \emph{root closure}, and \emph{complete integral closure} of $M$, respectively. It is not hard to verify that $M \subseteq M' \subseteq \widetilde{M} \subseteq \widehat{M} \subseteq \gp(M)$ for any monoid $M$. It has been recently proved by Geroldinger et al. that for Puiseux monoids the three closures coincide.

\begin{proposition} \cite[Proposition~3.1]{GGT19} \label{prop:closure of a PM}
	Let $M$ be a Puiseux monoid, and let $n = \gcd( \mathsf{n}(M^\bullet))$. Then
	\begin{equation} \label{eq:equality of closures of a PM}
		M' = \widetilde{M} = \widehat{M} = \gp(M) \cap \qq_{\ge 0} = n \langle 1/d : d \in \mathsf{d}(M^\bullet) \rangle.
	\end{equation}
\end{proposition}

A monoid $M$ is said to be \emph{root-closed} provided that $\widetilde{M} = M$. In addition, $M$ is called a \emph{Pr\"ufer monoid} if $M$ is the union of an ascending sequence of cyclic submonoids. 

\begin{corollary}
	For a Puiseux monoid $M$, the following statements are equivalent.
	\begin{enumerate}
		\item $M$ is root-closed.
		\vspace{2pt}
		
		\item $M = n \langle 1/d : d \in \mathsf{d}(M^\bullet) \rangle$, where $n = \gcd \mathsf{n}(M^\bullet)$.
		\vspace{2pt}
		
		\item $\emph{gp}(M) = M \cup -M$.
		\vspace{2pt}
		
		\item $M$ is a Pr\"ufer monoid.
	\end{enumerate}
\end{corollary}

\begin{proof}
	The equivalences (1) $\Leftrightarrow$ (2) $\Leftrightarrow$ (3) follow from Proposition~\ref{prop:closure of a PM}, while the equivalence (1) $\Leftrightarrow$ (4) follows from~\cite[Theorem~13.5]{rG84}.
\end{proof}

We now characterize finitely generated Puiseux monoids in terms of its root closures.

\begin{proposition} \label{prop:fg PM are NM}
	For a Puiseux monoid $M$ the following statements are equivalent.
	\begin{enumerate}
		\item $\widetilde{M} \cong (\nn_0,+)$.
		\vspace{2pt}
		
		\item $M$ is finitely generated.
		\vspace{2pt}
		
		\item $\mathsf{d}(M^\bullet)$ is finite.
		\vspace{2pt}
		
		\item $M$ is isomorphic to a numerical monoid.
	\end{enumerate}
\end{proposition}

\begin{proof}
	To prove (1) $\Rightarrow (2)$, suppose that $\widetilde{M} \cong (\nn_0,+)$. Proposition~\ref{prop:closure of a PM} ensures that $\mathsf{d}(M^\bullet)$ is finite. Now if $\ell := \text{lcm} \, \mathsf{d}(M^\bullet)$, then $\ell M$ is submonoid of $(\nn_0,+)$ that is isomorphic to $M$. Hence $M$ is finitely generated. To argue (2) $\Rightarrow$ (3), it suffices to notice that if $S$ is a finite generating set of $M$, then every element of $\mathsf{d}(M^\bullet)$ divides $\text{lcm} \, \mathsf{d}(S^\bullet)$. For (3) $\Rightarrow$ (4), let $\ell := \text{lcm} \, \mathsf{d}(M^\bullet)$. Then note that $\ell M$ is a submonoid of $(\nn_0,+)$ that is isomorphic to $M$. As a result, $M$ is isomorphic to a numerical monoid. To prove (4) $\Rightarrow$ (1), assume that $M$ is a numerical monoid and that $\gp(M)$ is a subgroup of $(\zz,+)$. By definition of $\widetilde{M}$, it follows that $\widetilde{M} \subseteq \nn_0$. On the other hand, the fact that $\nn_0 \setminus M$ is finite immediately implies that $\nn_0 \subseteq \widetilde{M}$. Thus, $\widetilde{M} = (\nn_0,+)$. 
\end{proof}

\begin{corollary} \label{cor:closure of a non-finitely generated PM is antimatter}
	A Puiseux monoid $M$ is not finitely generated if and only if $\widetilde{M}$ is antimatter.
\end{corollary}

\begin{proof}
	Suppose first that $M$ is not finitely generated. Set $n = \gcd(\mathsf{n}(M^\bullet))$. It follows from Proposition~\ref{prop:closure of a PM} that $\widetilde{M} = \langle n/d : d \in \mathsf{d}(M^\bullet) \rangle$. Fix $d \in \mathsf{d}(M^\bullet)$. Since $\mathsf{d}(M^\bullet)$ is an infinite set that is closed under taking least common multiples, there exists $d' \in \mathsf{d}(M^\bullet)$ such that $d'$ properly divides $d$. As a consequence, $n/d'$ properly divides $n/d$ in $\widetilde{M}$ and so $n/d \notin \mathcal{A}(\widetilde{M})$. As none of the elements in the generating set $\{n/d : d \in \mathsf{d}(M^\bullet)\}$ of $\widetilde{M}$ is an atom, $\widetilde{M}$ must be antimatter. The reverse implication is an immediate consequence of Proposition~\ref{prop:fg PM are NM}.
\end{proof}
\smallskip

The \emph{conductor} of a monoid $M$, denoted by $(M : \widehat{M})$, is defined to be
\begin{equation} \label{eq:conductor}
	(M : \widehat{M}) = \{ x \in \gp(M) : x + \widehat{M} \subseteq M \}.
\end{equation}
 For a numerical monoid $N$, the term `conductor' also refers to the number $\mathfrak{f}(N) + 1$, where $\mathfrak{f}(N)$ denotes the Frobenius number of~$N$. This does not generate ambiguity as the following example illustrates.

\begin{example}
	Let $N$ be a numerical monoid, and let $\mathfrak{f}(N)$ be the Frobenius number of $N$. It follows from~(\ref{eq:difference group of a PM}) that $\gp(N) = \zz$. Therefore Proposition~\ref{prop:closure of a PM} guarantees that $\widehat{N} = \nn_0$. For each $n \in N$ with $n \ge \mathfrak{f}(N) + 1$, it is clear that $n + \widehat{N} = n + \nn_0 \subseteq N$. On the other hand, for each $n \in \zz$ with $n \le \mathfrak{f}(N)$ the fact that $\mathfrak{f}(N) \in n + \widehat{N}$ implies that $n + \widehat{N} \nsubseteq N$. As a result,
	\begin{equation} \label{eq:Frobenius number and conductor}
		 (N : \widehat{N}) = \{ n \in \zz : n \ge \mathfrak{f}(N) + 1 \}.
	\end{equation}
	As the equality of sets~(\ref{eq:Frobenius number and conductor}) shows, the minimum of $(N : \widehat{N})$ is $\mathfrak{f}(N) + 1$, namely, the conductor number of $N$, as defined in the context of numerical monoids.
\end{example}

The conductor of a Puiseux monoid was first considered in~\cite{GGT19}, where the following result was established.

\begin{proposition} \label{prop:conductor of a PM} 
	Let $M$ be a Puiseux monoid. Then the following statements hold.
	\begin{enumerate}
		\item If $M$ is root-closed, then $(M : \widehat{M}) = \widehat{M} = M$.
		
		\item If $M$ is not root-closed and $\sigma = \sup \, \widehat{M} \setminus M$.
		\begin{enumerate} \label{part 2: conductor of PM}
			\item If $\sigma = \infty$, then $(M : \widehat{M}) = \emptyset$.
			\item If $\sigma < \infty$, then $(M : \widehat{M}) = M_{\ge \sigma}$.
		\end{enumerate}
	\end{enumerate}
\end{proposition}

\begin{remark}
	With notation as in Proposition~\ref{prop:conductor of a PM}.\ref{part 2: conductor of PM}, although $\widehat{M}_{> \sigma} = M_{> \sigma}$ holds, it can happen that $\widehat{M}_{\ge \sigma} \neq M_{\ge \sigma}$. For instance, consider the Puiseux monoid $\{0\} \cup \qq_{>1}$.
\end{remark}
\medskip

\subsection{Homomorphisms Between Puiseux Monoids} As we are about to show, homomorphisms between Puiseux monoids are given by rational multiplication.

\begin{proposition} \label{prop:homomorphisms between PM}
	The homomorphisms between Puiseux monoids are given by rational multiplication.
\end{proposition}

\begin{proof}
	Every rational-multiplication map is clearly a homomorphism. Suppose, on the other hand, that $\varphi \colon M \to M'$ is a homomorphism between Puiseux monoids. As the trivial homomorphism is multiplication by $0$, one can assume without loss of generality that $\varphi$ is nontrivial. Let $\{n_1, \dots, n_k\}$ be the minimal generating set of the additive monoid $N := M \cap \nn_0$. Since $\varphi$ is nontrivial, $k \ge 1$ and $\varphi(n_j) \neq 0$ for some $j \in \{ 1, \dots, k \}$. Set $q = \varphi(n_j)/n_j$ and then take $r \in M^\bullet$ and $c_1, \dots, c_k \in \nn_0$ satisfying that $\mathsf{n}(r) = c_1 n_1 + \dots + c_k n_k$. As $n_i \varphi(n_j) = \varphi(n_i n_j) = n_j \varphi(n_i)$ for every $i \in \{ 1, \dots, k \}$,
	\[
	\varphi(r) = \frac 1{\mathsf{d}(r)} \varphi(\mathsf{n}(r)) = \frac 1{\mathsf{d}(r)} \sum_{i=1}^k c_i \varphi(n_i) = \frac 1{\mathsf{d}(r)} \sum_{i=1}^k c_i n_i \frac{\varphi (n_j)}{n_j} = rq.
	\]
	Thus, the homomorphism $\varphi$ is multiplication by $q \in \mathbb{Q}_{>0}$.
\end{proof}

%% file: tex/ch2.tex
\chapter{Atomicity of Puiseux monoids} \label{chap:atomicity}


It is well known that in the class consisting of all monoids, the following chain of implications holds.
\begin{equation} \label{eq:monoid atomicity taxonomy}
\textbf{UFM} \ \Rightarrow \ \textbf{HFM} \stackrel{\text{reduced}}{\Rightarrow} \ \textbf{FFM} \ \Rightarrow \ \textbf{BFM} \ \Rightarrow \ \textbf{ACCP} \ \Rightarrow \ \textbf{atomic monoid}
\end{equation}
It is also known that, in general, none of the implications in~(\ref{eq:monoid atomicity taxonomy}) is reversible (even in the class of integral domains~\cite{AAZ90}). In this chapter, we provide various examples to illustrate that none of the above implications, except the first one, is reversible in the class of Puiseux monoids. We characterize the Puiseux monoids belonging to the first two classes of the chain of implications~(\ref{eq:monoid atomicity taxonomy}). For each of the last four classes, we find a family of Puiseux monoids belonging to such a class but not to the class right before.
\medskip

\section{A Class of Atomic Puiseux Monoids} We begin this section with an easy characterization of finitely generated Puiseux monoids in terms of their atomicity.

\begin{proposition}
	A Puiseux monoid $M$ is finitely generated if and only if $M$ is atomic and $\mathcal{A}(M)$ is finite.
\end{proposition}

\begin{proof}
	The direct implication follows immediately from the fact that finitely generated Puiseux monoids are isomorphic to numerical monoids via Proposition~\ref{prop:fg PM are NM}. The reverse implication is also obvious because the atomicity of $M$ means that $M$ is generated by $\mathcal{A}(M)$, which is finite.
\end{proof}

Corollary~\ref{cor:closure of a non-finitely generated PM is antimatter} yields, however, instances of non-finitely generated Puiseux monoids containing no atoms. As the next example shows, for every $n \in \nn$ there exists a non-finitely generated Puiseux monoid containing exactly $n$ atoms.

\begin{example}
	Let $m \in \nn$, and take distinct prime numbers $p$ and $q$ with $q > m$. Consider the Puiseux monoid $M = \big\langle \{ m, \dots, 2m-1 \} \cup \{qp^{-m-i} : i \in \nn\} \big\rangle$. To verify that $\mathcal{A}(M) = \{ m, \dots, 2m-1 \}$, write $a \in \{ m, \dots, 2m-1 \}$ as
	\begin{equation} \label{eq:monoid with m atoms}
		a = a' + \sum_{n=1}^N c_n \frac{q}{p^{m+n}},
	\end{equation}
	where $a' \in \{0\} \cup \{ m, \dots, 2m-1 \}$ and $c_n \in \nn_0$ for every $n \in \{ 1, \dots, N \}$. After cleaning denominators in~(\ref{eq:monoid with m atoms}), one finds that $q \mid a - a'$. So $a = a'$ and $c_1 =  \dots = c_N = 0$, which implies that $a \in \mathcal{A}(M)$. Thus, $\{ m, \dots, 2m-1 \} \subseteq \mathcal{A}(M)$. Clearly, $qp^{-m-i} \notin \mathcal{A}(M)$ for any $i \in \nn$. Hence $\mathcal{A}(M) = \{ m, \dots, 2m-1 \}$, and so $|\mathcal{A}(M)| = m$. As $\mathsf{d}(M^\bullet)$ is not finite, it follows from Proposition~\ref{prop:fg PM are NM} that $M$ is not finitely generated.
\end{example}

Perhaps the class of non-finitely generated Puiseux monoids that has been most thoroughly studied is the class of multiplicatively cyclic Puiseux monoids~\cite{CGG19}.

\begin{definition}
	For $r \in \qq_{>0}$, we call $M_r = \langle r^n : n \in \nn_0 \rangle$ a \emph{multiplicatively cyclic} Puiseux monoid. 
\end{definition}

Although $M_r$ is, indeed, a rational cyclic semiring, we shall only be concerned here with its additive structure. The atomicity of multiplicatively cyclic Puiseux monoids was first studied in \cite[Section~6]{GG17} while several factorization aspects were recently investigated in~\cite{CGG19}. 

\begin{proposition} \label{prop:atomic classification of multiplicative cyclic Puiseux monoids}
	For $r \in \qq_{> 0}$, consider the multiplicatively cyclic Puiseux monoid~$M_r$. Then the following statements hold.
	\begin{enumerate}
		\item If $r \ge 1$, then $M_r$ is atomic and \label{item:case r at least 1}
		\begin{itemize}
			\item either $r \in \nn$ and so $M_r = \nn_0$,
			
			\item or $r \notin \nn$ and so $\mathcal{A}(M_r) = \{r^n : n \in \nn_0\}$.
		\end{itemize}
		\item If $r < 1$, then \label{item:case r less than 1}
		\begin{itemize}
			\item either $\mathsf{n}(r) = 1$ and so $M_r$ is antimatter,
			
			\item or $\mathsf{n}(r) \neq 1$ and $M_r$ is atomic with $\mathcal{A}(M_r) = \{r^n : n \in \nn_0\}$.
		\end{itemize}
	\end{enumerate}
\end{proposition}

\begin{proof}
	To argue~(\ref{item:case r at least 1}), suppose that $r \ge 1$. If $r \in \nn$, then it easily follows that $M_r = \nn_0$. Then we assume that $r \notin \nn$. Clearly, $\mathcal{A}(M_r) \subseteq \{r^n : n \in \nn_0\}$. To check the reverse inequality, fix $j \in \nn_0$ and write $r^j = \sum_{i=0}^N \alpha_i r^i$ for some $N \in \nn_0$ and $\alpha_i \in \nn_0$ for every $i \in \{ 0, \dots, N \}$. As $(r^n)_{n \in \mathbb{N}_0}$ is an increasing sequence, one can assume that $N \le j$. Then, after cleaning denominators in $r^j = \sum_{i=0}^N \alpha_i r^i$ we obtain $N=j$ as well as $\alpha_j = 1$ and $\alpha_i = 0$ for every $i \neq j$. Hence $r^j \in \mathcal{A}(M_r)$ for every $j \in \nn_0$, which yields the second statement of~(\ref{item:case r at least 1}).
	
	Now suppose that $r < 1$. If $\mathsf{n}(r) = 1$ then $r^n = \mathsf{d}(r)r^{n+1}$ for every $n \in \nn_0$, and so $M_r$ is antimatter, which is the first statement of~(\ref{item:case r less than 1}). Finally, suppose that $\mathsf{n}(r) > 1$. Fix $j \in \nn$, and notice that $r^i \nmid_{M_r} r^j$ for any $i < j$. Then write  $r^j = \sum_{i=j}^{j+k} \beta_i r^i$, for some $k \in \nn_0$ and $\beta_i \in \nn_0$ for every $i \in \{ j, \dots, j+k \}$. Notice that $\beta_j \in \{0,1\}$. Suppose by a contradiction that $\beta_j = 0$. In this case, $k \ge 1$. Let $p$ be a prime dividing $\mathsf{n}(r)$, and let~$\alpha$ be the maximum power of $p$ dividing $\mathsf{n}(r)$. From $r^j = \sum_{i=j}^{j+k} \beta_i r^i$ one obtains
	\begin{equation} \label{eq:multiplicative cyclic 3}
		\alpha j = \pval\big( r^j \big) = \pval \bigg( \sum_{i=1}^k \beta_{j+i} r^{j+i} \bigg) \ge \min_{i \in \{ 1, \dots, k \} } \big\{ \pval \big(\beta_{j+i} r^{j+i} \big) \big\} \ge \alpha(j+m),
	\end{equation}
	where $m = \min\{i \in \{ 1, \dots, k \} : \beta_{j+i} \neq 0\}$. The inequality \eqref{eq:multiplicative cyclic 3} yields the desired contradiction. Hence $r^j \in \mathcal{A}(M_r)$ for every $j \in \nn_0$, which implies the second statement of~(\ref{item:case r less than 1}).
\end{proof}

\begin{corollary}\label{cor:atomic and not ACCP}
	For each $r \in \qq \cap (0,1)$ with $\mathsf{n}(r) \neq 1$, the monoid $M_r$ is an atomic monoid that does not satisfy the ACCP.
\end{corollary}

\begin{proof}
	Proposition~\ref{prop:atomic classification of multiplicative cyclic Puiseux monoids} guarantees the atomicity of $M_r$. To verify that $M_r$ does not satisfy the ACCP, consider the sequence of principal ideals $(\mathsf{n}(r)r^n + M_r)_{n \in \nn_0}$. Since
	\[
		\mathsf{n}(r)r^n = \mathsf{d}(r)r^{n+1} = (\mathsf{d}(r) - \mathsf{n}(r))r^{n+1} + \mathsf{n}(r)r^{n+1},
	\]
	$\mathsf{n}(r)r^{n+1} \mid_{M_r} \mathsf{n}(r)r^n$ for every $n \in \nn_0$. Therefore $(\mathsf{n}(r)r^n + M_r)_{n \in \nn_0}$ is an ascending chain of principal ideals. In addition, it is clear that such a chain of ideals does not stabilize. Hence $M_r$ does not satisfy the ACCP, which completes the proof.
\end{proof}
\medskip

\section{A Class of ACCP Puiseux Monoids} \label{sec:ACCP} We proceed to present a class of ACCP Puiseux monoids containing a subclass of monoids that are not BFMs.

\begin{definition}
	A Puiseux monoid $M$ is called \emph{prime reciprocal} provided that there exists a generating set $S$ of $M$ satisfying the following two conditions.
	\begin{itemize}
		\item $\mathsf{d}(S) \subseteq \pp$.
		
		\item $\mathsf{d}(a) = \mathsf{d}(a')$ implies that $a = a'$ for all $a,a' \in S$.
	\end{itemize}
\end{definition}
\smallskip

Let us prove that all prime reciprocal Puiseux monoids satisfy the ACCP.

\begin{theorem} \label{thm:a class of ACCP monoids}
	Every prime reciprocal Puiseux monoid satisfies the ACCP.
\end{theorem}

\begin{proof}
	Let $(p_n)_{n \in \nn}$ be a strictly increasing sequence of prime numbers, and let $(a_n)_{n \in \nn}$ be a sequence of positive integers such that $p_n \nmid a_n$. Now set $M = \langle a_n/p_n : n \in \nn \rangle$. It is not hard to verify that for all $x \in M$ there exist $k,n \in \nn_0$ and $\alpha_i \in \{ 0, \dots, p_i \}$ such that
	\begin{align} \label{eq:decomposition existence}
		x = n + \sum_{i=1}^k \alpha_i \frac{a_i}{p_i}.
	\end{align}
	Let us now check that the sum decomposition in (\ref{eq:decomposition existence}) is unique. To do this, take $k', m \in~\nn_0$ and $\beta_i \in \{ 0, \dots, p_i \}$ such that
	\begin{align} \label{eq:decomposition uniqueness}
		n + \sum_{i=1}^k \alpha_i \frac{a_i}{p_i} = m + \sum_{i=1}^{k'} \beta_i \frac{a_i}{p_i}.
	\end{align}
	Suppose, without loss of generality, that $k = k'$. After isolating $(\alpha_i - \beta_i)a_i/p_i$ in (\ref{eq:decomposition uniqueness}) and applying the $p_i$-adic valuation, we obtain that $p_i \mid \alpha_i - \beta_i$, which implies that $\alpha_i = \beta_i$ for each $i \in \{ 1, \dots, k \}$. As a consequence, $n = m$ and we can conclude that the uniqueness of the decomposition in~(\ref{eq:decomposition existence}) holds.
	
	Now suppose, by way of contradiction, that $M$ does not satisfy the ACCP. Then there exists a strictly decreasing sequence $(q_n)_{n \in \mathbb{N}}$ of elements in $M$ such that
	\begin{align} \label{eq:chain of ideals}
		q_1 + M \subsetneq q_2 + M \subsetneq \cdots.
	\end{align}
	In the unique sum decomposition of $q_1$ as in (\ref{eq:decomposition existence}), let $N$ be the integer and let $p_{n_1}, \dots, p_{n_k}$ be the distinct prime denominators of the atoms with nonzero coefficients. By~(\ref{eq:chain of ideals}), for every $n \in \nn$ there exists $r_{n+1} \in M$ such that $q_n = q_{n+1} + r_{n+1}$ and, therefore,
	\[
		\sum_{i=1}^n r_{i+1} \le q_{n+1} + \sum_{i=1}^n r_{i+1} = q_1
	\]
	for every $n \in \nn$. Thus, $\lim_{n \to \infty} r_n = 0$. Then, for any finite subset $A$ of $\mathcal{A}(M)$ there exist a large enough $\ell \in \nn$ and $a \in \mathcal{A}(M)$ such that $a \mid_M r_\ell$ and $a \notin A$. Hence we can find $t \in \nn$ and (possibly repeated) $p_{m_1}, \dots, p_{m_t} \in \pp$ satisfying that $|\{p_{m_1}, \dots, p_{m_t}\}| > N + k$ and $a_{m_i}/p_{m_i} \in \mathcal{A}(r_{i+1})$ for each $i \in \{ 1, \dots, t \}$. Take $z_i \in \mathsf{Z}(r_i)$ containing the atom $a_{m_i}/p_{m_i}$, and take $z_{t+1} \in \mathsf{Z}(q_{t+1})$. As $q_1 = q_{t+1} + \sum_{i=1}^t r_{i+1}$, we have $z = z_{t+1} + \sum_{i=1}^t z_i \in \mathsf{Z}(q_1)$. By the uniqueness of the sum decomposition in (\ref{eq:decomposition existence}), the factorization $z$ contains at least $\mathsf{d}(a)$ copies of each atom $a$ for which $\mathsf{d}(a) \in P := \{p_{m_1}, \dots, p_{m_t}\} \setminus \{p_{n_1}, \dots, p_{n_k}\}$. Since $|\{p_{m_1}, \dots, p_{m_t}\}| > N + k$, it follows that $|P| > N$. Thus,
	\[
		N \ge \sum_{a \in \mathcal{A}(M) \, \mid \, \mathsf{d}(a) \in P} \! \! \! \mathsf{d}(a) a \ge |P| > N,
	\]
	which is a contradiction. Hence $M$ satisfies the ACCP.
\end{proof}

\begin{corollary} \label{cor:prime reciprocal PM are atomic}
	Prime reciprocal Puiseux monoids are atomic.
\end{corollary}

\begin{proof}
	It follows from Theorem~\ref{thm:a class of ACCP monoids}, along with the chain of implications we exhibited at the beginning of this chapter.
\end{proof}

\begin{corollary} \label{cor:ACCP that is not BFM}
	There are Puiseux monoids satisfying the ACCP that are not BFMs.
\end{corollary}

\begin{proof}
	Consider the Puiseux monoid $M = \langle 1/p : p \in \pp \rangle$. We have seen in Theorem~\ref{thm:a class of ACCP monoids} that $M$ satisfies the ACCP. However, it is clear that $p \in \mathsf{L}(1)$ for every $p \in \pp$. As $|\mathsf{L}(1)| = \infty$, the monoid $M$ is not a BFM.
\end{proof}
\smallskip

\section{BFMs and FFMs} Our next main goal is to find a large class of Puiseux monoids that are BFMs. This amounts to prove the following result.

\begin{theorem} \label{thm:BF sufficient condition}
	Let $M$ be a Puiseux monoid. If $0$ is not a limit point of $M^\bullet$, then $M$ is a BFM.
\end{theorem}

\begin{proof}
	It is clear that $\mathcal{A}(M)$ consists of those elements of $M^\bullet$ that cannot be written as the sum of two positive elements of $M$. Since $0$ is not a limit point of $M$ there exists $\epsilon > 0$ such that $\epsilon < x$ for all $x \in M^\bullet$. Now we show that $M = \langle \mathcal{A}(M) \rangle$. Take $x \in M^\bullet$. Since $\epsilon$ is a lower bound for $M^\bullet$, the element $x$ can be written as the sum of at most $\lfloor x/\epsilon \rfloor$ elements of $M^\bullet$. Take the maximum natural $m$ such that $x = a_1 + \dots + a_m$ for some $a_1, \dots, a_m \in M^\bullet$. By the maximality of $m$, it follows that $a_i \in \mathcal{A}(M)$ for every $i \in \{ 1, \dots, m \}$, which means that $x \in \langle \mathcal{A}(M) \rangle$. Hence $M$ is atomic. We have already noticed that every element $x$ in $M^\bullet$ can be written as the sum of at most $\lfloor x/\epsilon \rfloor$ atoms, i.e., $|\mathsf{L}(x)| \le \lfloor x/\epsilon \rfloor$ for all $x \in M$. Thus, $M$ is a BFM.
\end{proof}

The converse of Theorem~\ref{thm:BF sufficient condition} does not hold. The following example sheds some light upon this observation.

\begin{example} \label{ex:BF PM with 0 as a limit point}
	Let $(p_n)_{n \in \nn}$ and $(q_n)_{n \in \nn}$ be two strictly increasing sequence of prime numbers such that $q_n > p_n^2$ for every $n \in \nn$. Set $M := \big \langle \frac{p_n}{q_n} : n \in \nn \big \rangle$. It follows from~\cite[Corollary~5.6]{GG17} that $M$ is atomic, and it is easy to verify that $\mathcal{A}(M) = \{p_n/q_n : n \in \nn\}$. To argue that $M$ is indeed a BFM, take $x \in M^\bullet$ and note that since both sequences $(p_n)_{n \in \nn}$ and $(q_n)_{n \in \nn}$ are strictly increasing, there exists $N \in \nn$ such that $q_n \nmid \mathsf{d}(x)$ and $p_n > x$ for every $n \ge N$. As a result, if $z \in \mathsf{Z}(x)$, then none of the atoms in $\{p_n/q_n : n > N\}$ can appear in $z$. From this, one can deduce that $\mathsf{Z}(x)$ is finite. Then $\mathsf{L}(x)$ is finite for any $x \in M$, and so $M$ is a BFM. However, $q_n > p_n^2$ for every $n \in \nn$ implies that $0$ is a limit point of $M^\bullet$.
\end{example}

As we have seen in Corollary~\ref{cor:ACCP that is not BFM}, not every ACCP Puiseux monoid is a BFM. However, under a mild assumption on conductors, each of these atomic conditions is equivalent to having~$0$ as a limit point. The next theorem, along with the next two examples, has been recently established in~\cite{GGT19}.

\begin{theorem} \cite[Theorem~3.4]{GGT19} \label{thm:strongly primary characterization when conductor non-empty}
	If $M$ is a nontrivial Puiseux monoid with nonempty conductor, then the following conditions are equivalent.
	\begin{enumerate}
		\item $0$ is not a limit point of $M^\bullet$.
		
		\item $M$ is a BFM.
		
		\item $M$ satisfies the ACCP.
	\end{enumerate}
\end{theorem}

As Corollary~\ref{cor:ACCP that is not BFM} and Example~\ref{ex:BF PM with 0 as a limit point} indicate, without the nonempty-conductor condition, none of the last two statements in Theorem~\ref{thm:strongly primary characterization when conductor non-empty} implies its predecessor. In addition, even inside the class of Puiseux monoids with nonempty conductor neither being atomic nor being an FFM is equivalent to being a BFM (or satisfying the ACCP).

\begin{example}
	Consider the Puiseux monoid $M := \{0\} \cup \qq_{\ge 1}$. It is clear that the conductor of $M$ is nonempty. In addition, it follows from Theorem~\ref{thm:strongly primary characterization when conductor non-empty} that $M$ is a BFM. Note that $\mathcal{A}(M) = [1,2)$. However, $M$ is far from being an FFM; for instance, the formal sum $(1 + 1/n) + (x - 1 - 1/n)$ is a length-$2$ factorization in $\mathsf{Z}(x)$ for all $x \in (2,3] \cap \qq$ and $n \ge \big\lceil \frac{1}{x-2} \big\rceil$, which implies that $|\mathsf{Z}(x)| = \infty$ for all $x \in M_{>2}$.
\end{example}

\begin{example}
	Now consider the Puiseux monoid $M =  \langle 1/p : p \in \pp \rangle \cup \qq_{\ge 1}$. Since the monoid $\langle 1/p : p \in \pp \rangle$ is atomic by Theorem~\ref{thm:a class of ACCP monoids}, it is not hard to check that $M$ is also atomic. It follows from Proposition~\ref{prop:conductor of a PM} that $M$ has nonempty conductor. Since $0$ is a limit point of~$M^\bullet$, Theorem~\ref{thm:strongly primary characterization when conductor non-empty} ensures that $M$ does not satisfy the ACCP.
\end{example}
\smallskip

\medskip
\section{Monotone Puiseux Monoids}
\label{sec:increasing PM}

\subsection{Increasing Puiseux Monoids} We are in a position now to begin our study of the atomic structure of Puiseux monoids generated by monotone sequences.

\begin{definition}
	A Puiseux monoid $M$ is said to be \emph{increasing} (resp., \emph{decreasing}) if it can be generated by an increasing (resp., decreasing) sequence. A Puiseux monoid is \emph{monotone} if it is either increasing or decreasing.
\end{definition}

Not every Puiseux monoid is monotone, as the next example illustrates.

\begin{example} \label{ex:bounded PM that is neither decreasing nor increasing}
	Let $p_1, p_2, \dots$ be an increasing enumeration of the set of prime numbers. Consider the Puiseux monoid $M = \langle A \cup B \rangle$, where
	\[
		A = \bigg\{ \frac{1}{p_{2n}} : n \in \nn \bigg\} \ \text{ and } \ B = \bigg\{ \frac{p_{2n-1} - 1}{p_{2n-1}} : n \in \nn \bigg\}.
	\]
	It follows immediately that both $A$ and $B$ belong to $\mathcal{A}(M)$. So $M$ is atomic, and $\mathcal{A}(M) =~A \cup B$. Every generating set of $M$ must contain $A \cup B$ and so will have at least two limit points, namely, $0$ and $1$. Since every monotone sequence of rationals can have at most one limit point in the real line, we conclude that $M$ is not monotone.
\end{example}

The following proposition describes the atomic structure of the family of increasing Puiseux monoids.

\begin{proposition} \label{prop:atoms of increasing monoids}
	Every increasing Puiseux monoid is atomic. Moreover, if $(r_n)_{n \in \mathbb{N}}$ is an increasing sequence of positive rationals generating a Puiseux monoid $M$, then $\mathcal{A}(M) = \{r_n : r_n \notin \langle r_1, \dots, r_{n-1} \rangle\}$.
\end{proposition}

\begin{proof}
	The fact that $M$ is atomic follows from observing that $r_1$ is a lower bound for~$M^\bullet$ and so $0$ is not a limit point of $M$. To prove the second statement, set
	\[
		A = \{r_n : r_n \notin \langle r_1, \dots, r_{n-1} \rangle\},
	\]
	and rename the elements of $A$ in a strictly increasing sequence (possibly finite), namely, $(a_n)_{n \in \mathbb{N}}$. Note that $M = \langle A \rangle$ and $a_n \notin \langle a_1, \dots, a_{n-1} \rangle$ for any $n \in \nn$. Since $a_1$ is the smallest nonzero element of $M$, we obtain that $a_1 \in \mathcal{A}(M)$. Suppose now that $n$ is a natural such that $2 \le n \le |A|$. Because $(a_n)_{n \in \mathbb{N}}$ is a strictly increasing sequence and $a_n \notin \langle a_1,\dots, a_{n-1} \rangle$, one finds that $a_n$ cannot be written as a sum of elements in $M$ in a non-trivial manner. Hence $a_n$ is an atom for every $n \in \nn$ and, therefore, we can conclude that $\mathcal{A}(M) = A$.
\end{proof}

Now we use Proposition~\ref{prop:atoms of increasing monoids} to show that every Puiseux monoid that is not isomorphic to a numerical monoid has an atomic submonoid with infinitely many atoms. Let us first prove the next lemma.

\begin{lemma} \label{lem:denominator set of like-NS PM}
	Let $M$ be a nontrivial Puiseux monoid. Then $\mathsf{d}(M^\bullet)$ is finite if and only if $M$ is isomorphic to a numerical monoid.
\end{lemma}

\begin{proof}
	Suppose first that $\mathsf{d}(M^\bullet)$ is finite. Since $M$ is not trivial, $M^\bullet$ is not empty. Take $a \in \nn$ to be the least common multiple of $\mathsf{d}(M^\bullet)$. Since $aM$ is a submonoid of $\nn_0$, it is isomorphic to a numerical monoid. Furthermore, the map $\varphi \colon M \to aM$ defined by $\varphi(x) = ax$ is a monoid isomorphism. Thus, $M$ is isomorphic to a numerical monoid. Conversely, suppose that $M$ is isomorphic to a numerical monoid. Because every numerical monoid is finitely generated, so is $M$. Hence $\mathsf{d}(M^\bullet)$ is finite.
\end{proof}

\begin{proposition}
	If $M$ is a nontrivial Puiseux monoid, then it satisfies exactly one of the following conditions:
	\begin{enumerate}
		\item $M$ is isomorphic to a numerical monoid;
		\vspace{3pt}
		\item $M$ contains an atomic submonoid with infinitely many atoms.
	\end{enumerate}
\end{proposition}

\begin{proof}
	Suppose that $M$ is not isomorphic to any numerical monoid. Take $r_1 \in M^\bullet$\!\!. By Lemma~\ref{lem:denominator set of like-NS PM}, the set $\mathsf{d}(M^\bullet)$ is not finite. Therefore $\mathsf{d}(\langle r_1 \rangle^\bullet)$ is strictly contained in $\mathsf{d}(M^\bullet)$. Take $r'_2 \in M^\bullet$ such that $\mathsf{d}(r'_2) \notin \mathsf{d}(\langle r_1 \rangle^\bullet)$. Let $r_2$ be the sum of $m_2$ copies of~$r'_2$, where $m_2$ is a natural number so that $\gcd(m_2, \mathsf{d}(r'_2)) = 1$ and $m_2 r'_2 > r_1$. Setting $r_2 = m_2r'_2$, we notice that $r_2 > r_1$ and $r_2 \notin \langle r_1 \rangle$. Now suppose that $r_1, \dots, r_n \in M$ have been already chosen so that $r_{i+1} > r_i$ and $r_{i+1} \notin \langle r_1, \dots, r_i \rangle$ for $i = 1,\dots,n-1$. Once again, by using Lemma~\ref{lem:denominator set of like-NS PM} we can guarantee that $\mathsf{d}(M^\bullet)\setminus \mathsf{d}(\langle r_1, \dots, r_n \rangle^\bullet)$ is not empty. Take $r'_{n+1} \in M^\bullet$ such that $\mathsf{d}(r'_{n+1}) \notin \mathsf{d}(\langle r_1, \dots, r_n \rangle^\bullet)$, and choose $m_{n+1} \in \nn$ so that $\gcd(m_{n+1}, \mathsf{d}(r'_{n+1})) = 1$ and $m_{n+1} r'_{n+1} > r_n$. Taking $r_{n+1} = m_{n+1}r'_{n+1}$, one finds that $r_{n+1} > r_n$ and $r_{n+1} \notin \langle r_1, \dots, r_n \rangle$. Using the method just described, we obtain an infinite sequence $(r_n)_{n \in \mathbb{N}}$ of elements in $M$ satisfying that $r_{n+1} > r_n$ and $r_{n+1} \notin \langle r_1, \dots, r_n \rangle$ for every $n \in \nn$. By Proposition~\ref{prop:atoms of increasing monoids}, the submonoid $N = \langle r_n : n \in \nn \rangle$ is atomic and $\mathcal{A}(N) = \{r_n : n \in \nn\}$. Hence $M$ has an atomic submonoid with infinitely many atoms, namely, $N$.
	
	Finally, note that conditions (1) and (2) exclude each other; this is because a submonoid of a numerical monoid is either trivial or isomorphic to a numerical monoid and so it must contain only finitely many atoms.
\end{proof}

Now we split the family of increasing Puiseux monoids into two fundamental subfamilies. We will see that these two subfamilies have different behavior. We say that a sequence of rationals is \emph{strongly increasing} if it increases to infinity. On the other hand, a bounded increasing sequence of rationals is called \emph{weakly increasing}.

\begin{definition}
	A Puiseux monoid is said to be \emph{strongly} (resp., \emph{weakly}) \emph{increasing} if it can be generated by a strongly (resp., weakly) increasing sequence.
\end{definition}

\begin{proposition}
	Every increasing Puiseux monoid is either strongly increasing or weakly increasing. A Puiseux monoid is both strongly and weakly increasing if and only if it is isomorphic to a numerical monoid.
\end{proposition}

\begin{proof}
	The first statement follows straightforwardly. For the second statement, suppose that $M$ is a Puiseux monoid that is both strongly and weakly increasing. By Proposition~\ref{prop:atoms of increasing monoids}, the monoid $M$ is atomic, and its set of atoms can be listed increasingly. Let $(a_n)_{n \in \mathbb{N}}$ be an increasing sequence with underlying set $\mathcal{A}(M)$. Suppose, by way of contradiction, that $\mathcal{A}(M)$ is not finite. Since $M$ is strongly increasing, $(a_n)_{n \in \mathbb{N}}$ must be unbounded. However, the fact that $M$ is weakly decreasing forces $(a_n)_{n \in \mathbb{N}}$ to be bounded, which is a contradiction. Hence $\mathcal{A}(M)$ is finite, which implies that $M$ is isomorphic to a numerical monoid.
	
	To prove the converse implication, take $M$ to be a Puiseux monoid isomorphic to a numerical monoid. So $M$ is finitely generated, namely, $M = \langle r_1, \dots, r_n \rangle$ for some $n \in \nn$ and $r_1 < \dots < r_n$. The sequence $(a_n)_{n \in \mathbb{N}}$ defined by $a_k = r_k$ if $k \le n$ and $a_k = kr_n$ if $k > n$ is an unbounded increasing sequence generating $M$. Similarly, the sequence $(b_n)_{n \in \mathbb{N}}$ defined by $b_k = r_k$ if $k \le n$ and $b_k = r_n$ if $k > n$ is a bounded increasing sequence generating $M$. Consequently, $M$ is both strongly and weakly increasing.
\end{proof}

We will show that the strongly increasing property is hereditary on the class of strongly increasing Puiseux monoids. We will require the following lemma.

\begin{lemma} \label{lem:no limit point implies strongly increasing subset}
	Let $R$ be an infinite subset of $\qq_{\ge 0}$. If $R$ does not contain any limit points, then it is the underlying set of a strongly increasing sequence.
\end{lemma}

\begin{proof}
	For every $r \in R$ and every subset $S$ of $R$, the interval $[0,r]$ must contain only finitely many elements of $S$; otherwise there would be a limit point of $S$ in $[0,r]$. Therefore every nonempty subset of $R$ has a minimum element. So the sequence $(r_n)_{n \in \mathbb{N}}$ recurrently defined by $r_1 = \min R$ and $r_n = \min R \! \setminus \! \{r_1, \dots, r_{n-1}\}$ is strictly increasing and has $R$ as its underlying set. Since $R$ is infinite and contains no limit points, the increasing sequence $(r_n)_{n \in \mathbb{N}}$ must be unbounded. Hence $R$ is the underlying set of the strongly increasing sequence $(r_n)_{n \in \mathbb{N}}$.
\end{proof}

\begin{theorem} \label{thm:strongly increasing iff super increasing}
	A nontrivial Puiseux monoid $M$ is strongly increasing if and only if every submonoid of $M$ is increasing.
\end{theorem}

\begin{proof}
	If $M$ is finitely generated, then it is isomorphic to a numerical monoid, and the statement of the theorem follows immediately. So we will assume for the rest of this proof that $M$ is not finitely generated. Suppose that $M$ is strongly increasing. Let us start by verifying that $M$ does not have any real limit points. By Proposition~\ref{prop:atoms of increasing monoids}, the monoid $M$ is atomic. As $M$ is atomic and non-finitely generated, $|\mathcal{A}(M)| = \infty$. Let $(a_n)_{n \in \mathbb{N}}$ be an increasing sequence with underlying set $\mathcal{A}(M)$. Since $M$ is strongly increasing and $\mathcal{A}(M)$ is an infinite subset contained in every generating set of $M$, the sequence $(a_n)_{n \in \mathbb{N}}$ is unbounded. Therefore, for every $r \in \rr$, the interval $[0,r]$ contains only finitely many elements of $(a_n)_{n \in \mathbb{N}}$, say $a_1, \dots, a_k$ for $k \in \nn$. Since $\langle a_1, \dots, a_k \rangle \cap [0,r]$ is a finite set, it follows that $M \cap [0,r]$ is finite as well. Because $|[0,r] \cap M| < \infty$ for all $r \in \rr$, it follows that $M$ does not have any limit points in $\rr$.
	
	Now suppose that $N$ is a nontrivial submonoid of $M$. Notice that, being a subset of $M$, the monoid $N$ cannot have any limit points in $\rr$. Thus, by Lemma~\ref{lem:no limit point implies strongly increasing subset}, the set~$N$ is the underlying set of a strongly increasing sequence of rationals. Hence $N$ is a strongly increasing Puiseux monoid, and the direct implication follows.
	
	For the converse implication, suppose that $M$ is not strongly increasing. We will check that, in this case, $M$ contains a submonoid that is not increasing. If $M$ is not increasing, then $M$ is a submonoid of itself that is not increasing. Suppose, therefore, that $M$ is increasing. By Proposition~\ref{prop:atoms of increasing monoids}, the monoid $M$ is atomic, and we can list its atoms increasingly. Let $(a_n)_{n \in \mathbb{N}}$ be an increasing sequence with underlying set $\mathcal{A}(M)$. Because $M$ is not strongly increasing, there exists a positive real $\ell$ that is the limit of the sequence $(a_n)_{n \in \mathbb{N}}$. Since $\ell$ is a limit point of $M$, which is closed under addition, it follows that $2\ell$ and $3\ell$ are both limit points of $M$. Let $(b_n)_{n \in \mathbb{N}}$ and $(c_n)_{n \in \mathbb{N}}$ be sequences in~$M$ having infinite underlying sets such that $\lim b_n = 2\ell$ and $\lim c_n = 3\ell$. Furthermore, assume that for each $n \in \nn$,
	\begin{equation}
		|b_n - 2\ell| < \frac{\ell}4 \ \text{ and } \ |c_n - 3\ell| < \frac{\ell}4. \label{eq:decreasing is not hereditary}
	\end{equation}
	Take $N$ to be the submonoid of $M$ generated by the set $A := \{b_n,c_n : n \in \nn\}$. Note that $A$ contains at least two limit points. Let us verify that $N$ is atomic with $\mathcal{A}(N) = A$. The inequalities \eqref{eq:decreasing is not hereditary} immediately imply that $A$ is bounded from above by $3\ell + \ell/4$. On the other hand, proving that $\mathcal{A}(N) = A$ amounts to showing that the sets $A$ and $A+A$ are disjoint. To verify this, it suffices to note that
	\begin{align*}
		\inf (A + A) &= \inf \big\{b_m + b_n, b_m + c_n, c_m + c_n : m,n \in \nn \big\} \\
		&\ge \min \bigg\{4\ell - \frac{\ell}{2}, \ 5\ell - \frac{\ell}{2}, \ 6\ell - \frac{\ell}{2} \bigg\} \\
		&> 3\ell + \frac{\ell}{4} \ge \sup A.
	\end{align*}
	Thus, $\mathcal{A}(N) = A$. Since every increasing sequence has at most one limit point in $\rr$, the set $A$ cannot be the underlying set of an increasing rational sequence. As every generating set of $N$ contains $A$, we conclude that $N$ is not an increasing Puiseux monoid, which completes the proof.
\end{proof}

As a direct consequence of Theorem \ref{thm:strongly increasing iff super increasing}, one obtains the following corollary.

\begin{corollary}
	Being atomic, increasing, and strongly increasing are hereditary properties on the class of strongly increasing Puiseux monoids.
\end{corollary}

The next result is important as it provides a large class of Puiseux monoids that are FFMs.

\begin{theorem} \label{thm:increasing PM of Archimedean fields are FF}
	Every increasing Puiseux monoid is an FFM.
\end{theorem}

\begin{proof}
	See \cite[Theorem~5.6]{fG19a}.
\end{proof}

On the other hand, the converse of Theorem~\ref{thm:increasing PM of Archimedean fields are FF} does not hold; the following example sheds some light upon this observation.

\begin{example} \label{ex:FF PM that is not increasing}
	Let $(p_n)_{n \in \nn}$ be a strictly increasing sequence of primes, and consider the Puiseux monoid of $\qq$ defined as follows:
	\begin{equation} \label{eq:FF does not implies increasing}
	M = \langle A \rangle, \ \text{where} \ A = \bigg\{\frac{p_{2n}^2 + 1}{p_{2n}}, \frac{p_{2n+1} + 1}{p_{2n+1}} : n \in \nn \bigg\}.
	\end{equation}
	Since $A$ is an unbounded subset of $\rr$ having $1$ as a limit point, it cannot be increasing. In addition, as $\mathsf{d}(a) \neq \mathsf{d}(a')$ for all $a,a' \in A$ such that $a \neq a'$, every element of $A$ is an atom of $M$. As each generating set of $M$ must contain $A$ (which is not increasing), $M$ is not an increasing Puiseux monoid.
	
	To verify that $M$ is an FFM, fix $x \in M$ and then take $D_x$ to be the set of prime numbers dividing $\mathsf{d}(x)$. Now choose $N \in \nn$ such that $N > \max \{x, \mathsf{d}(x)\}$. For each $a \in A$ with $\mathsf{d}(a) > N$, the number of copies $\alpha$ of the atom $a$ appearing in any $z \in \mathsf{Z}(x)$ must be a multiple of $\mathsf{d}(a)$ because $\mathsf{d}(a) \notin D_x$. Then $\alpha = 0$; otherwise, we would have that $x \ge \alpha a \ge \mathsf{d}(a)a > \mathsf{d}(a) > x$. Thus, if an atom $a$ divides $x$ in $M$, then $\mathsf{d}(a) \le N$. As a result, only finitely many elements of $\mathcal{A}(M)$ divide $x$ in $M$ and so $|\mathsf{Z}(x)| < \infty$. Hence $M$ is an FFM that is not increasing.
\end{example}

\begin{remark}
	For an ordered field $F$, a \emph{positive monoid} of $F$ is an additive submonoid of the nonnegative cone of $F$. As for Puiseux monoids, a positive monoid is \emph{increasing} if it can be generated by an increasing sequence. Increasing positive monoids are FFMs~\cite[Theorem~5.6]{fG19a}, but the proof of this general version of Theorem~\ref{thm:increasing PM of Archimedean fields are FF} is much more involved.
\end{remark}
\smallskip

\subsection{Decreasing Puiseux Monoids} \label{sec:decreasing PM}

Now that we have explored the structure of increasing Puiseux monoids, we will focus on the study of their decreasing counterpart.

A Puiseux monoid $M$ is said to be \emph{bounded} if $M$ can be generated by a bounded subset of rational numbers. Besides, we say that $M$ is \emph{strongly bounded} if $M$ can be generated by a set of rationals $R$ such that $\mathsf{n}(R)$ is bounded. If a Puiseux monoid is decreasing, then it is obviously bounded. On the other hand, there are bounded Puiseux monoids that are not even monotone; see Example~\ref{ex:bounded PM that is neither decreasing nor increasing}. However, every strongly bounded Puiseux monoid is decreasing, as we will show in Proposition~\ref{prop:SB PM are strongly decreasing}.

By contrast to the results we obtained in the previous section, the next proposition will show that being decreasing is almost never hereditary. In fact, we prove that being decreasing is hereditary only on those Puiseux monoids that are isomorphic to numerical monoids.

\begin{lemma} \label{lem:decreasing monoids have infinite limit points}
	If $M$ is a nontrivial decreasing Puiseux monoid, then exactly one of the following conditions holds:
	\begin{enumerate}
		\item $M$ is isomorphic to a numerical monoid;
		\vspace{3pt}
		\item $M$ contains infinitely many limit points in $\rr$.
	\end{enumerate}
\end{lemma}

\begin{proof}
	Suppose that $M$ is not isomorphic to a numerical monoid. Since $M$ is not trivial, it fails to be finitely generated. Therefore it can be generated by a strictly decreasing sequence $(a_n)_{n \in \mathbb{N}}$. The sequence $(a_n)_{n \in \mathbb{N}}$ must converge to a non-negative real number $\ell$. Since $(k a_n)_{n \in \mathbb{N}} \subseteq M$ converges to $k\ell$ for every $k \in \nn$, if $\ell \neq 0$, then every element of the infinite set $\{k\ell : k \in \nn\}$ is a limit point of $M$. On the other hand, if $\ell = 0$, then every term of the sequence $(a_n)_{n \in \mathbb{N}}$ is a limit point of $M$; this is because for every fixed $k \in \nn$ the sequence $(a_k + a_n)_{n \in \mathbb{N}} \subseteq M$ converges to $a_k$. Hence $M$ has infinitely many limit points in $\rr$.
	
	Now let us verify that at most one of the above two conditions can hold. For this, assume that $M$ is isomorphic to a numerical monoid. So $M$ is finitely generated, namely, $M = \langle r_1, \dots, r_n \rangle$, where $n \in \nn$ and $r_i \in \qq_{> 0}$ for $i=1,\dots, n$. For every $r \in \rr$ the interval $[0,r]$ contains only finitely many elements of $M$. Since $M \cap [0,r]$ is finite for all $r \in \rr$, it follows that $M$ cannot have any limit points in the real line.
\end{proof}

\begin{proposition}
	Let $M$ be a nontrivial decreasing Puiseux monoid. Then exactly one of the following conditions holds:
	\begin{enumerate}
		\item $M$ is isomorphic to a numerical monoid;
		\vspace{3pt}
		\item $M$ contains a submonoid that is not decreasing.
	\end{enumerate}
\end{proposition}

\begin{proof}
	Suppose that $M$ is not isomorphic to a numerical monoid. Let us construct a submonoid of $M$ that fails to be decreasing. Lemma~\ref{lem:decreasing monoids have infinite limit points} implies that $M$ has a nonzero limit point $\ell$. Since $M$ is closed under addition, $2\ell$ and $3\ell$ are both limit points of $M$. An argument as the one given in the proof of Theorem~\ref{thm:strongly increasing iff super increasing} will guarantee the existence of sequences $(a_n)_{n \in \mathbb{N}}$ and $(b_n)_{n \in \mathbb{N}}$ in $M$ having infinite underlying sets such that $(a_n)_{n \in \mathbb{N}}$ converges to $2\ell$, $(b_n)_{n \in \mathbb{N}}$ converges to $3\ell$, and the submonoid $N = \langle a_n, b_n : n \in \nn \rangle$ of $M$ is atomic with $\mathcal{A}(M) = \{a_n,b_n : n \in \nn\}$. Since every decreasing sequence of $\qq$ contains at most one limit point, $\mathcal{A}(M)$ cannot be the underlying set of a decreasing sequence of rationals. As every generating set of $N$ must contain $\mathcal{A}(M)$, we can conclude that $N$ is not decreasing. Hence at least one of the given conditions must hold.
	
	To see that both conditions cannot hold simultaneously, it suffices to observe that if~$M$ is isomorphic to a numerical monoid, then every nontrivial submonoid of~$M$ is also isomorphic to a numerical monoid and, therefore, decreasing.
\end{proof}

Similarly, as we did in the case of increasing Puiseux monoids, we will split the family of decreasing Puiseux monoids into two fundamental subfamilies, depending on whether $0$ is or is not a limit point. We say that a non-negative sequence of rationals is \emph{strongly decreasing} if it is decreasing and it converges to zero. A non-negative decreasing sequence of rationals converging to a positive real is called \emph{weakly decreasing}.

\begin{definition}
	A Puiseux monoid is \emph{strongly decreasing} if it can be generated by a strongly decreasing sequence of rational numbers. On the other hand, a Puiseux monoid is said to be \emph{weakly decreasing} if it can be generated by a weakly decreasing sequence of rationals.
\end{definition}

Observe that if a Puiseux monoid $M$ is weakly decreasing, then it has a generating sequence decreasing to a positive real number and, therefore, $0$ is not in the closure of~$M$. As a consequence, every weakly decreasing Puiseux monoid must be atomic. The next proposition describes those Puiseux monoids that are both strongly and weakly decreasing.

\begin{proposition}
	A decreasing Puiseux monoid is either strongly or weakly decreasing. A Puiseux monoid is both strongly and weakly decreasing if and only if it is isomorphic to a numerical monoid.
\end{proposition}

\begin{proof}
	As in the case of increasing Puiseux monoids, the first statement follows immediately. Now suppose that $M$ is a Puiseux monoid that is both strongly and weakly decreasing. Since $M$ is weakly decreasing, $0$ is not a limit point of $M^\bullet \!$. Let $(a_n)_{n \in \mathbb{N}}$ be a sequence decreasing to zero such that $M = \langle a_n : n \in \nn \rangle$. Because $0$ is not a limit point of $M^\bullet$, there exists $n_0 \in \nn$ such that $a_n = 0$ for all $n \ge n_0$. Hence $M$ is isomorphic to a numerical monoid. As in the increasing case, it is easily seen that every numerical monoid is both strongly and weakly decreasing.
\end{proof}

We mentioned at the beginning of this section that every strongly bounded Puiseux monoid is decreasing. Indeed, a stronger statement holds. 

\begin{proposition} \label{prop:SB PM are strongly decreasing}
	Every strongly bounded Puiseux monoid is strongly decreasing.
\end{proposition}

\begin{proof}
	Let $M$ be a strongly bounded Puiseux monoid. Since the trivial monoid is both strongly bounded and strongly decreasing, for this proof we will assume that $M \neq \{0\}$. Let $S \subset\qq_{> 0}$ be a generating set of $M$ such that $\mathsf{n}(S)$ is bounded. Since $\mathsf{n}(S)$ is finite, we can take $m$ to be the least common multiple of the elements of $\mathsf{n}(S)$. The map $x \mapsto \frac{1}{m}x$ is an order-preserving isomorphism from $M$ to $M' = \frac{1}{m}M$. Consequently, $M'$ is strongly decreasing if and only if $M$ is strongly decreasing. In addition, $S' = \frac{1}{m}S$ generates $M'$. Since $\mathsf{n}(S') = \{1\}$, it follows that $S'$ is the underlying set of a strongly decreasing sequence of rationals. Hence $M'$ is a strongly decreasing Puiseux monoid, which implies that $M$ is strongly decreasing as well.
\end{proof}

Recall that a Puiseux monoid $M$ is finite if $\pval(\mathsf{d}(M^\bullet)) = \{0\}$ for all but finitely many primes $p$. Strongly decreasing Puiseux monoids are not always strongly bounded, even if we require them to be finite. For example, if $r \in \qq$ such that $0 < r < 1$ and both $\mathsf{n}(r)$ and $\mathsf{d}(r)$ are different from $1$, then we have seen in Proposition~\ref{prop:atomic classification of multiplicative cyclic Puiseux monoids} that the Puiseux monoid $M_r = \langle r^n : n \in \nn_0 \rangle$ is atomic and $\mathcal{A}(M_r) = \{r^n : n \in \nn_0\}$. As a result, $M_r$ is finite and strongly decreasing. However, $M_r$ fails to be strongly bounded. On the other hand, not every bounded Puiseux monoid is decreasing, as illustrated in Example~\ref{ex:bounded PM that is neither decreasing nor increasing}.

Because numerical monoids are finitely generated, they are both increasing and decreasing Puiseux monoids. We end this section showing that numerical monoids are the only such Puiseux monoids.

\begin{proposition}
	A nontrivial Puiseux monoid $M$ is isomorphic to a numerical monoid if and only if $M$ is both increasing and decreasing.
\end{proposition}

\begin{proof}
	If $M$ is isomorphic to a numerical monoid, then it is finitely generated and, consequently, increasing and decreasing.
	
	Conversely, suppose that $M$ is a nontrivial Puiseux monoid that is increasing and decreasing. Proposition~\ref{prop:atoms of increasing monoids} implies that $M$ is atomic and, moreover, $\mathcal{A}(M)$ is the underlying set of an increasing sequence (because $\mathcal{A}(M) \neq \emptyset$). Suppose, by way of contradiction, that $\mathcal{A}(M)$ is not finite. In this case, $\mathcal{A}(M)$ does not contain a largest element. Since $M$ is decreasing, there exists $D = \{d_n : n \in \nn \} \subset \qq_{> 0}$ such that $d_1 > d_2 > \cdots$ and $M = \langle D \rangle$. Let $m = \min \{n \in \nn : d_n \in \mathcal{A}(M)\}$, which must exist because $\mathcal{A}(M) \subseteq D$. Since $\mathcal{A}(M)$ is contained in $D$, the minimality of $m$ implies that~$d_m$ is the largest element of $\mathcal{A}(M)$, which is a contradiction. Hence $\mathcal{A}(M)$ is finite. Since $M$ is atomic and $\mathcal{A}(M)$ is finite, $M$ is isomorphic to a numerical monoid.
\end{proof}
\medskip

\section{Factorial, Half-Factorial, and Other-Half-Factorial Monoids} The only Puiseux monoid that is a UFM (or even an HFM) is, up to isomorphism, $(\nn_0,+)$. The following proposition formalizes this observation.


\begin{proposition} \label{prop:HF PM characterization}
	For a nontrivial atomic Puiseux monoid $M$, the following statements are equivalent.
	\begin{enumerate}
		\item $M$ is a UFM.
		\vspace{2pt}
		
		\item $M$ is a HFM.
		\vspace{2pt}
		
		\item $M \cong (\nn_0,+)$.
		\vspace{2pt}
		
		\item $M$ contains a prime element.
	\end{enumerate}
\end{proposition}

\begin{proof}
	Clearly, (3) $\Rightarrow$ (1) $\Rightarrow$ (2). To argue (2) $\Rightarrow$ (3), assume that $M$ is an HFM. Since $M$ is an atomic nontrivial Puiseux monoid, $\mathcal{A}(M)$ is not empty. Let $a_1$ and~$a_2$ be two atoms of $M$. Then $z_1 := \mathsf{n}(a_2) \mathsf{d}(a_1) a_1$ and $z_2 :=  \mathsf{n}(a_1) \mathsf{d}(a_2) a_2$ are two factorizations of the element $\mathsf{n}(a_1) \mathsf{n}(a_2) \in M$. As $M$ is an HFM, it follows that $|z_1| = |z_2|$ and so $\mathsf{n}(a_2) \mathsf{d}(a_1) = \mathsf{n}(a_1)  \mathsf{d}(a_2)$. Then $a_1 = a_2$, which implies that $|\mathcal{A}(M)| = 1$. As a consequence, $M \cong (\nn_0,+)$. Because (3) $\Rightarrow$ (4) holds trivially, we only need to argue (4) $\Rightarrow$ (3). Fix a prime element $p \in M$ and take $a \in \mathcal{A}(M)$. Since $p \mid_M \mathsf{n}(p) \mathsf{d}(a) a$, one finds that $p \mid_M a$. This, in turn, implies that $a = p$. Hence $\mathcal{A}(M) = \{p\}$, and so $M \cong (\nn_0,+)$.
\end{proof}

\begin{example}
	For general monoids, the property of being an HFM is strictly weaker than that of being a UFM. For instance, the submonoid $\langle (1,n) \mid n \in \nn \rangle$ of $\nn_0^2$ is an HFM that is not a UFM (see \cite[Propositions~5.1 and~5.4]{fG19d} for more details).
\end{example}

A dual notion of being an HFM was introduced in~\cite{CS11} by Coykendall and Smith.

\begin{definition}
	An atomic monoid $M$ is an \emph{OHFM} (or an \emph{other-half-factorial monoid}) if for all $x \in M \setminus U(M)$ and $z, z' \in \mathsf{Z}(x)$ with $|z| = |z'|$, we have $z = z'$.
\end{definition}

Clearly, every UFM is an OHFM. Although the multiplicative monoid of an integral domain is a UFM if and only if it is an OHFM~\cite[Corollary~2.11]{CS11}, OHFMs are not always UFMs or HFMs, even in the class of Puiseux monoids.

\begin{proposition} \label{prop:OHF PM characterization}
	For a nontrivial atomic Puiseux monoid $M$, the following conditions are equivalent.
	\begin{enumerate}
		\item $M$ is an OHFM.
		\vspace{2pt}
		
		\item $|\mathcal{A}(M)| \le 2$.
		\vspace{2pt}
		
		\item $M$ is isomorphic to a numerical monoid with embedding dimension in $\{1,2\}$.
	\end{enumerate}
\end{proposition}

\begin{proof}
	To prove (1) $\Rightarrow$ (2), let $M$ be an OHFM. If $M$ is factorial, then $M \cong (\nn_0,+)$, and we are done. Then suppose that $M$ is not factorial. In this case, $|\mathcal{A}(M)| \ge 2$. Assume for a contradiction that $|\mathcal{A}(M)| \ge 3$. Take $a_1, a_2, a_3 \in \mathcal{A}(M)$ satisfying that $a_1 < a_2 < a_3$. Let $d = \mathsf{d}(a_1) \mathsf{d}(a_2) \mathsf{d}(a_3)$, and set $a'_i = d a_i$ for each $i \in \{ 1, \dots, 3 \}$. Since $a'_1, a'_2$, and $a'_3$ are integers satisfying that $a'_1 < a'_2 < a'_3$, there exist $m,n \in \nn$ such that
	\begin{equation} \label{eq:OHF}
	m(a'_2 - a'_1) = n(a'_3 - a'_2).
	\end{equation}
	Clearly, $z_1 := ma_1 + na_3$ and $z_2 := (m+n)a_2$ are two distinct factorizations in $\mathsf{Z}(M)$ satisfying that $|z_1| = m+n = |z_2|$. In addition, after dividing both sides of the equality~(\ref{eq:OHF}) by $d$, one obtains $ma_1 + na_3 = (m+n)a_2$, which means that $z_1$ and $z_2$ are factorizations of the same element. However, this contradicts that $M$ is an OHFM. Hence $|\mathcal{A}(M)| \le 2$, as desired.

	To show that (2) $\Rightarrow$ (3), suppose that $|\mathcal{A}(M)| \le 2$. By~Proposition~\ref{prop:fg PM are NM}, $M$ is isomorphic to a numerical monoid $N$. As $|\mathcal{A}(M)| \le 2$, the embedding dimension of $N$ belongs to $\{1,2\}$, as desired.
	
	To show (3) $\Rightarrow$ (1), suppose that either $M \cong (\nn_0,+)$ or $M \cong \langle a, b \rangle$ for $a, b \in \nn_{\ge 2}$ with $\gcd(a,b) = 1$. If $M \cong (\nn_0,+)$, then~$M$ is factorial and, in particular, an OHFM. On the other hand, if $M \cong \langle a, b \rangle$, then it is an OHFM by~\cite[Example~2.13]{CS11}.
\end{proof}

\begin{remark}
	There are Puiseux monoids that are FFMs but neither HFMs nor OHFMs. As a direct consequence  of Theorem~\ref{thm:increasing PM of Archimedean fields are FF} and Propositions~\ref{prop:HF PM characterization} and~\ref{prop:OHF PM characterization}, one finds that $\langle \frac{p-1}{p} : p \in \pp \rangle$ is one of such monoids.
\end{remark}

We started this chapter exhibiting the following chain of implications, which holds for all monoids:
\begin{equation*}
\textbf{UFM} \ \Rightarrow \ \textbf{HFM} \stackrel{\text{reduced}}{\Rightarrow} \ \textbf{FFM} \ \Rightarrow \ \textbf{BFM} \ \Rightarrow \ \textbf{ACCP} \ \Rightarrow \ \textbf{atomic monoid}
\end{equation*}
This chain was first introduced by Anderson et al. in the context of integral domains~\cite{AAZ90}. As we mentioned at the beginning of the chapter, none of the implications above is reversible in the context of integral domains. In the context of Puiseux monoids, we have established the following chain of implications:
\[
	(\textbf{UFM} \ \Leftrightarrow \ \textbf{HFM}) \ \Rightarrow \ \textbf{OHFM} \ \Rightarrow \ \textbf{FGM} \ \Rightarrow \ \textbf{FFM} \ \Rightarrow \ \textbf{BFM} \ \Rightarrow \ \textbf{ACCP} \ \Rightarrow \ \textbf{AM},
\]
where FGM and AM stand for finitely generated monoid and atomic monoid, respectively. We have also provided examples to illustrate that, except UFM $\Leftrightarrow$ HFM, none of the implications in the previous chain is reversible in the context of Puiseux monoids.
\smallskip


%% file: tex/ch3.tex
\chapter{Bounded and Finite Puiseux Monoids} \label{chap:density and p-adics}
\smallskip

\section{Bounded and Strongly Bounded Puiseux Monoids} We begin this chapter exploring two properties of Puiseux monoids: being bounded and being strongly bounded. Although these two properties are neither atomic nor algebraic, they often help to understand the structure of Puiseux monoids. For instance, they can be used to give another characterization of finitely generated Puiseux monoids (see Theorem~\ref{theo:bounded P-RM are NS or non-atomic} below).

Recall that a Puiseux monoid $M$ is said to be bounded if $M$ can be generated by a bounded subset of rational numbers and strongly bounded if $M$ can be generated by a set of rationals $R$ such that $\mathsf{n}(R)$ is bounded. There are Puiseux monoids that are not bounded, as we shall see below.
\smallskip

\subsection{Boundedness of Prime Reciprocal Puiseux Monoids} Let us proceed to study the boundedness of prime reciprocal Puiseux monoids. Obviously, every finitely generated Puiseux monoid is strongly bounded. However, there are strongly bounded Puiseux monoids that fail to be finitely generated. Clearly, the family of strongly bounded Puiseux monoids is strictly contained in that of bounded Puiseux monoids. The following example illustrates these observations.

\begin{example}
	Consider the Puiseux monoid
	\[
		M_1 = \langle A_1 \rangle, \ \text{ where } \ A_1 = \bigg\{\frac{p^2-1}{p}  : p \in \pp \bigg\}.
	\]
	An elementary divisibility argument will reveal that $\mathcal{A}(M_1) = A_1$. Because $A_1$ is an unbounded set, it follows that $M_1$ is not a bounded Puiseux monoid. On the other hand, let us consider the strongly bounded Puiseux monoid
	\[
		M_2 = \langle A_2 \rangle, \ \text { where } \ A_2 = \bigg\{\frac{1}{p} : p \in \pp \bigg\}.
	\]
	As in the previous case, it is not hard to verify that $\mathcal{A}(M_2) = A_2$. Therefore $M_2$ is a strongly bounded Puiseux monoid that is not finitely generated. Finally, consider the bounded Puiseux monoid
	\[
		M_3 = \langle A_3 \rangle, \ \text{ where } \ A_3 = \bigg\{\frac{p-1}{p} : p \in \pp \bigg\}.
	\]
	Once again, $\mathcal{A}(M_3) = A_3$ follows from an elementary divisibility argument. As a result,~$M_3$ cannot be strongly bounded.
\end{example}

Let $M$ be a Puiseux monoid, and let $N$ be a submonoid of $M$. If $M$ is finitely generated, then $N$ is also finitely generated. Thus, being finitely generated is hereditary on the class of finitely generated Puiseux monoids. As we should expect, not every property of a Puiseux monoid is inherited by all its submonoids. For example, being antimatter is not hereditary on the class of antimatter Puiseux monoids; for instance, $\qq_{\ge 0}$ is antimatter, but it contains the atomic additive submonoid $\nn_0$, which satisfies $\mathcal{A}(\nn_0) = \{1\}$. Moreover, as Corollary~\ref{cor:strongly bounded is not hereditary on primary PM} indicates, boundedness and strong boundedness are not hereditary, even on the class of prime reciprocal Puiseux monoids.

Let $S$ be a set of naturals. If the series $\sum_{s \in S} 1/s$ diverges, $S$  is said to be \emph{substantial}. If $S$ is not substantial, it is said to be \emph{insubstantial} (see~\cite[Chapter~10]{pC}). For example, it is well known that the set of prime numbers is substantial as it was first noticed by Euler that the series of reciprocal primes is divergent.

\begin{proposition} \label{prop:strongly bounded PM having a submonoid that is not even bounded}
	Let $P$ be a set of primes, and let $M$ be the prime reciprocal Puiseux monoid $\langle 1/p : p \in P \rangle$. If every submonoid of $M$ is bounded, then $P$ is insubstantial.
\end{proposition}

\begin{proof}
	Suppose, by way of contradiction, that $P$ is substantial. Then $P$ must contain infinitely many primes. Let $(p_n)_{n \in \mathbb{N}}$ be a strictly increasing enumeration of the elements in $P$. Take $N$ to be the submonoid of $M$ generated by $A = \{a_n : n \in \nn\}$, where
	\[
		a_n = \sum_{i=1}^n \frac{1}{p_i}.
	\]
	Since $P$ is substantial, $A$ is unbounded. We will show that $N$ fails to be bounded. For this purpose, we verify that $\mathcal{A}(N) = A$, which implies that every generating set of $N$ contains $A$ and, therefore, must be unbounded. Suppose that
	\begin{equation} \label{eq:boundedness in not hereditary 1}
		a_n = a_{n_1} + \dots + a_{n_{\ell}}
	\end{equation}
	for some $\ell, n, n_1, \dots, n_\ell \in \nn$ such that $n_1 \le \dots \le n_\ell$. Since $(a_n)_{n \in \mathbb{N}}$ is an increasing sequence, $n \ge n_\ell$. After multiplying the equation \eqref{eq:boundedness in not hereditary 1} by $m = p_1 \dots p_n$ and moving every summand but $m/p_n$ to the right-hand side, we obtain
	\begin{equation} \label{eq:boundedness in not hereditary 2}
		p_1 \dots p_{n-1} = \sum_{j=1}^\ell \sum_{i=1}^{n_j} \frac{m}{p_i} - \sum_{i=1}^{n-1} \frac{m}{p_i}.
	\end{equation}
	Now we observe that if $n$ were strictly greater than $n_\ell$, then $m/p_i$ would be an integer divisible by $p_n$ for each $i = 1, \dots, n_j$ and $j = 1, \dots, \ell$, which would imply that the right-hand side of \eqref{eq:boundedness in not hereditary 2} is divisible by $p_n$. This cannot be possible because $p_1 \dots p_{n-1}$ is not divisible by $p_n$. Thus, $n = n_\ell$ and so $\ell = 1$. Since $(a_n)_{n \in \mathbb{N}}$ is an increasing sequence satisfying that $a_n \notin \langle a_1, \dots, a_{n-1} \rangle$, Proposition~\ref{prop:atoms of increasing monoids} ensures that $\mathcal{A}(N) = A$. As a result, $M$ contains a submonoid that fails to be bounded; but this is a contradiction. Hence the set $P$ is insubstantial.
\end{proof}

For $m,n \in \nn_0$ such that $n > 0$ and $\gcd(m,n) = 1$, Dirichlet's theorem states that the set $P$ of all primes $p$ satisfying that $p \equiv m \pmod n$ is infinite. For a relatively elementary proof of Dirichlet's theorem, see \cite{hS50}. Furthermore, it is also known that the set $P$ is substantial; indeed, as indicated in \cite[page 156]{tA76}, there exists a constant $A$ for which
\begin{equation} \label{eq:prime in arithmetic progression are substantial}
	\sum_{p \in P, p \le x} \frac{1}{p} = \frac{1}{\varphi(n)} \log \log x + A + O\bigg(\frac{1}{\log x}\bigg),
\end{equation}
where $\varphi$ is the Euler totient function. In particular, the set comprising all primes of the form $4k+1$ (or $4k+ 3$) is substantial. The next corollary follows immediately from Proposition~\ref{prop:strongly bounded PM having a submonoid that is not even bounded} and equation~\eqref{eq:prime in arithmetic progression are substantial}.

\begin{corollary} \label{cor:strongly bounded is not hereditary on primary PM}
	Let $m,n \in \nn_0$ such that $n > 0$ and $\gcd(m,n) = 1$, and let $P$ be the set of all primes $p$ satisfying $p \equiv m \pmod n$. Then the prime reciprocal Puiseux monoid $M = \langle 1/p : p \in P \rangle$ contains an unbounded submonoid.
\end{corollary}
\smallskip

\subsection{Boundedness of Multiplicatively Cyclic Puiseux Monoids} 

Let us turn to study the boundedness of the multiplicatively cyclic Puiseux monoids.

\begin{proposition}
	For $r \in \qq_{> 0}$, let $M_r$ be the multiplicatively $r$-cyclic Puiseux monoid. Then the following statements hold.
	\begin{enumerate}
		\item If $\mathsf{n}(r)=1$ or $\mathsf{d}(r)=1$, then $M_r$ is strongly bounded.
		
		\item If $\mathsf{n}(r), \mathsf{d}(r) > 1$ and $r < 1$, then $M_r$ is bounded but not strongly bounded.

		\item If $\mathsf{n}(r), \mathsf{d}(r)  > 1$ and $r > 1$, then $M_r$ is not bounded.
	\end{enumerate}
\end{proposition}

\begin{proof}
	1. Set $A_r := \{r^n : n \in \nn_0\}$. If $\mathsf{n}(r) = 1$, then $\mathsf{n}(A_r) = \{1\}$ and, therefore, $M_r$ is strongly bounded. On the other hand, if $\mathsf{d}(r) = 1$, then $M_r$ is finitely generated, and so strongly bounded.
	
	 2. Suppose that $\mathsf{n}(r), \mathsf{d}(r) > 1$ and $r < 1$. It follows from Proposition~\ref{prop:atomic classification of multiplicative cyclic Puiseux monoids} that $M_r$ is atomic with $\mathcal{A}(M_r) = A_r$. Since $r < 1$ the generating set $A_r$ is bounded, which means that $M_r$ is bounded. However $M_r$ cannot be strongly bounded because every generating set of $M_r$ must contain the set $A_n$, whose numerator set $\mathsf{n}(A_r) = \{\mathsf{n}(r)^n : n \in \nn_0\}$ is unbounded.
	 
	 3. Finally, suppose that $\mathsf{n}(r), \mathsf{d}(r) > 1$ and $r > 1$. As in the previous part, it follows from Proposition~\ref{prop:atomic classification of multiplicative cyclic Puiseux monoids} that $M_r$ is atomic with $\mathcal{A}(M_r) = A_r$. As $A_r$ is unbounded and any generating set of $M_r$ must contain $A_r$, we obtain that $M_r$ is not bounded.
\end{proof}

As illustrated by Corollary~\ref{cor:strongly bounded is not hereditary on primary PM}, being bounded (or strongly bounded) is not hereditary on the class of prime reciprocal Puiseux monoids. Additionally, boundedness (resp., strong boundedness) is not hereditary on the class of bounded (resp., strongly bounded) multiplicatively cyclic Puiseux monoids.

\begin{example}
	Let $M$ be the multiplicatively (1/2)-cyclic Puiseux monoid, that is, $M = \langle 1/2^n : n \in \nn_0 \rangle$. It is strongly bounded, and yet its submonoid
	\begin{equation} \label{eq:strong boundedness is not hereditary even in multiplicatively cyclic PM}
	N = \bigg\langle \sum_{i=1}^n \frac{1}{2^i} : n \in \nn \bigg\rangle = \bigg\langle \frac{2^n - 1}{2^n} : n \in \nn \bigg\rangle
	\end{equation}
	is not strongly bounded; to see this, it is enough to verify that $\mathcal{A}(N) = S$, where $S$ is the generating set defining $N$ in \eqref{eq:strong boundedness is not hereditary even in multiplicatively cyclic PM}. Note that the sum of any two elements of the generating set $S$ is at least one, while every element of $S$ is less than one. Therefore each element of $S$ must be an atom of $N$, and so $\mathcal{A}(N) = S$.
\end{example}

We argue now that boundedness (resp., strong boundedness) is almost never hereditary on the class of bounded (resp., strongly bounded) multiplicatively cyclic Puiseux monoids.

\begin{proposition}
	For $r \in \qq_{> 0}$, let $M_r$ be the multiplicatively $r$-cyclic Puiseux monoid. Then every submonoid of $M$ is bounded (or strongly bounded) if and only if $M_r$ is isomorphic to a numerical monoid.
\end{proposition}

\begin{proof}
	Let $a$ and $b$ denote $\mathsf{n}(r)$ and $\mathsf{d}(r)$, respectively. To prove the direct implication, suppose, by way of contradiction, that $M_r$ is not isomorphic to a numerical monoid. In this case, $b > 1$. Consider the submonoid $N = \langle s_1, s_2, \dots \rangle$ of $M$, where
	\[
		s_n = \frac{(nb^n + 1)a^n}{b^n}
	\]
	for every natural $n$. Proving the forward implication amounts to verifying that $N$ is not bounded and, as a consequence, not strongly bounded. First, let us check that $\mathcal{A}(N) = \{s_n : n \in \nn\}$. Note that
	\[
		s_{n+1} = \frac{((n+1)b^{n+1} + 1)a^{n+1}}{b^{n+1}} > \frac{(nb^{n+1} + b)a^n}{b^{n+1}} = s_n
	\]
	for each $n \in \nn$, and so $(s_n)_{n \in \mathbb{N}}$ is an increasing sequence. Moreover, it is easy to see that $s_n > n$ for every $n$. Thus, $(s_n)_{n \in \mathbb{N}}$ is unbounded. Suppose that there exist $k, \alpha_k \in \nn$ and $\alpha_i \in \nn_0$ for every $i = 1, \dots, k-1$ such that $s_n = \alpha_1s_1 + \dots + \alpha_k s_k$. Since $(s_n)_{n \in \mathbb{N}}$ is increasing and $\alpha_k > 0$, we have $k \le n$. Let $p$ be a prime divisor of $b$, and let $m = \pval(b)$. The fact that $p \nmid (nb^n + 1)a^n$ for every natural $n$ implies $\pval(s_n) = -mn$. Therefore
	\[
		-mn = \pval(s_n) \ge \min_{1 \le i \le k}\{\pval(\alpha_is_i)\} \ge \min_{1 \le i \le k}\{\pval(s_i)\} = -mk,
	\]
	which implies that $k \ge n$. Thus, $k=n$ and then $\alpha_1 = \dots = \alpha_{n-1} = 0$ and $\alpha_n = 1$. So $s_n \notin \langle s_1, \dots, s_{n-1} \rangle$ for every $n \in \nn$ and, by Proposition~\ref{prop:atoms of increasing monoids}, $\mathcal{A}(N) = \{s_n : n \in \nn\}$. Since $\mathcal{A}(N)$ is unbounded, $N$ cannot be a bounded Puiseux monoid.
	
	On the other hand, if $M_r$ is isomorphic to a numerical monoid, then it is finitely generated and, hence, bounded and strongly bounded. This gives us the converse implication.
\end{proof}
\medskip

\section{Finite Puiseux Monoids}

\subsection{Connection with Strong Boundedness} A Puiseux monoid $M$ \emph{over} $P$ is a Puiseux monoid such that $\pval(m) \ge 0$ for every $m \in M$ and $p \notin P$. If $P$ is finite, then we say that $M$ is a \emph{finite} Puiseux monoid over $P$. The Puiseux monoid $M$ is said to be \emph{finite} if there exists a finite set of primes $P$ such that $M$ is finite over $P$. One can use strongly boundedness and finiteness to characterize finitely generated Puiseux monoids. The following result was established in~\cite{fG17}.

\begin{theorem} \cite[Theorem~5.8]{fG17} \label{theo:bounded P-RM are NS or non-atomic}
	For a Puiseux monoid $M$, the following statements are equivalent.
	\begin{enumerate}
		\item $M$ is atomic, strongly bounded, and finite.
		\item $M$ is finitely generated.
	\end{enumerate}
\end{theorem}

No two of the three conditions in statement~(1) of Theorem~\ref{theo:bounded P-RM are NS or non-atomic} imply the third one. In addition, the three conditions are required for the equivalence to hold. The following example illustrates these observations.

\begin{example} \label{ex:introductory example}
	For $q \in \pp$, consider the Puiseux monoid
	\[
		M_1 = \bigg\langle \frac 1{q^n} : n \in \nn \bigg\rangle.
	\]
	It follows immediately that $M_1$ is a nontrivial Puiseux monoid that is finite, strongly bounded, and antimatter. In particular, $M_1$ is not atomic. Since $M_1$ is not atomic, it cannot be finitely generated.

	Now consider the Puiseux monoid
	\[
		M_2 = \bigg\langle \bigg( \frac{q+1}{q} \bigg)^n : n \in \nn_0 \bigg\rangle.
	\]
	It is clear that $M_2$ is finite. It follows from Proposition~\ref{prop:atomic classification of multiplicative cyclic Puiseux monoids} that $M_2$ is an atomic Puiseux monoid with set of atoms $\mathcal{A}(M_2) = \big\{ \big( \frac{q+1}{q} \big)^n : n \in \nn_0 \big\}$. As a result, $M$ is not strongly bounded. Therefore $M_2$ is an atomic and finite Puiseux monoid that is not strongly bounded. Since $M_2$ is not strongly bounded, it is not finitely generated.

	Finally, we have already seen that the prime reciprocal Puiseux monoid
	\[
		M_3 = \bigg\langle \frac{1}{p} : p \in \pp \bigg\rangle
	\]
	is atomic with sets of atoms $\mathcal{A}(M_3) = \{ 1/p : p \in \pp \}$. Then $M_3$ is a strongly bounded atomic Puiseux monoid, but it is not finite. In particular, it is not finitely generated, which concludes our set of examples.
\end{example}
\medskip

\subsection{$p$-adic Puiseux Monoids} \label{sec:p-adic PM} Let $M$ be a Puiseux monoid. We have seen in Chapter~\ref{chap:atomicity} that when $0$ is not a limit point of $M^\bullet$, then $M$ is a BFM and, in particular, an atomic monoid. When $0$ is a limit point of $M^\bullet$, the situation is significantly more subtle, and there is not a general criterion to decide whether such Puiseux monoids are atomic (or BFMs). Some of the simplest representatives of this class are the Puiseux monoids $M_p = \langle 1/p^n : n \in \nn \rangle$ for $p \in \pp$, which happen to be antimatter. However, $M_p$ contain plenty of submonoids with a very diverse atomic structure. In this section we delve into the atomicity of submonoids of $M_p$.

\begin{definition}
	Let $p$ be a prime. We say that a Puiseux monoid $M$ is \emph{$p$-adic} if $\mathsf{d}(x)$ is a power of $p$ for all $x \in M^\bullet$.
\end{definition}

We use the term $p$-adic monoid as a short for $p$-adic Puiseux monoid. Throughout this section, every time that we define a $p$-adic monoid by specifying a sequence of generators $(r_n)_{n \in \mathbb{N}}$, we shall implicitly assume that $(\mathsf{d}(r_n))_{n \in \mathbb{N}}$ increases to infinity; this assumption comes without loss of generality because in order to generate a Puiseux monoid we only need to repeat each denominator finitely many times. On the other hand, $\lim \mathsf{d}(r_n) = \infty$ does not affect the generality of the results we prove in this section for if $(\mathsf{d}(r_n))_{n \in \mathbb{N}}$ is a bounded sequence, then the $p$-adic monoid generated by $(r_n)_{n \in \mathbb{N}}$ is finitely generated and, therefore, isomorphic to a numerical monoid.

Strongly bounded $p$-adic monoids happen to have only finitely many atoms (cf.~Theorem~\ref{theo:bounded P-RM are NS or non-atomic}), as revealed by the next proposition.

\begin{proposition} \label{prop:strongly bounded p-adic monoids has finitely many atoms}
	A strongly bounded $p$-adic monoid has only finitely many atoms.
\end{proposition}

\begin{proof}
	For $p \in \pp$, let $M$ be a strongly bounded $p$-adic monoid. Let $(r_n)_{n \in \mathbb{N}}$ be a generating sequence for $M$ with underlying set $R$ satisfying that $\mathsf{n}(R) = \{n_1, \dots, n_k\}$ for some $k, n_1, \dots, n_k \in \nn$. For each $i \in \{1,\dots, k\}$, take $R_i = \{r_n : \mathsf{n}(r_n) = n_i\}$ and $M_i = \langle R_i \rangle$. The fact that $R \subseteq M_1 \cup \dots \cup M_k$, along with $\mathcal{A}(M) \cap M_i \subseteq \mathcal{A}(M_i)$, implies that
	\[
		\mathcal{A}(M) \subseteq \bigcup_{i=1}^k \mathcal{A}(M_i).
	\]
	Thus, showing that $\mathcal{A}(M)$ is finite amounts to verifying that $|\mathcal{A}(M_i)| < \infty$ for each $i=1,\dots,k$. Fix $i \in \{1,\dots, k\}$. If $M_i$ is finitely generated, then $|\mathcal{A}(M_i)| < \infty$. Let us assume, therefore, that $M_i$ is not finitely generated. This means that there exists a strictly increasing sequence $(\alpha_n)_{n \in \mathbb{N}}$ such that $M_i = \langle n_i/p^{\alpha_n} : n \in \nn \rangle$. Because $n_i/p^{\alpha_n} = p^{\alpha_{n+1} - \alpha_n}(n_i/p^{\alpha_{n+1}})$, the monoid $M_i$ satisfies that $|\mathcal{A}(M_i)| = 0$. Hence we conclude that $\mathcal{A}(M)$ is finite.
\end{proof}

We are now in a position to give a necessary condition for the atomicity of $p$-adic monoids.

\begin{theorem} \label{thm:necessary condition for the atomicity of p-adic PM}
	Let $p \in \pp$, and let $M$ be an atomic $p$-adic monoid satisfying that $\mathcal{A}(M) = \{r_n : n\in \nn\}$. If $\lim r_n = 0$, then $\lim \mathsf{n}(r_n) = \infty$.
\end{theorem}

\begin{proof}
	Set $a_n = \mathsf{n}(r_n)$ and $p^{\alpha_n} = \mathsf{d}(r_n)$ for every natural $n$. Suppose, by way of contradiction, that $\lim a_n \neq \infty$. Then there exists $m \in \nn$ such that $a_n = m$ for infinitely many $n \in \nn$. For each positive divisor $d$ of $m$ we define the Puiseux monoid
	\[
		M_d = \langle S_d \rangle, \ \text{ where } \ S_d =  \bigg\{ \frac{a_{k_n}}{p^{\alpha_{k_n}}} : a_{k_n} = m \ \text{or} \ \gcd(m, a_{k_n}) = d \bigg\}.
	\]
	Observe that $\mathcal{A}(M)$ is included in the union of the $M_d$. On the other hand, the fact that $\mathcal{A}(M) \cap M_d \subseteq \mathcal{A}(M_d)$ for every $d$ dividing $m$ implies that
	\begin{equation} \label{eq:inclusion of atoms}
	\mathcal{A}(M) \subseteq \bigcup_{d \mid m} \mathcal{A}(M_d).
	\end{equation}
	Because $\mathcal{A}(M)$ contains infinitely many atoms, the inclusion \eqref{eq:inclusion of atoms} implies the existence of a divisor $d$ of $m$ such that $|\mathcal{A}(M_d)| = \infty$. Set $N_d = \frac 1d M_d$. Since $d$ divides $\mathsf{n}(q)$ for all $q \in M_d$, it follows that $N_d$ is also a $p$-adic monoid. In addition, the fact that $N_d$ is isomorphic to $M_d$ implies that $|\mathcal{A}(N_d)| = |\mathcal{A}(M_d)| = \infty$. After setting $b_n = a_{k_n}/d$ and $\beta_n = \alpha_{k_n}$ for every natural $n$ such that either $a_{k_n} = m$ or $\gcd(m, a_{k_n}) = d$, we have
	\[
		N_d = \bigg\langle \frac{b_n}{p^{\beta_n}} : \ n \in \nn \bigg\rangle.
	\]
	As $a_n = m$ for infinitely many $n \in \nn$, the sequence $(\beta_n)_{n \in \mathbb{N}}$ is an infinite subsequence of $(\alpha_n)_{n \in \mathbb{N}}$ and, therefore, it increases to infinity. In addition, as $\lim a_n/p^{\alpha_n} = 0$, it follows that $\lim b_n/p^{\beta_n} = 0$.

	Now we argue that $\mathcal{A}(N_d)$ is finite, which will yield the desired contradiction. Take $m' = m/d$. Since there are infinitely many $n \in \nn$ such that $b_n = m'$, it is guaranteed that $m'/p^n \in N_d$ for every $n \in \nn$. In addition, $\gcd(m', b_n) = 1$ for each $b_n \neq m'$. If $b_n \neq m'$ for only finitely many $n$, then $N_d$ is strongly bounded and Proposition~\ref{prop:strongly bounded p-adic monoids has finitely many atoms} ensures that $\mathcal{A}(N_d)$ is finite. Suppose otherwise that $\gcd(b_n, m') = 1$ (i.e., $b_n \neq m'$) for infinitely many $n \in \nn$. For a fixed $i$ with $b_i \neq m'$ take $j \in \nn$ satisfying that $\gcd(b_j, m') = 1$ and large enough so that $b_i p^{\beta_j - \beta_i} > b_j m'$; the existence of such an index $j$ is guaranteed by the fact that $\lim b_n/p^{\beta_n} = 0$. As $b_i p^{\beta_j - \beta_i} > b_j m' > \mathfrak{f}(\langle b_j, m'\rangle)$, there exist positive integers $x$ and $y$ such that $b_i p^{\beta_j - \beta_i} = x b_j + y m'$, that is
	\[
		\frac{b_i}{p^{\beta_i}} = x \frac{b_j}{p^{\beta_j}} + y \frac{m'}{p^{\beta_j}}.
	\]
	As $b_j/p^{\beta_j}, m'/p^{\beta_j} \in N_d^\bullet$, it follows that $b_i/p^{\beta_i} \notin \mathcal{A}(N_d)$. Because $i$ was arbitrarily taken, $N_d$ is antimatter. In particular, $\mathcal{A}(N_d)$ is finite, which leads to a contradiction.
\end{proof}

The conditions $\lim r_n = 0$ and $\lim \mathsf{n}(r_n) = \infty$ are not enough to guarantee that the non-finitely generated $p$-adic monoid $M$ is atomic. The next example sheds some light upon this observation.

\begin{example}
	For an odd prime $p$, consider the $p$-adic monoid
	\begin{equation} \label{eq:p-adic not atomic}
		M = \bigg\langle \frac{p^{2^n} - 1}{p^{2^{n+1}}}, \frac{p^{2^n} + 1}{p^{2^{n+1}}} : n \in \nn \bigg\rangle.
	\end{equation}
	Observe that the sequence of numerators $(p^{2^n} - 1, p^{2^n} + 1)_{n \in \mathbb{N}}$ increases to infinity while the sequence of generators of $M$ converges to zero. Also, notice that for every $n \in\nn$,
	\[
		\frac{2}{p^{2^n}} = \frac{p^{2^n} - 1}{p^{2^{n+1}}} + \frac{p^{2^n} + 1}{p^{2^{n+1}}} \in M.
	\]
	Now we can see that $M$ is not atomic; indeed, $M$ is antimatter, which immediately follows from the fact that
	\[
		\frac{p^{2^n} \pm 1}{p^{2^{n+1}}} = \frac{p^{2^n} \pm 1}2 \frac2{p^{2^{n+1}}}.
	\]
\end{example}
\vspace{4pt}
The next proposition yields a necessary and a sufficient condition for the atomicity of $p$-adic monoids having generating sets whose numerators are powers of the same prime.

\begin{proposition}\label{prop:atomicity p-adic}
	Let $p$ and $q$ be two different primes, and let $M = \langle r_n : n \in \nn \rangle$ be a $p$-adic monoid such that $\mathsf{n}(r_n)$ is a power of $q$ for every $n \in \nn$. Then
	\vspace{-3pt}
	\begin{enumerate}
		\item if $M$ is atomic, then $\lim \mathsf{n}(r_n) = \infty$;
		\vspace{2pt}
		\item if $\lim \mathsf{n}(r_n) = \infty$ and $(r_n)_{n \in \mathbb{N}}$ is decreasing, then $M$ is atomic.
	\end{enumerate}
\end{proposition}

\begin{proof}
	Define the sequences $(\alpha_n)_{n \in \mathbb{N}}$ and $(\beta_n)_{n \in \mathbb{N}}$ such that $p^{\alpha_n} = \mathsf{d}(r_n)$ and $q^{\beta_n} = \mathsf{n}(r_n)$. To check condition (1), suppose, by way of contradiction, that $\lim \mathsf{n}(r_n) \neq \infty$. Therefore there is a natural $j$ such that $\mathsf{n(r_n)} = q^j$ for infinitely many $n \in \nn$. This implies that $q^j/p^n \in M$ for every $n \in \nn$. Thus, for every $x \in M^\bullet$ such that $\mathsf{n}(x) = q^m \ge q^j$, one can write
	\[
		x = \frac{q^m}{\mathsf{d}(x)} = pq^{m-j} \frac{q^j}{p \mathsf{d}(x)} \notin \mathcal{A}(M).
	\]
	As a result, every $a \in \mathcal{A}(M)$ satisfies that $\mathsf{n}(a) < q^j$. This immediately implies that $\mathcal{A}(M)$ is finite. As $M$ is atomic with $|\mathcal{A}(M)| < \infty$, it must be finitely generated, which is a contradiction.

	Let us verify condition (2). Consider the subsequence $(k_n)_{n \in \mathbb{N}}$ of naturals satisfying that $\mathsf{n}(r_{k_n}) < \mathsf{n}(r_i)$ for every $i > k_n$. It follows immediately that the sequence $(\mathsf{n}(r_{k_n}))_{n \in \mathbb{N}}$ is increasing. We claim that $M = \langle r_{k_n} : n \in \nn \rangle$. Suppose that $j \notin (k_n)_{n \in \mathbb{N}}$. Because $\lim \mathsf{n}(r_n) = \infty$ there are only finitely many indices $i \in \nn$ such that $\mathsf{n}(r_i) \le \mathsf{n}(r_j)$, and it is easy to see that the maximum of such indices, say $m$, belongs to $(k_n)_{n \in \mathbb{N}}$. As $r_i = p^{\alpha_m - \alpha_i}q^{\beta_i - \beta_m} r_m$, it follows that $r_i \in \langle r_{k_n} : n \in \nn \rangle$. Hence $M = \langle r_{k_n} : n \in \nn \rangle$. Therefore it suffices to show that $r_{k_n} \in \mathcal{A}(M)$ for every $n \in \nn$. If
	\begin{equation} \label{eq:sufficient condition for atomicity of p-adic PM}
		\frac{q^{\beta_{k_n}}}{p^{\alpha_{k_n}}} = \sum_{i=1}^t c_i \frac{q^{\beta_{k_i}}}{p^{\alpha_{k_i}}},
	\end{equation}
	for some $t, c_1, \dots, c_t \in \nn_0$, then $t \ge n$, $c_1 = \dots = c_{n-1} = 0$, and $c_n \in \{0,1\}$. If $c_n = 0$, then by applying the $q$-adic valuation map to both sides of \eqref{eq:sufficient condition for atomicity of p-adic PM} we immediately obtain a contradiction. Thus, $c_n = 1$, which implies that $r_{k_n}$ is an atom. Hence $M$ is atomic.
\end{proof}

%% file: tex/ch4.tex
\chapter{Factorization invariants} \label{chap:factorization invariants}

\section{Sets of Lengths}

Let $M$ be an atomic monoid. Recall from Chapter~\ref{chap:algebraic background} that for each $x \in M$, the set of lengths of $x$ is defined by $\mathsf{L}(x) := \{|z| : z \in \mathsf{Z}(x)\}$. The \emph{system of sets of lengths} of $M$ is defined by
\[
	\mathcal{L}(M) := \{\mathsf{L}(x) : x \in M\}.
\]
In \cite{aG16} the interested reader can find a friendly introduction to sets of lengths and the role they play in factorization theory. In general, sets of lengths and systems of sets of lengths are arithmetic invariants of atomic monoids that have received significant attention in recent years (see, for instance, \cite{ACHP07,FG08,GS16}). 
\smallskip

\subsection{Full Systems of Sets of Lengths for BFMs}

If a monoid $M$ is a BFM, then it is not hard to verify that the set of lengths of each element of $M$ belongs to the collection
\[
	\mathcal{S} = \big\{\{0\}, \{1\}, S : S \subseteq \zz_{\ge 2} \ \text{and} \ |S| < \infty \big\},
\]
i.e., $\mathcal{L}(M) \subseteq \mathcal{S}$. We say that $M$ has \emph{full system of sets of lengths} provided that $\mathcal{L}(M) = \mathcal{S}$. The first class of BFMs with full system of sets of lengths was found by Kainrath \cite{fK99}; he proved that Krull monoids having infinite class groups, with primes in each class, have full systems of sets of lengths. On the other hand, Frisch~\cite{sF13} proved that the subdomain $\text{Int}(\zz)$ of $\qq[x]$ consisting of all integer-valued polynomials also has full system of sets of lengths; this result has been generalized for Dedekind domains~\cite{FNR17}.

In the context of numerical monoids, Geroldinger and Schmid have proved the following realization theorem for sets of lengths. 

\begin{theorem}\footnote{This is a simplified version of the original theorem, where the number of factorizations can be specified for each length.} \cite[Theorem~3.3]{GS18} \label{thm:simplified version of set of length realizability for NS}
	For every nonempty finite subset $S$ of $\zz_{\ge 2}$, there exists a numerical monoid $N$ and $x \in N$ such that $\mathsf{L}(x) = S$.
\end{theorem}

Theorem~\ref{thm:simplified version of set of length realizability for NS} was a crucial tool to construct the first Puiseux monoid with full system of sets of lengths.

\begin{theorem} \cite[Theorem~3.6]{fG19b} \label{thm:PM with full system of sets of lengths}
	There is an atomic Puiseux monoid with full system of sets of lengths.
\end{theorem}

\smallskip
\subsection{Puiseux Monoids with Arithmetic Sets of Lengths}

In this subsection we show that the set of lengths of each element in an atomic multiplicatively cyclic Puiseux monoid $M_r$ is an arithmetic sequence. First, we describe the minimum-length and maximum-length factorizations for elements of $M_r$. We start with the case where $0 < r < 1$.

\begin{lemma} \label{lem:factorization of extremal length II}
	Take $r \in (0,1) \cap \qq$ such that $M_r$ is atomic, and for $x \in M_r^\bullet$ consider the factorization $z = \sum_{i=0}^N \alpha_i r^i \in \mathsf{Z}(x)$, where $N \in \nn$ and $\alpha_0, \dots, \alpha_N \in \nn_0$. The following statements~hold.
	\begin{enumerate}
		\item $\min \mathsf{L}(x) = |z|$ if and only if $\alpha_i < \mathsf{d}(r)$ for $i \in \{ 1, \dots, N \}$.
		
		\item There exists exactly one factorization in $\mathsf{Z}(x)$ of minimum length.
		
		\item $\sup \mathsf{L}(x) = \infty$ if and only if $\alpha_i \ge \mathsf{n}(r)$ for some $i \in \{ 0, \dots, N \}$.
		
		\item $|\mathsf{Z}(x)| = 1$ if and only if $|\mathsf{L}(x)| = 1$, in which case, $\alpha_i < \mathsf{n}(r)$ for $i \in \{ 0, \dots, N \}$.
	\end{enumerate} 
\end{lemma}

\begin{proof}
	To verify the direct implication of~(1), we only need to observe that if $\alpha_i \ge \mathsf{d}(r)$ for some $i \in \{ 1,\dots, N \}$, then the identity $\alpha_i r^i = (\alpha_i - \mathsf{d}(r))r^i + \mathsf{n}(r)r^{i-1}$ would yield a factorization $z'$ in $\mathsf{Z}(x)$ with $|z'| < |z|$. To prove the reverse implication, suppose that $w := \sum_{i=0}^K \beta_i r^i \in \mathsf{Z}(x)$ has minimum length. By the implication already proved, $\beta_i < \mathsf{d}(r)$ for $i \in \{ 1, \dots, N \}$. Insert zero coefficients if necessary and assume that $K = N$. Suppose, by way of contradiction, that there exists $m \in \{ 1, \dots, N \}$ such that $\beta_m \neq \alpha_m$ and assume that such index $m$ is as large as possible. Since $z,w \in \mathsf{Z}(x)$ we can write
	\[
		(\alpha_m - \beta_m)r^m = \sum_{i=0}^{m-1} (\beta_i - \alpha_i) r^i.
	\]
	After multiplying the above equality by $\mathsf{d}(r)^m$, it is easy to see that $\mathsf{d}(r) \mid \alpha_m - \beta_m$, which contradicts the fact that $0 < |\alpha_m - \beta_m| \le \mathsf{d}(r)$. Hence $\beta_i = \alpha_i$ for $i \in \{ 0, \dots, N \}$ and, therefore, $w = z$. As a result, $|z| = |w| = \min \mathsf{L}(x)$. In particular, there exists only one factorization in $\mathsf{Z}(x)$ having minimum length, and~(2) follows.
	
	For the direct implication of~(3), take a factorization $w = \sum_{i=0}^N \beta_i r^i \in \mathsf{Z}(x)$ whose length is not the minimum of $\mathsf{L}(x)$; such a factorization exists because $\sup \mathsf{L}(x) =~\infty$. By part~(1), there exists $i \in \{ 1, \dots, N \}$ such that $\beta_i \ge \mathsf{d}(r)$. Now we can use the identity $\beta_i r^i = (\beta_i - \mathsf{d}(r))r^i + \mathsf{n}(r)r^{i-1}$ to obtain $w_1 \in \mathsf{Z}(x)$ with $|w_1| < |w|$. Notice that there is an atom (namely $r^{i-1}$) appearing at least $\mathsf{n}(r)$ times in $w_1$. In a similar way we can obtain factorizations $w = w_0, w_1, \dots, w_n$ in $\mathsf{Z}(x)$, where $w_n =: \sum_{i=0}^N \beta'_i r^i \in \mathsf{Z}(x)$ satisfies $\beta_i' < \mathsf{d}(r)$ for $i \in \{ 1, \dots, N \}$. By~(1) we have that $w_n$ is a factorization of minimum length and, therefore, $z = w_n$ by~(2). Hence $\alpha_i \ge \mathsf{n}(r)$ for some $i \in \{ 0, \dots, N \}$, as desired. For the reverse implication, it suffices to note that given a factorization $w = \sum_{i=0}^N \beta_i r^i \in \mathsf{Z}(x)$ with $\beta_i \ge \mathsf{n}(r)$ we can use the identity $\beta_i r^i = (\beta_i - \mathsf{n}(r)) r^i + \mathsf{d}(r) r^{i+1}$ to obtain another factorization $w' = \sum_{i=0}^{N+1} \beta'_i r^i \in \mathsf{Z}(x)$ (perhaps $\beta'_{N+1} = 0$) with $|w'| > |w|$ and satisfying $\beta_{i+1} > \mathsf{n}(r)$. 
	
	Finally, we argue the reverse implication of~(4) as the direct implication is trivial. To do this, assume that $\mathsf{L}(x)$ is a singleton. Then each factorization of $x$ has minimum length. By~(2) there exists exactly one factorization of minimum length in~$\mathsf{Z}(x)$. Thus, $\mathsf{Z}(x)$ is also a singleton. The last statement of~(4) is straightforward.
\end{proof}

We continue with the case of $r>1$.

\begin{lemma} \label{lem:factorization of extremal length I}
	Take $r \in \qq_{> 1} \setminus \nn$ such that $M_r$ is atomic, and for $x \in M_r^\bullet$ consider the factorization $z = \sum_{i=0}^N \alpha_i r^i \in \mathsf{Z}(x)$, where $N \in \nn$ and $\alpha_0, \dots, \alpha_N \in \nn_0$. The following statements hold.
	\begin{enumerate}
		\item $\min \mathsf{L}(x) = |z|$ if and only if $\alpha_i < \mathsf{n}(r)$ for $i \in \{ 0, \dots, N \}$.
		
		\item There exists exactly one factorization in $\mathsf{Z}(x)$ of minimum length.
		
		\item $\max \mathsf{L}(x) = |z|$ if and only if $\alpha_i < \mathsf{d}(r)$ for $i \in \{ 1, \dots, N \}$.
		
		\item There exists exactly one factorization in $\mathsf{Z}(x)$ of maximum length.
		
		\item $|\mathsf{Z}(x)| = 1$ if and only if $|\mathsf{L}(x)| = 1$, in which case $\alpha_0 < \mathsf{n}(r)$ and $\alpha_i < \mathsf{d}(r)$ for $i \in \{ 1, \dots, N \}$.
	\end{enumerate} 
\end{lemma}

\begin{proof}
	To argue the direct implication of~(1) it suffices to note that if $\alpha_i \ge \mathsf{n}(r)$ for some $i \in \{ 0, \dots, N \}$, then we can use the identity $\alpha_i r^i = (\alpha_i - \mathsf{n}(r))r^i + \mathsf{d}(r)r^{i+1}$ to obtain a factorization $z'$ in $\mathsf{Z}(x)$ satisfying $|z'| < |z|$. For the reverse implication, suppose that $w = \sum_{i=0}^K \beta_i r^i$ is a factorization in $\mathsf{Z}(x)$ of minimum length. There is no loss in assuming that $K = N$. Note that $\beta_i < \mathsf{n}(r)$ for each $i \in \{ 0, \dots, N \}$ follows from the direct implication. Now suppose for a contradiction that $w \neq z$, and let $m$ be the smallest nonnegative integer satisfying that $\alpha_m \neq \beta_m$. Then
	\begin{equation} \label{eq:set length of rational semirings II}
		(\alpha_m - \beta_m) r^m = \sum_{i = m+1}^N (\beta_i - \alpha_i) r^i.
	\end{equation}
	After clearing the denominators in~(\ref{eq:set length of rational semirings II}), it is easy to see that $\mathsf{n}(r) \mid \alpha_m - \beta_m$, which implies that $\alpha_m = \beta_m$, a contradiction. Hence $w = z$ and so $|z| = |w| = \min \mathsf{L}(x)$. We have also proved that there exists a unique factorization of $x$ of minimum length, which is~(2).
	
	For the direct implication of~(3), it suffices to observe that if $\alpha_i \ge \mathsf{d}(r)$ for some $i \in \{ 1, \dots, N \}$, then we can use the identity $\alpha_i r^i = \big( \alpha_i - \mathsf{d}(r) \big) r^i + \mathsf{n}(r) r^{i-1}$ to obtain a factorization $z'$ in $\mathsf{Z}(x)$ satisfying $|z'| > |z|$. For the reverse implication of~(3), take $w = \sum_{i=0}^K \beta_i r^i$ to be a factorization in $\mathsf{Z}(x)$ of maximum length ($M_r$ is a BFM because $0$ is not a limit point of $M_r^\bullet$). Once again, there is no loss in assuming that $K = N$. The maximality of $|w|$ now implies that $\beta_i < \mathsf{d}(r)$ for $i \in \{ 1, \dots, N \}$. Suppose, by way of contradiction, that $z \neq u$. Then take $m$ be the smallest index such that $\alpha_m \neq \beta_m$. Clearly, $m \ge 1$ and
	\begin{equation*}
		(\alpha_m - \beta_m) r^m = \sum_{i=0}^{m-1} (\beta_i - \alpha_i) r^i.
	\end{equation*}
	After clearing denominators, it is easy to see that $\mathsf{d}(r) \mid \alpha_m - \beta_m$, which contradicts that $0 < |\alpha_M - \beta_M| < \mathsf{d}(r)$. Hence $\alpha_i = \beta_i$ for each $i \in \{ 1, \dots, N \}$, which implies that $z = w$. Thus, $\max \mathsf{L}(x) = |z|$. In particular, there exists only one factorization of $x$ of maximum length, which is condition~(4).
	
	The direct implication of~(5) is trivial. For the reverse implication of~(5), suppose that $\mathsf{L}(x)$ is a singleton. Then any factorization in $\mathsf{Z}(x)$ is a factorization of minimum length. Since we proved in the first paragraph that $\mathsf{Z}(x)$ contains only one factorization of minimum length, we have that $\mathsf{Z}(x)$ is also a singleton. The last statement of~(5) is an immediate consequence of~(1) and~(3).
\end{proof}

We are in a position now to describe the sets of lengths of any atomic multiplicatively cyclic Puiseux monoid.

\begin{theorem} \label{thm:sets of lengths}
	Take $r \in \qq_{>0}$ such that $M_r$ is atomic.
	\begin{enumerate}
		\item If $r < 1$, then for each $x \in M_r$ with $|\mathsf{Z}(x)| > 1$,
		\[
			\mathsf{L}(x) = \big\{ \min \mathsf{L}(x) + k \big( \mathsf{d}(r) - \mathsf{n}(r)\big) : k \in \nn_0 \big\}.
		\]
		
		\item If $r \in \nn$, then $|\mathsf{Z}(x)| = |\mathsf{L}(x)| = 1$ for all $x \in M_r$.

		\item If $r \in \qq_{> 1} \setminus \nn$, then for each $x \in M_r$ with $|\mathsf{Z}(x)| > 1$,
		\[
			\mathsf{L}(x) = \bigg\{ \min \mathsf{L}(x) + k \big( \mathsf{n}(r) - \mathsf{d}(r) \big) : 0 \le k \le \frac{\max \mathsf{L}(x) - \min \mathsf{L}(x)}{\mathsf{n}(r) - \mathsf{d}(r)} \bigg\}.
		\]
	\end{enumerate}
	Thus, $\mathsf{L}(x)$ is an arithmetic progression with difference $|\mathsf{n}(r) - \mathsf{d}(r)|$ for all $x \in M_r$.
\end{theorem}

\begin{proof}
	To argue~(1), take $x \in M_r$ such that $|\mathsf{Z}(x)| > 1$. Let $z := \sum_{i=0}^N \alpha_i r^i$ be a factorization in $\mathsf{Z}(x)$ with $|z| > \min \mathsf{L}(x)$. Lemma~\ref{lem:factorization of extremal length II} guarantees that $\alpha_i \ge \mathsf{d}(r)$ for some $i \in \{ 1, \dots, N \}$. Then one can use the identity $\alpha_i r^i = (\alpha_i - \mathsf{d}(r))r^i + \mathsf{n}(r)r^{i-1}$ to find a factorization $z_1 \in \mathsf{Z}(x)$ with $|z_1| = |z| - (\mathsf{d}(r) - \mathsf{n}(r))$. Carrying out this process as many times as necessary, we can obtain a sequence $z_1, \dots, z_n \in \mathsf{Z}(x)$, where $z_n =: \sum_{i=0}^K \alpha'_i r^i$ satisfies that $\alpha'_i < \mathsf{d}(r)$ for $i \in \{ 1, \dots, K \}$ and $|z_j| = |z| - j(\mathsf{d}(r) - \mathsf{n}(r))$ for $j \in \{ 1, \dots, n \}$. By Lemma~\ref{lem:factorization of extremal length II}(1), the factorization $z_n$ has minimum length and, therefore, $|z| \in \{ \min \mathsf{L}(x) + k \big( \mathsf{d}(r) - \mathsf{n}(r)\big) : k \in \nn_0 \}$. Then
	\begin{align*}
		\mathsf{L}(x) \subseteq \big\{ \min \mathsf{L}(x) + k \big( \mathsf{d}(r) - \mathsf{n}(r)\big) : k \in \nn_0 \big\}.
	\end{align*}
	For the reverse inclusion, we check inductively that $\min \mathsf{L}(x) + k (\mathsf{d}(r) - \mathsf{n}(r)) \in \mathsf{L}(x)$ for every $k \in \nn_0$. Since $|\mathsf{Z}(x)| > 1$, Lemma~\ref{lem:factorization of extremal length II}(2) guarantees that $|\mathsf{L}(x)| > 1$. Then there exists a factorization of length strictly greater than $\min \mathsf{L}(x)$, and we have already seen that such a factorization can be connected to a minimum-length factorization of $\mathsf{Z}(x)$ by a chain of factorizations in $\mathsf{Z}(x)$ with consecutive lengths differing by $\mathsf{d}(r) - \mathsf{n}(r)$. Therefore $\min \mathsf{L}(x) + (\mathsf{d}(r) - \mathsf{n}(r)) \in \mathsf{L}(x)$. Suppose now that $z = \sum_{i=0}^N \beta_i r^i$ is a factorization in $\mathsf{Z}(x)$ with length $\min \mathsf{L}(x) + k(\mathsf{d}(r) - \mathsf{n}(r))$ for some $k \in \nn$. Then by Lemma~\ref{lem:factorization of extremal length II}(1), there exists $i \in \{ 1, \dots, N \}$ such that $\beta_i \ge \mathsf{d}(r) > \mathsf{n}(r)$. Now using the identity $\beta_i r^i = (\beta_i - \mathsf{n}(r))r^i + \mathsf{d}(r)r^{i+1}$, one can produce a factorization $z' \in \mathsf{Z}(x)$ such that $|z'| = \min \mathsf{L}(x) + (k+1)(\mathsf{d}(r) - \mathsf{n}(r))$. Hence the reverse inclusion follows by induction.
	
	Clearly, statement~(2) is a direct consequence of the fact that $r \in \nn$ implies that $M_r = (\nn_0,+)$.
	
	To prove~(3), take $x \in S^\bullet_r$. Since $M_r$ is a BFM, there exists $z \in \mathsf{Z}(x)$ such that $|z| = \max \mathsf{L}(x)$. Take $N \in \nn$ and $\alpha_0, \dots, \alpha_N \in \nn_0$ such that $z = \sum_{i=0}^N \alpha_i r^i$. If $\alpha_i \ge \mathsf{n}(r)$ for some $i \in \{ 0, \dots, N \}$, then we can use the identity $\alpha_i r^i = (\alpha_i - \mathsf{n}(r))r^i + \mathsf{d}(r)r^{i+1}$ to find a factorization $z_1 \in \mathsf{Z}(x)$ such that $|z_1| = |z| - (\mathsf{n}(r) - \mathsf{d}(r))$. Carrying out this process as many times as needed, we will end up with a sequence $z_1, \dots, z_n \in \mathsf{Z}(x)$, where $z_n =: \sum_{i=0}^K \beta_i r^i$ satisfies that $\beta_i < \mathsf{n}(r)$ for $i \in \{ 0, \dots, K \}$ and $|z_j| = |z| - j(\mathsf{n}(r) - \mathsf{d}(r))$ for $j \in \{ 1, \dots, n \}$. Lemma~\ref{lem:factorization of extremal length I}(1) ensures that $|z_n| = \min \mathsf{L}(x)$. Then
	\begin{align} \label{eq:sets of lengths are arithmetic sequences 2}
		\bigg\{ \min \mathsf{L}(x) + j( \mathsf{n}(r) - \mathsf{d}(r) ) : 0 \le j \le \frac{ \max \mathsf{L}(x) - \min \mathsf{L}(x) }{ \mathsf{n}(r) - \mathsf{d}(r) } \bigg\} \subseteq \mathsf{L}(x).
	\end{align}
	On the other hand, we can connect any factorization $w \in \mathsf{Z}(x)$ to the minimum-length factorization $w' \in \mathsf{Z}(x)$ by a chain $w = w_1, \dots, w_t = w'$ of factorizations in $\mathsf{Z}(x)$ so that $|w_i|- |w_{i+1}| = \mathsf{n}(r) - \mathsf{d}(r)$. As a result, both sets involved in the inclusion~(\ref{eq:sets of lengths are arithmetic sequences 2}) are indeed equal.
\end{proof}

We conclude this section collecting some immediate consequences of Theorem~\ref{thm:sets of lengths}.

\begin{corollary} \label{bigcor} Take $r \in \qq_{>0}$ such that $M_r$ is atomic.
	\begin{enumerate}
		\item $M_r$ is a BFM if and only if $r \ge 1$.
		
		\item If $r \in \nn$, then $M_r = \nn_0$ and, as a result, $\Delta(x) = \emptyset$. 
		
		\item If $r \notin \nn$, then $\Delta(x) = \{ |\mathsf{n}(r) - \mathsf{d}(r)| \}$ for all $x \in M_r$ such that $|\mathsf{Z}(x)| > 1$. Therefore $\Delta(M_r) = \{ |\mathsf{n}(r) - \mathsf{d}(r)| \}$.
		
	\end{enumerate}
\end{corollary}


\medskip
\section{Elasticity}
\label{sec:elasticity}

Similar to the system of sets of lengths, the elasticity is another arithmetical invariant used to measure up to what extent factorizations in monoids (or domains) fail to be unique. The elasticity was introduced by R.~Valenza~\cite{rV90} as a tool to measure the phenomenon of non-unique factorizations in the context of algebraic number theory. The elasticity of numerical monoids has been successfully studied in~\cite{CHM06}. In addition, the elasticity of atomic monoids naturally generalizing numerical monoids has received substantial attention in the literature in recent years (see, for instance,~\cite{fG19c, GO19, mG19, qZ18}).

\begin{definition}
	The \emph{elasticity} $\rho(M)$ of an atomic monoid $M$ is given by
	\[
		\rho(M) = \sup \{\rho(x) : x \in M\}, \ \text{where} \ \rho(x) = \frac{\sup \mathsf{L}(x)}{\min \mathsf{L}(x)}.
	\]
\end{definition}

The following formula for the elasticity of an atomic Puiseux monoid in terms of the infimum and supremum of its set of atoms was established in~\cite{GO19}.

\begin{theorem} \cite[Theorem~3.2]{GO19} \label{thm:elasticity of PM}
	Let $M$ be an atomic Puiseux monoid. If $0$ is a limit point of $M^\bullet$, then
	$\rho(M) = \infty$. Otherwise,
	\[
		\rho(M) = \frac{\sup \mathcal{A}(M)}{\inf \mathcal{A}(M)}.
	\]
\end{theorem}
\smallskip

Note that if $\mathcal{A}(M)$ is finite, then $M$ is a numerical monoid; in this case, Theorem~\ref{thm:elasticity of PM} coincides with the elasticity formula given in \cite[Theorem~2.1]{CHM06}.

The elasticity of $M$ is accepted if there exists $x \in M$ with $\rho(x) = \rho(M)$.

\begin{theorem} \cite[Theorem~3.4]{GO19} \label{thm:accepted elasticity of PM}
	For any atomic Puiseux monoid $M$ such that $\rho(M) < \infty$, the elasticity of $M$ is accepted if and only if $\mathcal{A}(M)$ has both a maximum and a minimum.
\end{theorem}

For an atomic monoid $M$ the set
\[
	R(M) = \{ \rho(x) : x \in M\}
\]
is called the \emph{set of elasticities} of $M$, and $M$ is called \emph{fully elastic} if $R(M) = \qq \cap [1, \rho(M)]$ when $\infty \notin R(M)$ and $R(M) \setminus \{\infty\} = \qq \cap [1, \infty)$ when $\infty \in R(M)$. The sets of elasticities of prime reciprocal Puiseux monoids was described in~\cite[Section~4]{GO19}.

\smallskip
\subsection{Union of Sets of Lengths and Local Elasticity} For a nontrivial reduced monoid $M$ and $k \in \nn$, we let $\mathcal{U}_k(M)$ denote the union of sets of lengths containing $k$, that is, $\mathcal{U}_k(M)$ is the set of $\ell \in \nn$ for which there exist atoms $a_1, \dots, a_k, b_1, \dots, b_\ell$ such that $a_1 \dots a_k = b_1 \dots b_\ell$. The set $\mathcal{U}_k(M)$ is known as the \emph{union of sets of lengths} of $M$ containing~$k$. In addition, we set
\[
	\lambda_k(M) := \min \, \mathcal{U}_k(M) \quad \text{and} \quad  \rho_k(M) := \sup \, \mathcal{U}_k(M),
\]
and we call $\rho_k(M)$ the \emph{k-th local elasticity} of $M$. Unions of sets of lengths have received a great deal of attention in recent literature; see, for example, \cite{BS18,BGG11,FGKT17}. By~\cite[Section~1.4]{GH06}, the elasticity of an atomic monoid can be expressed in terms of its local elasticities as follows
\[
	\rho(M) = \sup \bigg\{ \frac{\rho_k(M)}{k} : k \in \nn \bigg\} = \lim_{k \to \infty} \frac{\rho_k (M)}{k}.
\]
For $k \in \mathbb{N}_0$, we define $\mathsf{L}^{-1}(k) := \{x \in M : k \in \mathsf{L}(x)\}$. It is easy to verify that
\[
	\mathcal{U}_k(M) = \{|z| : z \in \mathsf{Z}(x) \ \text{for some} \ x \in \mathsf{L}^{-1}(k) \}.
\]

For a numerical monoid $N$ with minimal generating set $A$, it was proved in \cite[Section~2]{CHM06} that the elasticity of $N$ is given by $\max A/ \min A$. On the other hand, it is not hard to verify that $\mathcal{U}_n(N)$ is bounded and, therefore, every local elasticity of $N$ is finite. In the next two sections, we will generalize this fact in two different ways to Puiseux monoids.

Now we will propose a sufficient condition under which most of the local elasticities of an atomic Puiseux monoid have infinite cardinality. On the other hand, we will describe a subclass of Puiseux monoids (containing isomorphic copies of each numerical monoid) whose local $k$-elasticities are finite.

If $P$ is a Puiseux monoid, then we say that $a_0 \in \mathcal{A}(P)$ is \emph{stable} provided that the set $\{a \in \mathcal{A}(P) : \mathsf{n}(a) = \mathsf{n}(a_0)\}$ is infinite.

\begin{proposition} \label{prop:union of sets of lengths: infinite case}
	Let $P$ be an atomic Puiseux monoid. If $P$ contains a stable atom, then $\rho_k(P)$ is infinite for all sufficiently large $k$.
\end{proposition}

\begin{proof}
	Suppose that for some $m \in \mathbb{N}$ the set $A := \{a \in \mathcal{A}(P) : \mathsf{n}(a) = m\}$ contains infinitely many elements. Let $(a_n)_{n \in \mathbb{N}}$ be an enumeration of the elements of $A$. Because the elements of $A$ have the same numerator, namely $m$, we can assume that the sequence $(a_n)_{n \in \mathbb{N}}$ is decreasing. Setting $d = \mathsf{d}(a_1)$, we can easily see that $d a_1 = m = \mathsf{d}(a_j) a_j$ for each $j \in \mathbb{N}$. Therefore $\mathsf{d}(a_j) \in \mathcal{U}_d(P)$ for each $j \in \mathbb{N}$. As $\mathsf{d}(A)$ is an infinite set so is $\mathcal{U}_d(P)$. The fact that $|\mathcal{U}_d(P)| = \infty$ immediately implies that $|\mathcal{U}_k(P)| = \infty$ for all $k \ge d$. Hence $\rho_k(P) = \sup \, \mathcal{U}_k(P) = \infty$ for every $k \ge d$.
\end{proof}

Recall that a Puiseux monoid $P$ is strongly bounded if it can be generated by a set of rationals $A$ whose numerator set $\mathsf{n}(A)$ is bounded. As a direct consequence of Proposition~\ref{prop:union of sets of lengths: infinite case} we obtain the following result.

\begin{corollary}
	If $P$ is a non-finitely generated strongly bounded atomic Puiseux monoid, then $\rho_k(P)$ is infinite for all $k$ sufficiently large.
\end{corollary}

In contrast to the previous proposition, the next result gives a condition under which Puiseux monoids have finite $k$-elasticity for each $k \in \mathbb{N}$.

\begin{proposition} \label{prop:union of sets of lengths: finite case}
	Let $P$ be a Puiseux monoid that does not contain $0$ as a limit point. If $P$ is bounded, then $\rho_k(P) < \infty$ for every $k \in \mathbb{N}$.
\end{proposition}

\begin{proof}
	Because $0$ is not a limit point of $P^\bullet$, the Puiseux monoid $P$ is atomic As $P$ is a bounded Puiseux monoid, $\mathcal{A}(P)$ is a bounded set of rational numbers. Take $q, Q \in \mathbb{Q}$ such that $0 < q < a < Q$ for all $a \in \mathcal{A}(P)$. Now fix $k \in \mathbb{N}$, and suppose that $\ell \in \mathcal{U}_k(P)$. Then there exists $x \in \mathsf{L}^{-1}(k)$ such that $\ell \in \mathsf{L}(x)$. Because $x$ has a factorization of length $k$, it follows that $x < kQ$. Taking $a_1, \dots, a_\ell \in \mathcal{A}(P)$ such that $x = a_1 + \dots + a_\ell$, we find that
	\[
		q \ell < a_1 + \dots + a_\ell = x < kQ.
	\]
	Therefore $\ell < kQ/q$. Because neither $q$ nor $Q$ depends on the choice of $x$, one obtains that $\mathcal{U}_k(P)$ is bounded from above by $kQ/q$. Hence $\rho_k(P) = \sup \mathcal{U}_k(P)$ is finite, and the proof follows.
\end{proof}

With the following two examples, we shall verify that the conditions of containing a stable atom and not having $0$ as a limit point are not superfluous in Proposition~\ref{prop:union of sets of lengths: infinite case} and Proposition~\ref{prop:union of sets of lengths: finite case}, respectively.

\begin{example}
	Let $(p_n)_{n \in \mathbb{N}}$ be a strictly increasing enumeration of the prime numbers, and consider the following Puiseux monoid:
	\[
		P = \langle A \rangle, \ \text{ where } \ A = \bigg\{ \frac{p_n - 1}{p_n} \ : \ n \in \mathbb{N} \bigg \}.
	\]
	As the denominators of elements in $A$ are pairwise distinct primes, it immediately follows that $\mathcal{A}(P) = A$. Therefore $P$ is atomic. Clearly, $P$ does not contain stable atoms. Because $A$ is bounded so is $P$ (as a Puiseux monoid). On the other hand, $0$ is not a limit point of $P$. Thus, it follows by Proposition~\ref{prop:union of sets of lengths: finite case} that $\rho_k(P)$ is finite for every $k \in \mathbb{N}$.
	Notice also that
	\begin{enumerate}
		\item if $q \in P$ has at least two factorizations with no atoms in common, then $q \in \mathbb{N}$;
		\item by Proposition \ref{prop:union of sets of lengths: finite case}, we have both a lower and an upper bound for any $q \in ~\mathsf{L}^{-1}(k)$.
	\end{enumerate}

	Using the previous two observations, we have created an 
	\href{https://www.github.com/marlycormar/find\_u\_k}{R-script} that generates the sets $\mathcal{U}_k$ for $k \in \{1, \dots, 15\}$. Each $\mathcal{U}_k$ appears as the $k$-th column in Table~\ref{fig:U_k}.
\end{example}

\begin{table}[ht] 
	\centering
	\caption{$\mathcal{U}_k$ for $k \in \{1, \dots, 15\}$.} 
	\label{fig:U_k}
	\begin{tabular}{rlllllllllllllll}
		\hline
		$\mathcal{U}_{1}$ & $\mathcal{U}_{2}$ & $\mathcal{U}_{3}$ & $\mathcal{U}_{4}$ & $\mathcal{U}_{5}$ & $\mathcal{U}_{6}$ & $\mathcal{U}_{7}$ & $\mathcal{U}_{8}$ & $\mathcal{U}_{9}$ & $\mathcal{U}_{10}$ & $\mathcal{U}_{11}$ & $\mathcal{U}_{12}$ & $\mathcal{U}_{13}$ & $\mathcal{U}_{14}$ & $\mathcal{U}_{15} $ \\ 
		\hline
		1 & 2 & 3 & 3 & 4 & 5 & 5 & 5 & 6 & 7 & 7 & 7 & 8 & 9 & 10 \\ 
		&   & 4 & 4 & 5 & 6 & 6 & 6 & 7 & 8 & 8 & 8 & 9 & 10 & 11 \\ 
		&   &   & 5 & 6 & 7 & 7 & 7 & 8 & 9 & 9 & 9 & 10 & 11 & 12 \\ 
		&   &   &   & 7 & 8 & 8 & 8 & 9 & 10 & 10 & 10 & 11 & 12 & 13 \\ 
		&   &   &   & 8 & 9 & 9 & 9 & 10 & 11 & 11 & 11 & 12 & 13 & 14 \\ 
		&   &   &   &   &   & 10 & 10 & 11 & 12 & 12 & 12 & 13 & 14 & 15 \\ 
		&   &   &   &   &   & 11 & 11 & 12 & 13 & 13 & 13 & 14 & 15 & 16 \\ 
		&   &   &   &   &   & 12 & 12 & 13 & 14 & 14 & 14 & 15 & 16 & 17 \\ 
		&   &   &   &   &   &   & 13 & 14 & 15 & 15 & 15 & 16 & 17 & 18 \\ 
		&   &   &   &   &   &   &   &   & 16 & 16 & 16 & 17 & 18 & 19 \\ 
		&   &   &   &   &   &   &   &   &   & 17 & 17 & 18 & 19 & 20 \\ 
		&   &   &   &   &   &   &   &   &   & 18 & 18 & 19 & 20 & 21 \\ 
		&   &   &   &   &   &   &   &   &   & 19 & 19 & 20 & 21 & 22 \\ 
		&   &   &   &   &   &   &   &   &   & 20 & 20 & 21 & 22 & 23 \\ 
		&   &   &   &   &   &   &   &   &   &   & 21 & 22 & 23 & 24 \\ 
		&   &   &   &   &   &   &   &   &   &   &   & 23 & 24 & 25 \\ 
		&   &   &   &   &   &   &   &   &   &   &   & 24 & 25 & 26 \\ 
		\hline
	\end{tabular}
\end{table}

\begin{example}
	Let $(p_n)_{n \in \mathbb{N}}$ be an enumeration of the prime numbers, and consider the Puiseux monoid $P = \big\langle 1/p_n : n \in \mathbb{N} \big\rangle$. It is not difficult to argue that $P$ is atomic with $\mathcal{A}(P) = \{1/p_n : n \in \mathbb{N}\}$. As $\mathcal{A}(P)$ is a bounded subset of positive rationals, the Puiseux monoid $P$ is bounded. Notice, however, that $0$ is a limit point of $P$. By Proposition~\ref{prop:union of sets of lengths: infinite case}, it follows that the local elasticities $\rho_k(P)$ are infinite for all $k$ sufficiently large. 
\end{example}

The condition of boundedness on Proposition~\ref{prop:union of sets of lengths: finite case} is also required, as shown by the following proposition.

\begin{proposition} \label{prop:PM with all its local elasticities infinite}
	There exist infinitely many non-isomorphic Puiseux monoids without $0$ as a limit point that have no finite local elasticities.
\end{proposition}

\begin{proof}
	Let $\mathcal{P} = \{S_n : n \in \mathbb{N} \}$ be a family of disjoint infinite sets of odd prime numbers. For each set $S_n$, we will construct an atomic Puiseux monoid $M_n$. Then we will show that $M_i \cong M_j$ implies $i = j$.

	Fix $j \in \mathbb{N}$ and take $p \in S_j$. To construct the Puiseux monoid $M_j$, let us inductively create a sequence $(A_n)_{n \in \mathbb{N}}$ of finite subsets of positive rationals with $A_1 \subsetneq A_2 \subsetneq \cdots$ such that, for each $k \in \mathbb{N}$, the following three conditions hold:
	\begin{enumerate}
		\item $\mathsf{d}(A_k)$ consists of odd prime numbers;

		\item $\mathsf{d}(\max A_k) = \max \, \mathsf{d}(A_k)$;

		\item $A_k$ minimally generates the Puiseux monoid $P_k = \langle A_k \rangle$.
	\end{enumerate}
	Take $A_1 = \{1/p\}$, with $p$ an odd prime number, and assume we have already constructed the sets $A_1, \dots, A_n$ for some $n \in \mathbb{N}$ satisfying our three conditions. To construct $A_{n+1}$, we take $a = \max A_n$ and let
	\[
		b_1 = \frac{\mathsf{n}(a) \lfloor q/2 \rfloor}{q} \ \text{ and } \ b_2 = \frac{\mathsf{n}(a)\big(q - \lfloor q/2 \rfloor \big)}{q},
	\]
	where $q$ is an odd prime in $S_j$ satisfying $q > \max \mathsf{d}(A_n)$ and $q \nmid \mathsf{n}(a)$. Using the fact that $q \ge 5$ and $\mathsf{d}(a) \ge 3$, one obtains that
	\[
		b_2 > b_1 = \frac{\lfloor q/2 \rfloor}{q} \mathsf{n}(a) > \frac 13 \mathsf{n}(a) \ge a.	
	\]
	Now set $A_{n+1} = A_n \cup \{b_1, b_2\}$. Notice that $b_1 + b_2 = \mathsf{n}(a)$. Clearly, $A_n \subsetneq A_{n+1}$, and condition~(1) is an immediate consequence of our inductive construction. In addition,
	\[
		\mathsf{d}(\max A_{n+1}) = \mathsf{d}(b_2) = q = \max \mathsf{d}(A_{n+1}),
	\]
	which is condition (2). Therefore it suffices to verify that $A_{n+1}$ minimally generates $P_{n+1} = \langle A_{n+1} \rangle$. Because both $b_1$ and $b_2$ are greater than every element in $A_n$, we only need to check that $b_1 \notin P_n$ and $b_2 \notin \langle A_n \cup \{b_1\} \rangle$. Let $d$ be the product of all the elements in $\mathsf{d}(A_n)$. Assuming that $b_1 = a_1 + \dots + a_r$ for some $a_1, \dots, a_r \in A_n$, and multiplying both sides of the same equality by $qd$, we would obtain that $q \mid \mathsf{n}(b_1)$, which contradicts that $q \nmid \mathsf{n}(a)$. Hence $b_1 \notin P_n$. Similarly, one finds that $b_2 \notin P_n$. Suppose, again by contradiction, that $b_2 \in \langle A_n \cup \{b_1\} \rangle$. Then there exist $a'_1 , \dots, a'_s \in A_n$ and $m \in \mathbb{N}$ such that $b_2 = mb_1 + a'_1 + \dots + a'_s$. Notice that $2b_1 = \mathsf{n}(a) (q-1)/q > b_2$, which implies that $m \le 1$. As $b_2 \notin P_n$, it follows that $m=1$. Then we can write
	\begin{align} \label{eq:b_2}
		\frac{\mathsf{n}(a)}q = b_2 - b_1 = \sum_{i=1}^s a'_i. 
	\end{align}
	Once again, we can multiply the extreme parts of the equality~(\ref{eq:b_2}) by $q \, \mathsf{d}(\{a'_1, \dots, a'_s\})$, to obtain that $q \mid \mathsf{n}(a)$, a contradiction. As a result, condition~(3) follows.

	Now set $M_j := \cup_{n \in \mathbb{N}} P_n$. As $P_1 \subsetneq P_2 \subsetneq \dots$, the set $M_j$ is, indeed, a Puiseux monoid. We can easily see that $M_j$ is generated by the set $A := \cup_{n \in \mathbb{N}} A_n$. Let us verify now that $\mathcal{A}(M_j) = A$. It is clear that $\mathcal{A}(M_j) \subseteq A$. To check the reverse inclusion, suppose that $a \in A$ is the sum of atoms $a_1, \dots, a_r \in \mathcal{A}(M_j)$. Take $t \in \mathbb{N}$ such that $a, a_1, \dots, a_r \in A_t$. Because $A_t$ minimally generates $P_t$ it follows that $r=1$ and $a = a_1$ and, therefore, that $a \in \mathcal{A}(M_j)$. Hence $\mathcal{A}(M_j) = A$, which implies that $M_j$ is an atomic monoid.

	To disregard $0$ as a limit point of $M_j$, it is enough to observe that $\min \mathcal{A}(M_j) = 1/p$. We need to show then that $\rho_k(M_j) = \infty$ for $k \ge 2$. Set $a_n = \max A_n$. When constructing the sequence $(A_n)_{n \in \mathbb{N}}$, we observed that $\mathsf{n}(a_n) = b_{n_1} + b_{n_2}$, where $\{b_{n_1}, b_{n_2}\} = A_{n+1} \setminus A_n$. Because $\mathsf{n}(a_n) \in M_j$ and
	\[
		b_{n_1} + b_{n_2} =\mathsf{n}(a_n) = \mathsf{d}(a_n) a_n,
	\]
	one has that the factorizations $z = b_{n_1} + b_{n_2}$ and $z' = \mathsf{d}(a_n) a_n$ are both in $\mathsf{Z}(\mathsf{n}(a_n))$. Since $|z| = 2$ and $|z'| = \mathsf{d}(a_n)$ it follows that $\mathsf{d}(a_n) \in \mathcal{U}_2(M_j)$. By condition~(2) above, $\mathsf{d}(a_n) =\mathsf{d}(\max A_n) = \max \mathsf{d}(A_n)$. This implies that the set $\{\mathsf{d}(a_n) : n \in \mathbb{N} \}$ contains infinitely many elements. As $\{\mathsf{d}(a_n) : n \in \mathbb{N} \} \subseteq \mathcal{U}_2(M_j)$, we obtain that $\rho_2(M_j) = \infty$. Hence $\rho_k(M_j) = \infty$ for all $k \ge 2$.

	We have just constructed an infinite family $\mathcal{F} := \{M_n : n \in \mathbb{N}\}$ of atomic Puiseux monoids with infinite $k$-elasticities. Let us show now that the monoids in $\mathcal{F}$ are pairwise non-isomorphic. To do this we use the fact that the only homomorphisms between Puiseux monoids are given by rational multiplication \cite[Lemma~3.3]{GGP16}. Take $i,j \in \mathbb{N}$ such that $M_i \cong M_j$. Then there exists $r \in \mathbb{Q}$ such that $M_i = rM_j$. Let $m \in M_j$ such that $\mathsf{d}(m) = p$ and $p \nmid \mathsf{n}(r)$ for some prime $p$ in $S_j$. Since the element $rm \in M_i$ and $p \mid \mathsf{d}(rm)$, we must have that the prime $p$ belongs to $S_i$. Because the sets in $\mathcal{P}$ are pairwise disjoint, we conclude that $i = j$. This completes the proof.
\end{proof}

Proposition~\ref{prop:union of sets of lengths: infinite case} (respectively, Proposition~\ref{prop:union of sets of lengths: finite case}) establishes sufficient conditions under which a Puiseux monoid has most of its local elasticities infinite (respectively, finite). In addition, we have verified that such conditions are not necessary. For the sake of completeness, we now exhibit a Puiseux monoid that does not satisfy the conditions of either of the propositions above and has no finite $k$-elasticity for any $k \ge 2$.

\begin{example}
	Consider the Puiseux monoid
	\[
		P = \left\langle \left(\frac{2}{3} \right)^n : \, n \in \mathbb{N} \right \rangle.
	\]
	It was proved in \cite[Theorem~6.2]{GG17} that $P$ is atomic and $\mathcal{A}(P) = \{(2/3)^n : n \in \mathbb{N}\}$. In addition, it is clear that $P$ is bounded, has $0$ as a limit point, and does not contain any stable atoms. So neither Proposition~\ref{prop:union of sets of lengths: infinite case} nor Proposition~\ref{prop:union of sets of lengths: finite case} applies to $P$. Now we argue that $\rho_k(P) = \infty$ for each $k \in \mathbb{N}$ such that $k \ge 2$.

	Take $k \ge 2$ and set $x = k\frac{2}{3} \in P$. Notice that, by definition, $x \in \mathsf{L}^{-1}(k)$. We can conveniently rewrite $x$ as
	\[
		x = \big((k - 2) + 2\big) \frac{2}{3} = (k - 2)\frac{2}{3} + 3\cdot  \left(\frac{2}{3}\right)^2\! \!,
	\]
	which reveals that $z = (k-2)\frac 23 + 3(\frac 23)^2$ is a factorization of $x$ with $|z| = k+1$. Taking $k' = 3$ to play the role of $k$ and repeating this process as many times as needed, one can obtain factorizations of $x$ of lengths as large as one desires. The fact that $k$ was chosen arbitrarily implies now that $\rho_k(P) = \infty$ for each $k \ge 2$. \\
\end{example}

\smallskip
\subsection{Prime Reciprocal Puiseux Monoids}
\label{sec:primary case}

We proceed to study the local elasticity of prime reciprocal Puiseux monoids. Recall from Section~\ref{sec:ACCP} that a Puiseux monoid is said to be prime reciprocal if it can be generated by a subset of positive rational numbers whose denominators are pairwise distinct primes. We have seen before that every prime reciprocal Puiseux monoid is atomic. 

In Proposition~\ref{prop:union of sets of lengths: finite case}, we established a sufficient condition on Puiseux monoids to ensure that all their local $k$-elasticities are finite. Here we restrict our study to the case of prime reciprocal Puiseux monoids, providing two more sufficient conditions to guarantee the finiteness of all the local $k$-elasticities.

\begin{theorem}\label{theo:sufficient conditions for finite elasticity in ppm}
	For a prime reciprocal Puiseux monoid $P$, the following two conditions hold.
	\begin{enumerate}
		\item If $0$ is not a limit point of $P$, then $\rho_k(P) < \infty$ for every $k \in \mathbb{N}$. \vspace{3pt}
		\item If $P$ is bounded and has no stable atoms, then $\rho_k(P) < \infty$ for every $k \in \mathbb{N}$.
	\end{enumerate}
\end{theorem}

\begin{proof}
	Because every finitely generated Puiseux monoid is isomorphic to a numerical monoid, and numerical monoids have finite $k$-elasticities, we can assume, without loss of generality, that $P$ is not finitely generated.

	To prove condition~(1), suppose, by way of contradiction, that $\rho_k(P) = \infty$ for some $k \in \mathbb{N}$. Because $0$ is not a limit point of $P$ there exists $q \in \mathbb{Q}$ such that $0 < q < a$ for each $a \in \mathcal{A}(P)$. Let
	\[
		\ell = \min \{n \in \mathbb{N} : |\mathcal{U}_n(P)| = \infty\}.
	\]
	Clearly, $\ell \ge 2$. Let $m = \max \, \mathcal{U}_{\ell - 1}(P)$. Now take $N \in \mathbb{N}$ sufficiently large such that, for each $a \in \mathcal{A}(P)$, $a > N$ implies that $\mathsf{d}(a) > \ell$. As $\mathcal{U}_\ell(P)$ contains infinitely many elements, there exists $k \in \mathcal{U}_\ell(P)$ such that
	\[
		k > \max\bigg\{\frac{\ell}{q}N, \, m + 1 \bigg\}.
	\]
	In particular, $k-1$ is a strict upper bound for $\mathcal{U}_{\ell - 1}(P)$. As $k \in \mathcal{U}_\ell(P)$, we can choose an element $x \in P$ such that $\{k,\ell\} \subseteq \mathsf{L}(x)$. Take $A = \{a_1, \dots, a_k\} \subsetneq \mathcal{A}(P)$ and $B = \{b_1, \dots, b_\ell\} \subsetneq \mathcal{A}(P)$ with
	\begin{align} \label{eq:different length factorizations 1}
		a_1 + \dots + a_k = x = b_1 + \dots + b_\ell.
	\end{align}
	Observe that the sets $A$ and $B$ must be disjoint, for if $a \in A \cap B$, canceling $a$ in (\ref{eq:different length factorizations 1}) would yield that $\{\ell - 1, k - 1\} \subseteq \mathsf{L}(x - a)$, which contradicts that $k-1$ is a strict upper bound for $\mathcal{U}_{\ell - 1}(P)$. Because $k > (\ell/q)N$, it follows that
	\[
		x > kq > \ell N.
	\]
	Therefore $b := \max\{b_1, \dots, b_\ell\} > N$, which implies that $p = \mathsf{d}(b) > \ell$. Since $a_i \neq b$ for each $i = 1, \dots, k$, it follows that $p \notin \mathsf{d}(\{a_1, \dots, a_k\})$. We can assume, without loss of generality, that there exists $j \in \{1, \dots, \ell\}$ such that $b_i \neq b$ for every $i \le j$ and $b_{j+1} = \dots = b_\ell = b$. This allows us to rewrite (\ref{eq:different length factorizations 1}) as
	\begin{align} \label{eq:different length factorization 2}
		(\ell - j)b = \sum_{i=1}^k a_i - \sum_{i=1}^j b_i.
	\end{align}
	After multiplying \ref{eq:different length factorization 2} by $p$ times the product $d$ of all the denominators of the atoms $\{a_1, \dots, a_k, b_1, \dots, b_j\}$, we find that $p$ divides $d(\ell - j)b$. As $\gcd(p,d) = 1$ and $\ell - j < p$, it follows that $p$ divides $\mathsf{n}(b)$, which is a contradiction. Hence we conclude that $\rho_k(P) < \infty$ for every $k \in \mathbb{N}$. \vspace{4pt}

	Now we argue the second condition. Let $(a_n)_{n \in \mathbb{N}}$ be an enumeration of the elements of $\mathcal{A}(P)$ such that $(\mathsf{d}(a_n))_{n \in \mathbb{N}}$ is an increasing sequence. Set $p_n = \mathsf{d}(a_n)$. Since $P$ has no stable atoms, $\lim \mathsf{n}(a_n) = \infty$. Let $B$ be an upper bound for $\mathcal{A}(P)$.

	Suppose, by way of contradiction, that $\rho_n(P) = \infty$ for some $n \in \mathbb{N}$. Let $k$ be the smallest natural number such that $|\mathcal{U}_k(P)| = \infty$. Now take $\ell \in \mathcal{U}_k(P)$ large enough such that $\ell - 1 > \max \, \mathcal{U}_{k-1}(P)$ and for each $a \in \mathcal{A}(P)$ satisfying $a \le Bk/\ell$ we have that $\mathsf{n}(a) > Bk$. Take $x \in \mathsf{L}^{-1}(k)$ such that $a_1 + \dots + a_k = x = b_1 + \dots + b_\ell$ for some $a_1, \dots, a_k, b_1, \dots, b_\ell \in \mathcal{A}(P)$. Now set $b = \min\{b_1, \dots, b_\ell\}$. Then
	\[
		b \le \frac{b_1 + \dots + b_\ell}{\ell} = \frac{a_1 + \dots + a_k}{\ell} \le \frac{Bk}{\ell}.
	\]
	Therefore $\mathsf{n}(b) > Bk$. We claim that $\mathsf{d}(b) \notin \mathsf{d}(\{a_1, \dots, a_k\})$. Suppose by contradiction that this is not the case. Then $b = a_i$ for some $i \in \{1, \dots, k\}$. This implies that $\{k - 1, \ell - 1 \} \subseteq \mathsf{L}(x-b)$, contradicting that $\ell - 1 > \max \, \mathcal{U}_{k-1}(P)$. Hence $\mathsf{d}(b) \notin \mathsf{d}(\{a_1, \dots, a_k\})$. Now assume, without loss of generality, that there exists $j \in \{1, \dots, \ell\}$ such that $b_i \neq b$ for each $i \le j$ and $b_{j+1} = \dots = b_\ell = b$. Write
	\begin{align} \label{eq:different length factorization 3}
		(\ell - j) b = \sum_{i=1}^k a_i - \sum_{i=1}^j b_i.
	\end{align}
	From (\ref{eq:different length factorization 3}) we obtain that $p_\ell$ divides $\ell - j$. As a consequence,
	\[
		Bk \ge \sum_{i=1}^k a_i \ge \frac{\ell - j}{p_\ell} \mathsf{n}(b) \ge \mathsf{n}(b) > Bk,
	\]
	which is a contradiction. Hence $\rho_k(P) < \infty$ for every $k \in \mathbb{N}$.
\end{proof}

The sufficient conditions in part (1) of Theorem~\ref{theo:sufficient conditions for finite elasticity in ppm} and the condition of boundedness in part (2) of Theorem~\ref{theo:sufficient conditions for finite elasticity in ppm} are not necessary, as the following example illustrates.

\begin{example} \hfill
	\begin{enumerate}
		\item Consider the prime reciprocal Puiseux monoid
		\[
			P = \left \langle \frac{n}{p_n} : n \in \mathbb{N} \right \rangle,
		\]
		where $(p_n)_{n \in \mathbb{N}}$ is the increasing sequence of all prime numbers. Since $\mathcal{A}(P) = \{n/p_n : n \in \mathbb{N}\}$, it follows that $P$ does not contain any stable atom. It is well known that the sequence $(n/p_n)_{n \in \mathbb{N}}$ converges to $0$, which implies that $P$ is bounded. Hence part~(2) of Theorem~\ref{theo:sufficient conditions for finite elasticity in ppm} ensures that $\rho_k(P) < \infty$ for all $k \in \mathbb{N}$. Thus, the reverse implication of part~(1) in Theorem~\ref{theo:sufficient conditions for finite elasticity in ppm} does not hold. \vspace{4pt}

		\item Consider now the Puiseux monoid
		\[
			P = \left \langle \frac{p_n^2 - 1}{p_n} : n \in \mathbb{N} \right \rangle,
		\]
		where $(p_n)_{n \in \mathbb{N}}$ is any enumeration of the prime numbers. Since $0$ is not a limit point of $P$, we can apply part~(1) of Theorem~\ref{theo:sufficient conditions for finite elasticity in ppm} to conclude that $\rho_k(P) < \infty$ for all $k \in \mathbb{N}$. Notice, however, that $P$ is not bounded. Therefore, the boundedness in part~(2) of Theorem \ref{theo:sufficient conditions for finite elasticity in ppm} is not a necessary condition.
	\end{enumerate}
\end{example}
\smallskip

\subsection{Multiplicatively Cyclic Puiseux Monoids}

On this subsection, we focus on the elasticity of multiplicatively cyclic Puiseux monoids.

\begin{proposition} \label{prop:elasticity of rational semirings}
	Take $r \in \qq_{>0}$ such that $M_r$ is atomic. Then the following statements are equivalent.
	\begin{enumerate}
		\item $r \in \nn$.
		
		\item $\rho(M_r) = 1$.
		
		\item $\rho(M_r)<\infty$.
	\end{enumerate}
	Hence, if $M_r$ is atomic, then either $\rho(M_r)=1$ or $\rho(M_r)=\infty$.
\end{proposition}

\begin{proof}
	To prove that (1) implies (2), suppose that $r \in \nn$. In this case, $M_r \cong \nn_0$. Since~$\nn_0$ is a factorial monoid, $\rho(M_r) = \rho(\nn_0) = 1$.  Clearly, (2) implies (3). Now assume (3) and that $r \notin \nn$. If $r < 1$, then $0$ is a limit point of $M_r^\bullet$ as $\lim_{n \to \infty} r^n = 0$. Therefore it follows by Theorem~\ref{thm:elasticity of PM} that $\rho(M_r) = \infty$. If $r > 1$, then $\lim_{n \to \infty} r^n = \infty$ and, as a result, $\sup \mathcal{A}(M_r) = \infty$. Then Theorem~\ref{thm:elasticity of PM} ensures that $\rho(M_r) = \infty$. Thus, (3) implies~(1). The final statement now easily follows.
\end{proof}

Recall that the elasticity of an atomic monoid $M$ is said to be accepted if there exists $x \in M$ such that $\rho(M) = \rho(x)$.

\begin{proposition}\label{accepted}
	Take $r \in \qq_{> 0}$ such that $M_r$ is atomic. Then the elasticity of $M_r$ is accepted if and only if $r \in \nn$ or $r < 1$.
\end{proposition}

\begin{proof}
	For the direct implication, suppose that $r \in \qq_{> 1} \setminus \nn$. Proposition~\ref{prop:elasticity of rational semirings} ensures that $\rho(M_r) = \infty$. However, as $0$ is not a limit point of $M_r^\bullet$, it follows from Theorem~\ref{thm:BF sufficient condition} that~$M_r$ is a BFM, and, therefore, $\rho(x) < \infty$ for all $x \in M_r$. As a result, $M_r$ cannot have accepted elasticity
	
	For the reverse implication, assume first that $r \in \nn$ and, therefore, that $M_r = \nn_0$. In this case, $M_r$ is a factorial monoid and, as a result, $\rho(M_r) = \rho(1) = 1$. Now suppose that $r < 1$. Then it follows by Proposition~\ref{prop:elasticity of rational semirings} that $\rho(M_r) = \infty$. In addition, for $x = \mathsf{n}(r) \in M_r$ Lemma~\ref{lem:factorization of extremal length II}(1) and Theorem~\ref{thm:sets of lengths}(1) guarantee that
	\[
		\mathsf{L}(x) = \big\{ \mathsf{n}(r) + k \big( \mathsf{d}(r) - \mathsf{n}(r) \big) : k \in \nn_0 \big\}.
	\]
	Because $\mathsf{L}(x)$ is an infinite set, we have $\rho(M_r) = \infty = \rho(x)$. Hence $M_r$ has accepted elasticity, which completes the proof.
\end{proof}

Let us proceed to describe the sets of elasticities of atomic multiplicatively cyclic Puiseux monoids.

\begin{proposition} \label{prop:set of elasticities}
	Take $r \in \qq_{> 0}$ such that $M_r$ is atomic.
	\begin{enumerate}
		\item If $r < 1$, then $R(M_r) = \{1, \infty\}$ and, therefore, $M_r$ is not fully elastic.
		
		\item If $r \in \nn$, then $R(M_r) = \{1\}$ and, therefore, $M_r$ is fully elastic.
		
		\item If $r \in \qq_{> 0} \setminus \nn$ and $\mathsf{n}(r) = \mathsf{d}(r) + 1$, then $M_r$ is fully elastic, in which case $R(M_r) = \qq_{\ge 1}$.
	\end{enumerate}
\end{proposition}

\begin{proof}
	First, suppose that $r < 1$. Take $x \in M_r$ such that $|\mathsf{Z}(x)| > 1$. It follows by Theorem~\ref{thm:sets of lengths}(1) that $\mathsf{L}(x)$ is an infinite set, which implies that $\rho(x) = \infty$. As a result, $\rho(M_r) = \{1,\infty\}$ and then $M_r$ is not fully elastic.
	
	To argue~(2), it suffices to observe that $r \in \nn$ implies that $M_r = (\nn_0,+)$ is a factorial monoid and, therefore, $\rho(M_r) = \{1\}$.
	
	Finally, let us argue that $M_r$ is fully elastic when $\mathsf{n}(r) = \mathsf{d}(r) + 1$. To do so, fix $q \in \qq_{>1}$. Take $m \in \nn$ such that $m \mathsf{d}(q) > \mathsf{d}(r)$, and set $k = m \big( \mathsf{n}(q) - \mathsf{d}(q) \big)$. Let $t = m \mathsf{d}(q) - \mathsf{d}(r)$, and consider the factorizations $z = \mathsf{d}(r) r^k + \sum_{i=1}^{t} r^{k+i} \in \mathsf{Z}(M_r)$ and $z' = \mathsf{d}(r) \cdot 1 + \sum_{i=0}^{k-1} r^i  + \sum_{i=1}^{t} r^{k+i} \in \mathsf{Z}(M_r)$. Since $\mathsf{n}(r) = \mathsf{d}(r) + 1$, it can be easily checked that $\frac{1}{r-1} = \mathsf{d}(r)$. As
	\[
		\mathsf{d}(r) + \sum_{i=0}^{k-1} r^i + \sum_{i=1}^{t} r^{k+i}= \mathsf{d}(r) + \frac{r^k - 1}{r-1} + \sum_{i=1}^{t} r^{k+i} = \mathsf{d}(r)r^k + \sum_{i=1}^{t} r^{k+i},
	\]
	there exists $x \in M_r$ such that $z,z' \in \mathsf{Z}(x)$. By Lemma~\ref{lem:factorization of extremal length I} it follows that $z$ is a factorization of $x$ of minimum length and $z'$ is a factorization of $x$ of maximum length. Thus, 
	\[
		\rho(x) = \frac{|z'|}{|z|} = \frac{\mathsf{d}(r) + k + t}{\mathsf{d}(r) + t} = \frac{m \, \mathsf{n}(q)}{m \, \mathsf{d}(q)} = q.
	\]
	As $q$ was arbitrarily taken in $\qq_{>1}$, it follows that $R(M_r) = \qq_{\ge 1}$. Hence $M_r$ is fully elastic when $\mathsf{n}(r) = \mathsf{d}(r) + 1$.
\end{proof}
\medskip

We were unable to determine in Proposition~\ref{prop:set of elasticities} whether $M_r$ is fully elastic when $r \in \qq_{> 1} \setminus \nn$ with $\mathsf{n}(r) \neq \mathsf{d}(r) + 1$. However, we proved in Proposition \ref{prop:the set of elasticity of S_r is dense when r > 1} that the set of elasticities of $M_r$ is dense in $\rr_{\ge 1}$.

\begin{proposition} \label{prop:the set of elasticity of S_r is dense when r > 1}
	If $r \in \qq_{>1} \setminus \nn$, then the set $R(M_r)$ is dense in $\rr_{\ge 1}$.
\end{proposition}

\begin{proof}
	Since $\sup \mathcal{A}(M_r) = \infty$, it follows by Theorem~\ref{thm:elasticity of PM} that $\rho(M_r) = \infty$. This, along with the fact that $M_r$ is a BFM (because of Theorem~\ref{thm:BF sufficient condition}, ensures the existence of a sequence $(x_n)_{n \in \mathbb{N}}$ of elements of $M_r$ such that $\lim_{n \to \infty} \rho(x_n) = \infty$. Then it follows by~\cite[Lemma~5.6]{GO19} that the set
	\[
		S := \bigg\{ \frac{\mathsf{n}(\rho(x_n)) + k}{\mathsf{d}(\rho(x_n)) + k} : n,k \in \nn \bigg\}
	\]
	is dense in $\rr_{\ge 1}$. Fix $n,k \in \nn$. Take $m \in \nn$ such that $r^m$ is the largest atom dividing $x_n$ in $M_r$. Now take $K := k \gcd(\min \mathsf{L}(x_n), \max \mathsf{L}(x_n))$. Consider the element $y_{n,k} := x_n + \sum_{i=1}^K r^{m + i} \in M_r$. It follows by Lemma~\ref{lem:factorization of extremal length I} that $x_n$ has a unique minimum-length factorization and a unique maximum-length factorization; let them be $z_0$ and $z_1$, respectively. Now consider the factorizations $w_0 := z_0 + \sum_{i=1}^K r^{m + i} \in \mathsf{Z}(y_{n,k})$ and $w_1 := z_1 + \sum_{i=1}^K r^{m + i} \in \mathsf{Z}(y_{n,k})$. Once again, we can appeal to Lemma~\ref{lem:factorization of extremal length I} to ensure that $w_0$ and $w_1$ are the minimum-length and maximum-length factorizations of $y_{n,k}$. Therefore $\min \mathsf{L}(y_{n,k}) = \min \mathsf{L}(x_n) + K$ and $\max \mathsf{L}(y_{n,k}) = \max \mathsf{L}(x_n) + K$. Then we have
	\[
		\rho(y_{n,k}) = \frac{\max \mathsf{L}(y_{n,k})}{\min \mathsf{L}(y_{n,k})} = \frac{\max \mathsf{L}(x_n) + K}{\min \mathsf{L}(x_n) + K} = \frac{\mathsf{n}(\rho(x_n)) + k}{\mathsf{d}(\rho(x_n)) + k}.
	\]
	Since $n$ and $k$ were arbitrarily taken, it follows that $S$ is contained in $R(M_r)$. As $S$ is dense in $\rr_{\ge 1}$ so is $R(M_r)$, which concludes our proof.
\end{proof}

\begin{corollary}
	The set of elasticities of $M_r$ is dense in $\rr_{\ge 1}$ if and only if $r \in \qq_{> 1} \setminus \nn$.
\end{corollary}

\begin{remark}
	Proposition~\ref{prop:the set of elasticity of S_r is dense when r > 1} contrasts with the fact that the elasticity of a numerical monoid is always nowhere dense in $\rr$~\cite[Corollary~2.3]{CHM06}.
\end{remark}
\smallskip

Let us conclude this section studying the unions of sets of lengths and the local elasticities of atomic multiplicatively cyclic Puiseux monoids.

\begin{proposition}\label{local}
	Take $r \in \qq_{> 0}$ such that $M_r$ is atomic. Then $\mathcal{U}_k(M_r)$ is an arithmetic progression containing $k$ with distance $|\mathsf{n}(r) - \mathsf{d}(r)|$ for every $k \in \nn$. More specifically, the following statements hold.
	\begin{enumerate}
		\item If $r < 1$, then
		\begin{itemize}
			\item $\mathcal{U}_k(M_r) = \{k\}$ if $k < \mathsf{n}(r)$,
			
			\item $\mathcal{U}_k(M_r) = \{ k + j (\mathsf{d}(r) - \mathsf{n}(r)) : j \in \nn_0 \}$ if $\mathsf{n}(r) \le k < \mathsf{d}(r)$, and
			
			\item $\mathcal{U}_k(M_r) = \{ k + j (\mathsf{d}(r) - \mathsf{n}(r)) : j \in \zz_{\ge \ell} \}$ for some $\ell \in \zz_{< 0}$ if $k \ge \mathsf{d}(r)$.
		\end{itemize}
		\vspace{3pt}
		
		\item If $r \in \qq_{> 1} \setminus \nn$, then
		\begin{itemize}
			\item $\mathcal{U}_k(M_r) = \{k\}$ if $k < \mathsf{d}(r)$,
			
			\item $\mathcal{U}_k(M_r) = \{ k + j (\mathsf{n}(r) - \mathsf{d}(r)) : j \in \nn_0 \}$ if $\mathsf{d}(r) \le k < \mathsf{n}(r)$, and
			
			\item $\mathcal{U}_k(M_r) = \{ k + j (\mathsf{n}(r) - \mathsf{d}(r)) : j \in \zz_{\ge \ell} \}$ for some $\ell \in \zz_{< 0}$ if $k \ge \mathsf{n}(r)$.
		\end{itemize}
		
		\item If $r \in \nn$, then $\mathcal{U}_k(M_r) = \{k\}$ for every $k \in \nn$.
	\end{enumerate}
\end{proposition}

\begin{proof}
	That $\mathcal{U}_k(M_r)$ is an arithmetic progression containing $k$ with distance $| \mathsf{n}(r) - \mathsf{d}(r)|$ is an immediate consequence of Theorem~\ref{thm:sets of lengths}.
	
	To show~(1), assume that $r < 1$. Suppose first that $k < \mathsf{n}(r)$. Take $L \in \mathcal{L}(M_r)$ with $k \in L$, and take $x \in M_r$ such that $L = \mathsf{L}(x)$. Choose $z = \sum_{i=0}^N \alpha_i r^i \in \mathsf{Z}(x)$ with $\sum_{i=0}^N \alpha_i = k$. Since $\alpha_i \le k < \mathsf{n}(r)$ for $i \in \{ 0, \dots, N \}$, Lemma~\ref{lem:factorization of extremal length II} ensures that $|\mathsf{Z}(x)| = 1$, which yields $L = \mathsf{L}(x) = \{k\}$. Thus, $\mathcal{U}_k(M_r) = \{k\}$. Now suppose that $\mathsf{n}(r) \le k < \mathsf{d}(r)$. Notice that the element $k \in M_r$ has a factorization of length $k$, namely, $k \cdot 1 \in \mathsf{Z}(k)$. Now we can use Lemma~\ref{lem:factorization of extremal length II}(3) to conclude that $\sup \mathsf{L}(k) = \infty$. Hence $\rho_k(M_r) = \infty$. On the other hand, let $x$ be an element of $M_r$ having a factorization of length $k$. Since $k < \mathsf{d}(r)$, it follows by Lemma~\ref{lem:factorization of extremal length II}(1) that any length-$k$ factorization in $\mathsf{Z}(x)$ is a factorization of~$x$ of minimum length. Hence $\lambda_k(M_r) = k$ and, therefore,
	\[
		\mathcal{U}_k(M_r) = \{ k + j (\mathsf{d}(r) - \mathsf{n}(r)) : j \in \nn_0 \}.
	\]
	Now assume that $k \ge \mathsf{d}(r)$. As $k \ge \mathsf{n}(r)$, we have once again that $\rho_k(M_r) = \infty$. Also, because $k \ge \mathsf{d}(r)$ one finds that $(k - \mathsf{d}(r))r + \mathsf{n}(r) \cdot 1$ is a factorization in $\mathsf{Z}(kr)$ of length $k - (\mathsf{d}(r) - \mathsf{n}(r))$. Then there exists $\ell \in \zz_{< 0}$ such that
	\[
		\mathcal{U}_k(M_r) = \{ k + j (\mathsf{d}(r) - \mathsf{n}(r)) : j \in \zz_{\ge \ell} \}.
	\]
	
	Suppose now that $r \in \qq_{> 1} \setminus \nn$. Assume first that $k < \mathsf{d}(r)$. Take $L \in \mathcal{L}(M_r)$ containing $k$ and $x \in M_r$ such that $L = \mathsf{L}(x)$. If $z = \sum_{i=0}^N \alpha_i r^i \in \mathsf{Z}(x)$ satisfies $|z| = k$, then $\alpha_i \le k < \mathsf{d}(r)$ for $i \in \{ 0, \dots, N \}$, and Lemma~\ref{lem:factorization of extremal length I} implies that $L = \mathsf{L}(x) = \{k\}$. As a result, $\mathcal{U}_k(M_r) = \{k\}$. Suppose now that $\mathsf{d}(r) \le k < \mathsf{n}(r)$. In this case, for each $n > k$, we can consider the element $x_n = k r^n \in M_r$ and set $L_n := \mathsf{L}(x_n)$. It is not hard to check that
	\[
		z_n := \mathsf{n}(r) \cdot 1 + \bigg( \sum_{i=1}^{n-1} \big( \mathsf{n}(r) - \mathsf{d}(r)\big) r^i \bigg) + \big( k - \mathsf{d}(r) \big) r^n
	\]
	is a factorization of $x_n$. Therefore $|z_n| = k + n( \mathsf{n}(r) - \mathsf{d}(r)) \in L_n$. Since $k \in L_n$ for every $n \in \nn$, it follows that $\rho_k(M_r) = \infty$. On the other hand, it follows by Lemma~\ref{lem:factorization of extremal length I}(1) that any factorization of length $k$ of an element $x \in M_r$ must be a factorization of minimum length in $\mathsf{Z}(x)$. Hence $\lambda_k(M_r) = k$, which implies that
	\[
		\mathcal{U}_k(M_r) = \{ k + j (\mathsf{n}(r) - \mathsf{d}(r)) : j \in \nn_0 \}.
	\]
	Assume now that $k \ge \mathsf{n}(r)$. As $k \ge \mathsf{d}(r)$ we still obtain $\rho_k(M_r) = \infty$. In addition, because $k \ge \mathsf{n}(r)$, we have that $(k - \mathsf{n}(r)) \cdot 1 + \mathsf{d}(r)r$ is a factorization in $\mathsf{Z}(k)$ having length $k - (\mathsf{n}(r) - \mathsf{d}(r))$. Thus, there exists $\ell \in \zz_{< 0}$ such that
	\[
		\mathcal{U}_k(M_r) = \{ k + j (\mathsf{n}(r) - \mathsf{d}(r)) : j \in \zz_{\ge \ell} \}.
	\]
	
	Finally, condition~(3) follows directly from the fact that $M_r = (\nn_0,+)$ when $r \in \nn$ and, therefore, for every $k \in \nn$ there exists exactly one element in $M_r$ having a length-$k$ factorization, namely $k$.
\end{proof}

\begin{corollary} \label{cor:local elastiicity}
	Take $r \in \qq_{>0}$ such that $M_r$ is atomic. Then $\rho(M_r) < \infty$ if and only if $\rho_k(M_r) < \infty$ for every $k \in \nn$.
\end{corollary}

\begin{proof}
	It follows from \cite[Proposition~1.4.2(1)]{GH06} that $\rho_k(M_r) \le k \rho(M_r)$, which yields the direct implication. For the reverse implication, we first notice that, by Proposition~\ref{local}, if $r \notin \nn$ and $k > \max \{\mathsf{n}(r), \mathsf{d}(r)\}$, then $\rho_k(M_r) = \infty$. Hence the fact that $\rho_k(M_r) < \infty$ for every $k \in \nn$ implies that $r \in \nn$. In this case $\rho(M_r) = \rho(\nn_0) = 1$, and so $\rho(M_r) < \infty$.
\end{proof}

As~\cite[Proposition~1.4.2(1)]{GH06} holds for every atomic monoid, the direct implication of Corollary~\ref{cor:local elastiicity} also holds for any atomic monoid. However, the reverse implication of the same corollary is not true even in the context of Puiseux monoids.

\begin{example}
	Let $(p_n)_{n \in \mathbb{N}}$ be a strictly increasing sequence of primes, and consider the Puiseux monoid
	\[
		M := \bigg\langle \frac{p_n^2 + 1}{p_n} : n \in \nn \bigg \rangle.
	\]
	It is not hard to verify that the monoid $M$ is atomic with set of atoms given by the displayed generating set. Then it follows from \cite[Theorem~3.2]{GO19} that $\rho(M_r) = \infty$. However, \cite[Theorem~4.1(1)]{mG19} guarantees that $\rho_k(M) < \infty$ for every $k \in \nn$.
\end{example}

\medskip
\section{Tame Degree}
\label{sec:tame degree}

As the elasticity, the tameness is an arithmetic tool to measure how far is an atomic monoid from being a UFM. Although the tameness of many classes of atomic monoids has been studied in the past (see \cite{CGL09}, \cite{CCMMP14}, \cite{GH08}), no systematic investigation of the tameness has been carried out for Puiseux monoids. For the special class of strongly primary Puiseux monoids, recent results have been achieved in~\cite[Section~3]{GGT19}. In this section, we study the tameness of the multiplicatively cyclic Puiseux monoids.

\subsection{Omega Primality} Let $M$ be a reduced atomic monoid. The \emph{omega function} $\omega \colon M \to \nn_0 \cup \{\infty\}$ is defined as follows: for each $x \in M^\bullet$ we take $\omega(x)$ to be the smallest $n \in \nn$ satisfying that whenever $x \mid_M \sum_{i=1}^t a_i$ for some $a_1, \dots, a_t \in \mathcal{A}(M)$, there exists $T \subseteq \{ 1, \dots, t \}$ with $|T| \le n$ such that $x \mid_M \sum_{i \in T} a_i$.  If no such $n$ exists, then $\omega(x) = \infty$. In addition, we define $\omega(0) = 0$. Then we define
\[
	\omega(M) := \sup\{\omega(a) : a \in \mathcal{A}(M)\}.
\]
Notice that $\omega(x) = 1$ if and only if $x$ is prime in $M$. The omega function was introduced by Geroldinger and Hassler in~\cite{GH08} to measure how far in an atomic monoid an element is from being prime.

Before proving the main results of this section, let us collect two technical lemmas.

\begin{lemma} \label{lem:element divisible by 1}
	If $r \in \qq_{> 1}$, then $1 \mid_{M_r} \mathsf{d}(r) r^k$ for every $k \in \nn_0$.
\end{lemma}

\begin{proof}
	If $r \in \nn$, then $M_r = (\nn_0,+)$ and the statement of the lemma follows straightforwardly. Then we assume that $r \in \qq_{> 1} \setminus \nn$. For $k=0$, the statement of the lemma holds trivially. For $k \in \nn$, consider the factorization $z_k := \mathsf{d}(r) \, r^k \in \mathsf{Z}(M_r)$. The factorization
	\[
		z := \mathsf{n}(r) + \sum_{i=1}^{k-1} (\mathsf{n}(r) - \mathsf{d}(r)) r^i
	\]
	belongs to $\mathsf{Z}(\phi(z_k))$ (recall that $\phi \colon \mathsf{Z}(M_r) \to M_r$ is the factorization homomorphism of~$M_r$). This is because
	\begin{align*}
		\mathsf{n}(r) + \sum_{i=1}^{k-1} (\mathsf{n}(r) - \mathsf{d}(r)) r^i
		&= \mathsf{n}(r) + \sum_{i=1}^{k-1} \mathsf{n}(r) r^i - \sum_{i=1}^{k-1} \mathsf{d}(r) r^i \\ &= \mathsf{n}(r) + \sum_{i=1}^{k-1} \mathsf{n}(r) r^i - \sum_{i=1}^{k-1} \mathsf{n}(r) r^{i-1} = \mathsf{d}(r) r^k.
	\end{align*}
	Hence $1 \mid_{M_r} \mathsf{d}(r) r^k$
\end{proof}

\begin{lemma} \label{lem:1 dividies constant coefficient of min-length factorization}
	Take $r \in \qq \cap (0,1)$ such that $M_r$ is atomic, and let $\sum_{i=0}^N \alpha_i r^i$ be the factorization in $\mathsf{Z}(x)$ of minimum length. Then $\alpha_0 \ge 1$ if and only if $1 \mid_{M_r} x$.
\end{lemma}

\begin{proof}
	The direct implication is straightforward. For the reverse implication, suppose that $1 \mid_{M_r} x$. Then there exists a factorization $z' := \sum_{i=0}^K \beta_i r^i \in \mathsf{Z}(x)$ such that $\beta_0 \ge 1$. If $\beta_i \ge \mathsf{d}(r)$ for some $i \in \{ 1, \dots, K \}$, then we can use the identity $\mathsf{d}(r) r^i = \mathsf{n}(r) r^{i-1}$ to find another factorization $z'' \in \mathsf{Z}(x)$ such that $|z''| < |z'|$. Notice that the atom $1$ appears in~$z''$. Then we can replace $z'$ by $z''$. After carrying out such a replacement as many times as possible, we can guarantee that $\beta_i < \mathsf{d}(r)$ for $i \in \{ 1, \dots, K \}$. Then Lemma~\ref{lem:factorization of extremal length II}(1) ensures that $z'$ is a minimum-length factorization of $x$. Now Lemma~\ref{lem:factorization of extremal length II}(2) implies that $z' = z$. Finally, $\alpha_0 = \beta_0 \ge 1$ follows from the fact that the atom $1$ appears in~$z'$.
\end{proof}

\begin{proposition} \label{prop:omega primality}
	Take $r \in \qq_{> 0}$ such that  $M_r$ is atomic.
	\begin{enumerate}
		\item If $r<1$, then $\omega(1) = \infty$.
		
		\item If $r \in \nn$, then $\omega(1) = 1$.
		
		\item If $r \in \qq_{> 1} \setminus \nn$, then $\omega(1) = \mathsf{d}(r)$.
	\end{enumerate}
\end{proposition}

\begin{proof}
	To verify~(1), suppose that $r < 1$. Then set $x = \mathsf{n}(r) \in M_r$ and note that $1 \mid_{M_r} x$. Fix an arbitrary $N \in \nn$. Take now $n \in \nn$ such that $\mathsf{d}(r) + n( \mathsf{d}(r) - \mathsf{n}(r)) \ge N$. It is not hard to check that
	\[
		z := \mathsf{d}(r) r^{n+1} + \sum_{i=1}^n (\mathsf{d}(r) - \mathsf{n}(r) ) r^i
	\]
	is a factorization in $\mathsf{Z}(x)$. Suppose that $z' = \sum_{i=1}^K \alpha_i r^i$ is a sub-factorization of $z$ such that $1 \mid_{M_r} x' := \phi(z')$. Now we can move from $z'$ to a factorization $z''$ of $x'$ of minimum length by using the identity $\mathsf{d}(r)r^{i+1} = \mathsf{n}(r)r^i$ finitely many times. As $1 \mid_{M_r} x'$, it follows by Lemma~\ref{lem:1 dividies constant coefficient of min-length factorization} that the atom $1$ appears in $z''$. Therefore, when we obtained $z''$ from~$z'$ (which does not contain $1$ as a formal atom), we must have applied the identity $\mathsf{d}(r)r = \mathsf{n}(r) \cdot 1$ at least once. As a result $z''$ contains at least $\mathsf{n}(r)$ copies of the atom~$1$. This implies that $x' = \phi(z'') \ge \mathsf{n}(r) = x$. Thus, $x' = x$, which implies that $z'$ is the whole factorization~$z$. As a result, $\omega(1) \ge |z| \ge N$. Since $N$ was arbitrarily taken, we can conclude that $\omega(1) = \infty$, as desired.
	
	Notice that~(2) is a direct consequence of the fact that $1$ is a prime element in $M_r = (\nn_0,+)$.
	
	Finally, we prove~(3). Take $z = \sum_{i=0}^N \alpha_i r^i \in \mathsf{Z}(x)$ for some $x \in M_r$ such that $1 \mid_{M_r} x$. We claim that there exists a sub-factorization $z'$ of $z$ such that $|z'| \le \mathsf{d}(r)$ and $1 \mid_{M_r} \phi(z')$, where $\phi$ is the factorization homomorphism of $M_r$. If $\alpha_0 > 0$, then $1$ is one of the atoms showing in $z$ and our claim follows trivially. Therefore assume that $\alpha_0 = 0$. Since $1 \mid_{M_r} x$ and $1$ does not show in $z$, we have that $|\mathsf{Z}(x)| > 1$. Then conditions~(1) and~(3) in Lemma~\ref{lem:factorization of extremal length I} cannot be simultaneously true, which implies that $\alpha_i \ge \mathsf{d}(r)$ for some $i \in \{ 1, \dots, N \}$. Lemma~\ref{lem:element divisible by 1} ensures now that $1 \mid_{M_r} \phi(z')$ for the sub-factorization $z' := \mathsf{d}(r)r^i$ of $z$. This proves our claim and implies that $\omega(1) \le \mathsf{d}(r)$. On the other hand, take $w$ to be a strict sub-factorization of $\mathsf{d}(r) \, r$. Note that the atom $1$ does not appear in $w$. In addition, it follows by Lemma~\ref{lem:factorization of extremal length I} that $|\mathsf{Z}(\phi(w))| = 1$. Hence $1 \nmid_{M_r} \phi(w)$. As a result, we have that $\omega(1) \ge \mathsf{d}(r)$, and~(3) follows.
\end{proof}


\subsection{Tameness} For an atom $a \in \mathcal{A}(M)$, the {\it local tame degree} $\mathsf{t}(a) \in \nn_0$ is the smallest $n \in~\nn_0 \cup \{\infty\}$ such that in any given factorization of $x \in a + M$  at most $n$ atoms have to be replaced by at most $n$ new atoms to obtain a new factorization of $x$ that contains $a$. More specifically, it means that $\mathsf{t}(a)$ is the smallest $n \in \nn_0 \cup \{\infty\}$ with the following property: if $\mathsf{Z}(x) \cap (a + \mathsf{Z}(M)) \ne \emptyset$ and $z \in \mathsf{Z}(x)$, then there exists a $z' \in \mathsf{Z}(x) \cap (a + \mathsf{Z}(M))$ such that $\mathsf{d}(z,z') \le n$.

\begin{definition}
	An atomic monoid $M$ is said to be {\it locally tame} provided that $\mathsf{t}(a) < \infty$ for all $a \in \mathcal{A}(M)$.
\end{definition}

Every factorial monoid is locally tame (see \cite[Theorem~1.6.6 and Theorem~1.6.7]{GH06}). In particular, $(\nn_0,+)$ is locally tame. The tame degree of numerical monoids was first considered in~\cite{CGL09}. 
The factorization invariant $\tau \colon M \to \nn_0 \cup \{\infty\}$, which was introduced in~\cite{GH08}, is defined as follows: for $k \in \nn$ and $b \in M$, we take
\[
	\mathsf{Z}_{\text{min}}(k,b) := \bigg\{ \sum_{i=1}^j a_i \in \mathsf{Z}(M) : j \le k, \ b \mid_M \sum_{i=1}^j a_i, \, \text{ and } \, b \nmid_M \sum_{i \in I} a_i \ \text{ for any } \ I \subsetneq \{ 1, \dots, j \} \bigg\}
\]
and then we set
\[
	\tau(b) = \sup_k \sup_z \big\{ \min \mathsf{L}\big(\phi(z) - b \big) : z \in \mathsf{Z}_{\text{min}}(k,b)\big\}.
\]
The monoid $M$ is called {\it (globally) tame} provided that the \emph{tame degree}
\[
	\mathsf{t}(M) = \sup \{\mathsf{t}(a) : a \in \mathcal{A}(M)\} < \infty.
\]

The following result will be used in the proof of Theorem~\ref{thm:cyclic semirings are no locally tame}.

\begin{theorem} \cite[Theorem~3.6]{GH08} \label{thm:characterization of locally tame monoids}
	Let $M$ be a reduced atomic monoid. Then $M$ is locally tame if and only if $\omega(a) < \infty$ and $\tau(a) < \infty$ for all $a \in \mathcal{A}(M)$.
\end{theorem}

We conclude this section by characterizing the multiplicatively cyclic Puiseux monoids that are locally tame.

\begin{theorem} \label{thm:cyclic semirings are no locally tame}
	Take $r \in \qq_{>0}$ such that $M_r$ is atomic. Then the following conditions are equivalent:
	\begin{enumerate}
		\item $r \in \nn$;
		
		\item $\omega(M_r) < \infty$;
		
		\item $M_r$ is globally tame;
		
		\item $M_r$ is locally tame.
	\end{enumerate}
\end{theorem}

\begin{proof}
	That (1) implies (2) follows from Proposition~\ref{prop:omega primality}(2). Now suppose that (2) holds. Then~\cite[Proposition~3.5]{GK10} ensures that $\mathsf{t}(M_r) \le \omega(M_r)^2 < \infty$, which implies~(3). In addition, (3) implies (4) trivially. 
	
	To prove that (4) implies (1) suppose, by way of contradiction, that $r \in \qq_{> 0} \setminus \nn$. Let us assume first that $r < 1$. In this case, $\omega(1) = \infty$ by Proposition~\ref{prop:omega primality}(3). Then it follows by Theorem~\ref{thm:characterization of locally tame monoids} that $M_r$ is not locally tame, which is a contradiction. For the rest of the proof, we assume that $r \in \qq_{> 1} \setminus \nn$.
	
	We proceed to show that $\tau(1) = \infty$. For $k \in \nn$ such that $k \ge \mathsf{d}(r)$, consider the factorization $z_k = \mathsf{d}(r) r^k \in \mathsf{Z}(M_r)$. Since any strict sub-factorization $z'_k$ of $z_k$ is of the form $\beta r^k$ for some $\beta < \mathsf{d}(r)$, it follows by Lemma~\ref{lem:factorization of extremal length I} that $|\mathsf{Z}(z'_k)| = 1$. On the other hand, $1 \mid_{M_r} \mathsf{d}(r) r^k$ by Lemma~\ref{lem:element divisible by 1}. Therefore $z_k \in \mathsf{Z}_{\text{min}}(k, 1)$. Now consider the factorization
	\[
		z'_k := (\mathsf{n}(r) - 1) \cdot 1 + \sum_{i=1}^{k-1} (\mathsf{n}(r) - \mathsf{d}(r)) r^i.
	\]
	Proceeding as in the proof of Lemma~\ref{lem:element divisible by 1}, one can verify that $\phi(z'_k) = \mathsf{d}(r)r^k - 1$. In addition, the coefficients of the atoms $1, \dots, r^{k-1}$ in $z'_k$ are all strictly less than~$\mathsf{n}(r)$. Then it follows from Lemma~\ref{lem:factorization of extremal length I}(1) that $z'_k$ is a factorization of $\mathsf{d}(r)r^k - 1$ of minimum length. Because $|z'_k| = k(\mathsf{n}(r) - \mathsf{d}(r)) + \mathsf{d}(r) - 1$, one has that
	\begin{align*}
		\tau(1) 
		&= \sup_k \sup_z \big\{ \min \mathsf{L}\big(\phi(z) - 1 \big) : z \in \mathsf{Z}_{\text{min}}(k,1)\big\} \\
		&\ge \sup_k \min \mathsf{L}\big( \phi(z_k) - 1 \big) = \sup_k |z'_k| \\
		&= \lim_{k \to \infty} k(\mathsf{n}(r) - \mathsf{d}(r)) + \mathsf{d}(r) - 1 \\
		&= \infty.
	\end{align*}
	Hence $\tau(1) = \infty$. Then it follows by Theorem~\ref{thm:characterization of locally tame monoids} that $M_r$ is not locally tame, which contradicts condition~(3). Thus, (3) implies (1), as desired.
\end{proof}

%% file: tex/ch5.tex
\chapter{Factorial Elements of Puiseux Monoids} \label{conclusion}

\section{Introduction}
\label{sec:intro}

The elements having exactly one factorization are crucial in the study of factorization theory of commutative cancellative monoids and integral domains. Aiming to avoid repeated long descriptions, we call such elements \emph{molecules}. Molecules were first studied in the context of algebraic number theory by W.~Narkiewicz and other authors in the 1960's. For instance, in~\cite{wN66} and~\cite{wN66a} Narkiewicz studied some distributional aspects of the molecules of quadratic number fields. In addition, he gave an asymptotic formula for the number of (non-associated) integer molecules of any algebraic number field~\cite{wN72}. In this chapter, we study the molecules of submonoids of~$(\qq_{\ge 0},+)$, including numerical monoids, and the molecules of their corresponding monoid algebras.

\begin{figure}
	\centering
	\includegraphics[width = 11cm]{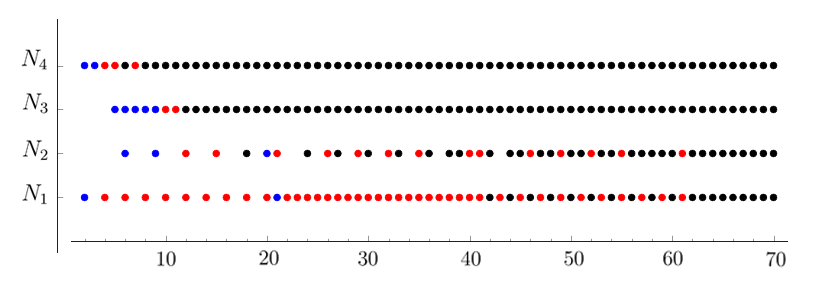}
	\caption[Atoms, molecules, and non-molecules of the numerical monoids $N_1 = \langle 2, 21 \rangle$, $N_2 = \langle 6,9,20 \rangle$, $N_3 = \langle 5,6,7,8,9 \rangle$, and $N_4 = \langle 2,3 \rangle$.]{The dots on the horizontal line labeled by $N_i$ represent the nonzero elements of the numerical monoid $N_i$; here we are setting $N_1 = \langle 2, 21 \rangle$, $N_2 = \langle 6,9,20 \rangle$, $N_3 = \langle 5,6,7,8,9 \rangle$, and $N_4 = \langle 2,3 \rangle$.
		Atoms are represented in blue, molecules that are not atoms in red, and non-molecules in black.}
	\label{fig:molecule of four NS}
\end{figure}

If a numerical monoid $N$ satisfies that $N \neq \nn_0$, then it contains only finitely many molecules. Notice, however, that every positive integer is a molecule of $(\nn_0,+)$. Figure~\ref{fig:molecule of four NS} shows the distribution of the sets of molecules of four numerical monoids. We begin Section~\ref{sec:molecules of NS} pointing out how the molecules of numerical monoids are related to the Betti elements. Then we show that each element in the set $\nn_{\ge 4} \cup \{\infty\}$ (and only such elements) can be the number of molecules of a numerical monoid. We conclude our study of molecules of numerical monoids exploring the possible cardinalities of the sets of reducible molecules (i.e., molecules that are not atoms).

The class of Puiseux monoids, on the other hand, contains members having infinitely many atoms and, consequently, infinitely many molecules. In Section~\ref{sec:molecules of PM}, we study the sets of molecules of Puiseux monoids, finding infinitely many non-isomorphic Puiseux monoids all whose molecules are atoms (in contrast to the fact that the set of molecules of a numerical monoid always differs from its set of atoms). 

We conclude with Section~\ref{sec:mol prime reciprocals}, where we construct infinitely many non-isomorphic Puiseux monoids having infinitely many molecules that are not atoms (in contrast to the fact that the set of molecules of a nontrivial numerical monoid is always finite). Special attention is given in this section to prime reciprocal Puiseux monoids and a characterization of their molecules.
\medskip

\section{Molecules of Numerical Monoids}
\label{sec:molecules of NS}

In this section we study the sets of molecules of numerical monoids, putting particular emphasis on their possible cardinalities. 

\subsection{Atoms and Molecules}

As one of the main purposes of this chapter is to study elements with exactly one factorization in Puiseux monoids (in particular, numerical monoids), we introduce the following definition.

\begin{definition}
	Let $M$ be a monoid. We say that an element $x \in M \setminus U(M)$ is a \emph{molecule} provided that $|\mathsf{Z}(x)| = 1$. The set of all molecules of $M$ is denoted by $\mathcal{M}(M)$.
\end{definition}

It is clear that the set of atoms of any monoid is contained in the set of molecules. However, such an inclusion might be proper (consider, for instance, the additive monoid $\nn_0$). In addition, for any atomic monoid $M$ the set $\mathcal{M}(M)$ is \emph{divisor-closed} in the sense that if $x \in \mathcal{M}(M)$ and $x' \mid_M x$ for some $x' \in M \setminus U(M)$, then $x' \in \mathcal{M}(M)$. If the condition of atomicity is dropped, then this observation is not necessarily true (see Example~\ref{ex:PM whose set of molecules is not divisor-closed}).
\medskip

%

\begin{example} \label{ex:NS with embedding dimension two}
	For $k \ge 1$, consider the numerical monoid $N_1 = \langle 2, 21 \rangle$, whose molecules are depicted in Figure~\ref{fig:molecule of four NS}. It is not hard to see that $x \in N_1^\bullet$ is a molecule if and only if every factorization of $x$ contains at most one copy of $21$. Therefore
	\[
	\mathcal{M}(N_1) = \big\{2m + 21n : 0 \le m < 21, n \in \{0,1\}, \, \text{and} \, (m, n) \neq (0, 0)\big\}.
	\]
	In addition, if $2m + 21n = 2m' + 21n'$ for some $m,m' \in \{0,\dots, 20\}$ and $n,n' \in \{0,1\}$, then one can readily check that $m = m'$ and $n = n'$. Hence $|\mathcal{M}(N_1)| = 41$.
\end{example}
\smallskip

\subsection{Betti Elements} Let $N = \langle a_1, \dots, a_n \rangle$ be a minimally generated numerical monoid. We always represent an element of $\mathsf{Z}(N)$ with an $n$-tuple $z = (c_1, \dots, c_n) \in \nn^n_0$, where the entry $c_i$ specifies the number of copies of $a_i$ that appear in $z$. Clearly, $|z| = c_1 + \dots + c_n$. Given factorizations $z = (c_1, \dots, c_n)$ and $z' = (c'_1, \dots, c'_n)$, we define
\[
	\gcd(z,z') = (\min\{c_1,c'_1\}, \dots, \min\{c_n,c'_n\}).
\]
The \emph{factorization graph} of $x \in N$, denoted by $\nabla\!_x(N)$ (or just $\nabla\!_x$ when no risk of confusion exists), is the graph with vertices $\mathsf{Z}(x)$ and edges between those $z, z' \in \mathsf{Z}(x)$ satisfying that $\gcd(z,z') \neq 0$. The element $x$ is called a \emph{Betti element} of $N$ provided that $\nabla\!_x$ is disconnected. The set of Betti elements of $N$ is denoted by $\text{Betti}(N)$.

\begin{example}
	Take $N$ to be the numerical monoid $\langle 14, 16, 18, 21, 45 \rangle$. A computation in SAGE using the \texttt{numericalsgps GAP} package immediately reveals that $N$ has nine Betti elements. In particular, $90 \in \text{Betti}(N)$. In Figure~\ref{fig:factorization graphs} one can see the disconnected factorization graph of the Betti element $90$ on the left and the connected factorization graph of the non-Betti element $84$ on the right. 
	\begin{figure}
		\centering     
		\subfigure[\label{figa}]{\includegraphics[width = 5.4cm]{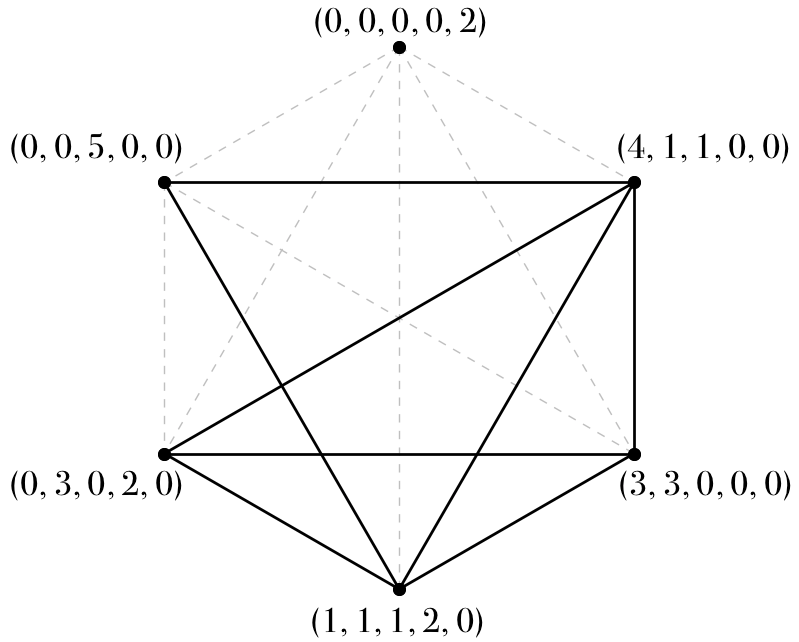}}
		\hspace{10pt}
		\subfigure[]{\includegraphics[width = 5.4cm]{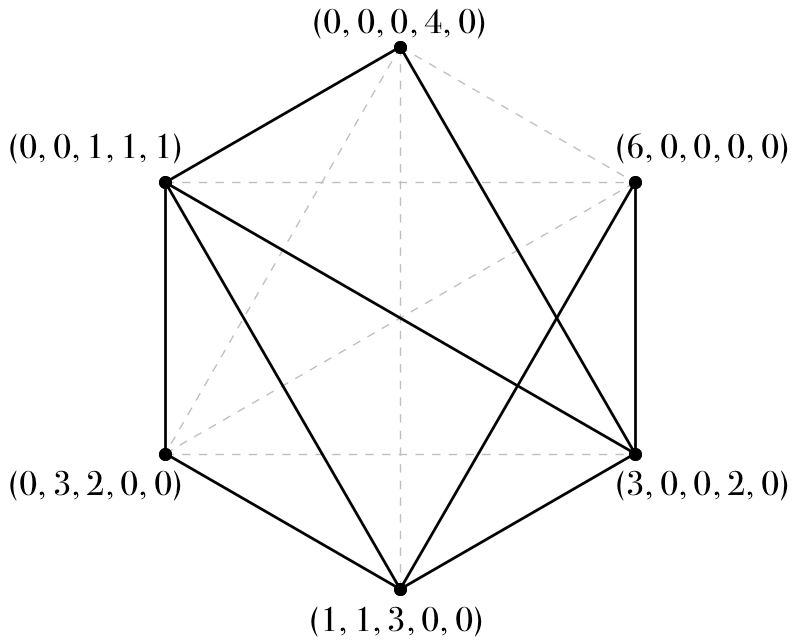}}
		\caption[The factorization graphs of $90 \in \text{Betti}(N)$ and $84 \notin \text{Betti}(N)$, where $N$ is the numerical monoid $\langle 14,16,18,21,45 \rangle$.]{Factorization graphs of one Betti element and one non-Betti element in the numerical monoid $N = \langle 14,16,18,21,45 \rangle$: A) the element $90 \in \text{Betti}(N)$ \\and B) the element $84 \notin \text{Betti}(N)$. }
		\label{fig:factorization graphs}
	\end{figure}
\end{example}

Observe that $0 \notin \text{Betti}(N)$ since $|\mathsf{Z}(0)| = 1$. It is well-known that every numerical monoid has finitely many Betti elements. Betti elements play a fundamental role in the study of uniquely-presented numerical monoids \cite{GO10} and the study of delta sets of BFMs~\cite{CGLMS12}. In a numerical monoid, Betti elements and molecules are closely related.

\begin{remark}
	Let $N$ be a numerical monoid. An element $m \in N$ is a molecule if and only if $\beta \nmid_N m$ for any $\beta \in \text{Betti}(N)$.
\end{remark}

\begin{proof}
	For the direct implication, suppose that $m$ is a molecule of $N$ and take $\alpha \in N$ such that $\alpha \mid_N m$. As the set of molecules is closed under division, $|\mathsf{Z}(\alpha)| = 1$. This implies that $\nabla \!_\alpha$ is connected and, therefore, $\alpha$ cannot be a Betti element. The reverse implication is just a rephrasing of \cite[Lemma~1]{GO10}.
\end{proof}
\medskip

\subsection{On the Sizes of the Sets of Molecules} Obviously, for every $n \in \nn$ there exists a numerical monoid having exactly $n$ atoms. The next proposition answers the same realization question replacing the concept of an atom by that one of a molecule. Recall that $\mathcal{N}^\bullet$ denotes the class of all nontrivial numerical monoids.

\begin{proposition} \label{prop:realization for the size of the sets of molecules}
	$\{|\mathcal{M}(N)| : N \in \mathcal{N}^\bullet\} = \nn_{\ge 4}$.
\end{proposition}

\begin{proof}
	Let $N$ be a nontrivial numerical monoid. Then $N$ must contain at least two atoms. Let $a$ and $b$ denote the two smallest atoms of $N$, and assume that $a < b$. Note that $2a$ and $a+b$ are distinct molecules that are not atoms. Hence $|\mathcal{M}(N)| \ge 4$. As a result, $\{|\mathcal{M}(N)| : N \in \mathcal{N}^\bullet\} \subseteq \nn_{\ge 4} \cup \{ \infty \}$. Now take $x \in \nn$ with $x > \mathfrak{f}(N) + a b$. Since $x' := x - a b > \mathfrak{f}(N)$, we have that $x' \in N$ and, therefore, $\mathsf{Z}(x')$ contains at least one factorization, namely $z$. So we can find two distinct factorizations of $x$ by adding to $z$ either $a$ copies of $b$ or $b$ copies of $a$. Thus, $\mathfrak{f}(N) + a b$ is an upper bound for $\mathcal{M}(N)$, which means that $|\mathcal{M}(N)| \in \nn_{\ge 4}$. Thus, $\{|\mathcal{M}(N)| : N \in \mathcal{N}^\bullet\} \subseteq \nn_{\ge 4}$.
	
	To argue the reverse inclusion, suppose that $n \in \nn_{\ge 4}$, and let us find $N \in \mathcal{N}$ with $|\mathcal{M}(N)| = n$. For $n = 4$, we can take the numerical monoid $\langle 2, 3 \rangle$ (see Figure~\ref{fig:molecule of four NS}). For $n > 4$, consider the numerical monoid
	\[
		N = \langle n-2, n-1, \dots, 2(n-2)-1 \rangle.
	\]
	It follows immediately that $\mathcal{A}(N) = \{n-2, n-1, \dots, 2(n-2)-1\}$. In addition, it is not hard to see that $2(n-2), 2(n-2)+1 \in \mathcal{M}(N)$ while $k \notin \mathcal{M}(N)$ for any $k > 2(n-2)+1$. Consequently, $\mathcal{M}(N) = \mathcal{A}(N) \cup \{2(n-2), 2(n-2)+1\}$, which implies that $|\mathcal{M}(N)| = n$. Therefore $\{|\mathcal{M}(N)| : N \in \mathcal{N}\} \supseteq \nn_{\ge 4}$, which completes the proof.
\end{proof}

\begin{corollary}
	The monoid $(\nn_0,+)$ is the only numerical monoid having infinitely many molecules.
\end{corollary}

In Proposition~\ref{prop:realization for the size of the sets of molecules} we have fully described the set $\{|\mathcal{M}(N)| : N \in \mathcal{N} \}$. A full description of the set $\{|\mathcal{M}(N) \setminus \mathcal{A}(N)| : N \in \mathcal{N} \}$ seems to be significantly more involved. However, the next theorem offers some evidence to believe that
\[
	\{|\mathcal{M}(N) \setminus \mathcal{A}(N)| : N \in \mathcal{N} \} = \nn_{\ge 2} \cup \{\infty\}.
\]

\begin{theorem} \label{thm:more molecules than atoms in any NS}
	The following statements hold.
	\begin{enumerate}
		\item $\{|\mathcal{M}(N) \setminus \mathcal{A}(N)| : N \in \mathcal{N}^\bullet \} \subseteq \nn_{\ge 2}$.
		\vspace{3pt}
		\item $|\mathcal{M}(N) \! \setminus \! \mathcal{A}(N)| = 2$ for infinitely many numerical monoids $N$.
		\vspace{3pt}
		\item For each $k \in \nn$, there is a numerical monoid $N_k$ with $|\mathcal{M}(N) \! \setminus \! \mathcal{A}(N)| > k$.
	\end{enumerate}
\end{theorem}

\begin{proof}
	To prove~(1), take $N \in \mathcal{N}^\bullet$. Then we can assume that $N$ has embedding dimension $n$ with $n \ge 2$. Take $a_1, \dots, a_n \in \nn$ such that $a_1 < \dots < a_n$ such that $N = \langle a_1, \dots, a_n \rangle$. Since $a_1 < a_2 < a_j$ for every $j = 3, \dots,n$, the elements $2a_1$ and $a_1 + a_2$ are two distinct molecules of $N$ that are not atoms. Hence $\mathcal{M}(N) \setminus \mathcal{A}(N) \subseteq \nn_{\ge 2} \cup \{\infty\}$. On the other hand, Proposition~\ref{prop:realization for the size of the sets of molecules} guarantees that $|\mathcal{M}(N)| < \infty$, which implies that $|\mathcal{M}(N) \setminus \mathcal{A}(N)| < \infty$. As a result, the statement~(1) follows.
	
	To verify the statement~(2), one only needs to consider for every $n \in \nn$ the numerical monoid $N_n := \{0\} \cup \nn_{\ge n-2}$. The minimal set of generators of $N_n$ is the $(n-2)$-element set $\{n-2, n-1, \dots, 2(n-2)-1 \}$ and, as we have already argued in the proof of Proposition~\ref{prop:realization for the size of the sets of molecules}, the set $\mathcal{M}(N_n) \! \setminus \! \mathcal{A}(N_n)$ consists precisely of two elements.
	
	Finally, let us prove condition~(3). We first argue that for any $a,b \in ~\nn_{\ge 2}$ with $\gcd(a,b) = 1$ the numerical monoid $\langle a, b \rangle$ has exactly $ab-1$ molecules (cf. Example~\ref{ex:NS with embedding dimension two}). Assume $a < b$, take $N := \langle a, b \rangle $, and set
	\[
		\mathcal{M} = \{ma + nb : 0 \le m < b, \, 0 \le n < a, \, \text{and} \ \, (m, n) \neq (0, 0)\}.
	\]
	Now take $x \in N$ to be a molecule of $N$. As $|\mathsf{Z}(x)| = 1$, the unique factorization $z := ~(c_1, c_2) \in \mathsf{Z}(x)$ (with $c_1,c_2 \in \nn_0$) satisfies that $c_1 < b$; otherwise, we could exchange $b$ copies of the atom $a$ by $a$ copies of the atom $b$ to obtain another factorization of $x$. A similar argument ensures that $c_2 < a$. As a consequence, $\mathcal{M}(N) \subseteq \mathcal{M}$. On the other hand, if $ma + nb = m'a + n'b$ for some $m,m',n,n' \in \nn_0$, then $\gcd(a,b) = 1$ implies that $b \mid m-m'$ and $a \mid n-n'$. Because of this observation, the element $(b-1)a + (a-1)b$ has only the obvious factorization, namely $(b-1,a-1)$. Since $(b-1)a + (a-1)b$ is a molecule satisfying that $y \mid_N (b-1)a + (a-1)b$ for every $y \in \mathcal{M}$, the inclusion $\mathcal{M} \subseteq \mathcal{M}(N)$ holds. Hence $|\mathcal{M}(N)| = |\mathcal{M}| = ab-1$. To argue the statement~(3) now, it suffices to take $N_k := \langle 2, 2k+1 \rangle$.
\end{proof}
\medskip

%

\section{Molecules of Generic Puiseux Monoids}
\label{sec:molecules of PM}
\smallskip

 In this section we study the sets of molecules of the general class of Puiseux monoids. We will argue that there are infinitely many non-finitely generated atomic Puiseux monoids~$P$ such that $|\mathcal{M}(P) \setminus \mathcal{A}(P)| = \infty$. On the other hand, we will prove that, unlike the case of numerical monoids, there are infinitely many non-isomorphic atomic Puiseux monoids all whose molecules are, indeed, atoms.

In Section~\ref{sec:molecules of NS} we mentioned that the set of molecules of an atomic monoid is divisor-closed. The next example indicates that this property might not hold for non-atomic monoids.

\begin{example} \label{ex:PM whose set of molecules is not divisor-closed}
	Consider the Puiseux monoid
	\[
		P = \bigg\langle \frac 25, \frac 35, \frac 1{2^n} \, : \, n \in \nn \bigg\rangle.
	\]
	First, observe that $0$ is not a limit point of $P^\bullet$, and so $P$ cannot be finitely generated. After a few easy verifications, one can see that $\mathcal{A}(P) = \{2/5, 3/5\}$. On the other hand, it is clear that $1/2 \notin \langle 2/5,3/5 \rangle$, so $P$ is not atomic. Observe now that~$\mathsf{Z}(1)$ contains only one factorization, namely $2/5 + 3/5$. Therefore $1 \in \mathcal{M}(P)$. Since $\mathsf{Z}(1/2)$ is empty, $1/2$ is not a molecule of $P$. However, $1/2 \mid_P 1$. As a result, $\mathcal{M}(P)$ is not divisor-closed.
\end{example}

Although the additive monoid $\nn_0$ contains only one atom, it has infinitely many molecules. The next result implies that $\nn_0$ is basically the only atomic Puiseux monoid having finitely many atoms and infinitely many molecules. 

\begin{proposition} \label{prop:atoms-molecules cardinality for PM}
	Let $P$ be a Puiseux monoid. Then $|\mathcal{M}(P)| \in \nn_{\ge 2}$ if and only if $|\mathcal{A}(P)| \in \nn_{\ge 2}$. 
\end{proposition}

\begin{proof}
	Suppose first that $|\mathcal{M}(P)| \in \nn_{\ge 2}$. As every atom is a molecule, $\mathcal{A}(P)$ is finite. Furthermore, note that if $\mathcal{A}(P) = \{a\}$, then every element of the set $S = \{na : n \in \nn\}$ would be a molecule, which is not possible as $|S| = \infty$. As a result, $|\mathcal{A}(P)| \in \nn_{\ge 2}$. Conversely, suppose that $|\mathcal{A}(P)| \in \nn_{\ge 2}$. Since the elements in $P \! \setminus \! \langle \mathcal{A}(P) \rangle$ have no factorizations, $\mathcal{M}(P) = \mathcal{M}(\langle \mathcal{A}(P) \rangle)$. Therefore there is no loss in assuming that $P$ is atomic. As $1 < |\mathcal{A}(P)| < \infty$, the monoid $P$ is isomorphic to a nontrivial numerical monoid. The proposition now follows from the fact that nontrivial numerical monoids have finitely many molecules.
\end{proof}

\begin{corollary}
	If $P$ is a Puiseux monoid, then $|\mathcal{M}(P)| \neq 1$.
\end{corollary}

The set of atoms of a numerical monoid is always strictly contained in its set of molecules. However, there are many atomic Puiseux monoids which do not satisfy such a property. Before proceeding to formalize this observation, recall that if two Puiseux monoids $P$ and $P'$ are isomorphic, then there exists $q \in \qq_{>0}$ such that $P' = qP$; this is a consequence of~Proposition~\ref{chap:algebraic background}.\ref{prop:homomorphisms between PM}.

\begin{theorem}[cf. Theorem~\ref{thm:more molecules than atoms in any NS}(1)] \label{thm:F(M) equals A(M)}
	There are infinitely many non-isomorphic atomic Puiseux monoids $P$ satisfying that $\mathcal{M}(P) = \mathcal{A}(P)$.
\end{theorem}

\begin{proof}
	Let $\mathcal{S} = \{S_n : n \in \nn\}$ be a collection of infinite and pairwise-disjoint sets of primes. Now take $S = S_n$ for some arbitrary $n \in \nn$, and label the primes in $S$ strictly increasingly by $p_1, p_2, \dots$. Recall that $\mathsf{D}_S(r)$ denotes the set of primes in $S$ dividing $\mathsf{d}(r)$ and that $\mathsf{D}_S(R) = \cup_{r \in R} \mathsf{D}_S(r)$ for $R \subseteq \qq_{> 0}$. We proceed to construct a Puiseux monoid $P_S$ satisfying that $\mathsf{D}_S(P_S) = S$.
	
	Take $P_1 := \langle 1/p_1 \rangle$ and $P_2 := \langle P_1, 2/(p_1p_2) \rangle$. In general, suppose that $P_k$ is a finitely generated Puiseux monoid such that $\mathsf{D}_S(P_k) \subset S$, and let $r_1, \dots, r_{n_k}$ be all the elements in $P_k$ which can be written as a sum of two atoms. Clearly, $n_k \ge 1$. Because $|S| = \infty$, one can take $p'_1, \dots, p'_{n_k}$ to be primes in $S \! \setminus \! \mathsf{D}_S(P_k)$ satisfying that $p'_i \nmid \mathsf{n}(r_i)$. Now consider the following finitely generated Puiseux monoid
	\[
		P_{k+1} := \bigg\langle P_k \cup \bigg\{ \frac{r_1}{p'_1}, \dots, \frac{r_{n_k}}{p'_{n_k}} \bigg\} \bigg\rangle.
	\]
	For every $i \in \{1,\dots,n_k\}$, there is only one element in $P_k \cup \{r_1/p'_1, \dots, r_{n_k}/p'_{n_k}\}$ whose denominator is divisible by $p'_i$, namely $r_i/p'_i$. Therefore $r_i/p'_i \in \mathcal{A}(P_{k+1})$ for $i=1, \dots, n_k$. To check that $\mathcal{A}(P_k) \subset \mathcal{A}(P_{k+1})$, fix $a \in \mathcal{A}(P_k)$ and take
	\begin{equation} \label{eq:FM equals AM}
		z := \sum_{i=1}^m \alpha_i a_i + \sum_{i=1}^{n_k} \beta_i \frac{r_i}{p'_i} \in \mathsf{Z}_{P_{k+1}}(a),
	\end{equation}
	where $a_1, \dots, a_k$ are pairwise distinct atoms in $\mathcal{A}(P_{k+1}) \cap P_k$ and $\alpha_i, \beta_j$ are nonnegative coefficients for $i = 1,\dots,m$ and $j = 1,\dots, n_k$. In particular, $a_1, \dots, a_k \in \mathcal{A}(P_k)$. For each $i=1,\dots,n_k$, the fact that the $p'_i$-adic valuation of $a$ is nonnegative implies that $p'_i \mid \beta_i$. Hence
	\[
		a = \sum_{i=1}^m \alpha_i a_i + \sum_{i=1}^{n_k} \beta'_i r_i,
	\]
	where $\beta'_i = \beta_i/p'_i \in \nn_0$ for $i = 1,\dots,n_k$. Since $r_i \in \mathcal{A}(P_k) + \mathcal{A}(P_k)$ and $(\beta_i/p'_i)r_i \mid_{P_k} a$ for every $i=1,\dots,n_k$, one obtains that $\beta_1 = \dots = \beta_{n_k} = 0$. As a result, $a = \sum_{i=1}^m \alpha_i a_i$. Because $a \in \mathcal{A}(P_k)$, the factorization $\sum_{i=1}^m \alpha_i a_i$ in $\mathsf{Z}_{P_k}(a)$ must have length $1$, i.e, $\sum_{i=1}^m \alpha_i = 1$. Thus, $\sum_{i=1}^m \alpha_i + \sum_{i=1}^{n_k} \beta_i =~1$, which means that $z$ has length $1$ and so $a \in \mathcal{A}(P_{k+1})$. As a result, the inclusion $\mathcal{A}(P_k) \subseteq \mathcal{A}(P_{k+1})$ holds. Observe that because $n_k \ge 1$, the previous containment must be strict. Now set
	\[
		P_S = \bigcup_{k \in \nn} P_k.
	\]
	Let us verify that $P_S$ is an atomic monoid satisfying that $\mathcal{A}(P_S) = \cup_{k \in \nn} \mathcal{A}(P_k)$. Since $P_k$ is atomic for every $k \in \nn$, the inclusion chain $\mathcal{A}(P_1) \subset \mathcal{A}(P_2) \subset \dots$ implies that $P_1 \subset P_2 \subset \dots$. In addition, if $a_0 = a_1 + \dots + a_m$ for $m \in \nn$ and $a_0, a_1, \dots, a_m \in P_S$, then $a_0 = a_1 + \dots + a_m$ will also hold in $P_k$ for some $k \in \nn$ large enough. This immediately implies that $\cup_{k \in \nn} \mathcal{A}(P_k) \subseteq \mathcal{A}(P_S)$. Since the reverse inclusion follows trivially, $\mathcal{A}(P_S) = \cup_{k \in \nn} \mathcal{A}(P_k)$. To check that $P_S$ is atomic, take $x \in P_S^\bullet$. Then there exists $k \in \nn$ such that $x \in P_k$ and, because $P_k$ is atomic, $x \in \langle \mathcal{A}(P_k) \rangle \subseteq \langle \mathcal{A}(P_S) \rangle$. Hence $P_S$ is atomic.
	
	To check that $\mathcal{M}(P_S) = \mathcal{A}(P_S)$, suppose that $m$ is a molecule of $P_S$, and then take $K \in \nn$ such that $m \in P_k$ for every $k \ge K$. Since $\mathcal{A}(P_k) \subset \mathcal{A}(P_{k+1}) \subset \dots$, we have that $\mathsf{Z}_{P_k}(m) \subseteq \mathsf{Z}_{P_{k+1}}(m) \subseteq \dots$. Moreover, $\cup_{k \ge K} \mathcal{A}(P_k) = \mathcal{A}(P_S)$ implies that $\cup_{k \ge K} \mathsf{Z}_{P_k}(m) = \mathsf{Z}_{P_S}(m)$. Now suppose for a contradiction that $m = \sum_{j=1}^i a_j$ for $i \in \nn_{\ge 2}$, where $a_1, \dots, a_i \in \mathcal{A}(P_S)$. Take $j \in \nn_{\ge K}$ such that $a_1, \dots, a_i \in \mathcal{A}(P_j)$. Then the way in which $P_{j+1}$ was constructed ensures that $|\mathsf{Z}_{P_{j+1}}(a_1 + a_2)| \ge 2$ and, therefore, $|\mathsf{Z}_{P_{j+1}}(m)| \ge 2$. As $\mathsf{Z}_{P_{j+1}}(m) \subseteq \mathsf{Z}_{P_S}(m)$, it follows that $|\mathsf{Z}_{P_S}(m)| \ge 2$, which contradicts that $m$ is a molecule. Hence $\mathcal{M}(P_S) = \mathcal{A}(P_S)$.
	
	Finally, we argue that the monoids constructed are not isomorphic. Let $S$ and $S'$ be two distinct members of the collection $\mathcal{S}$ and suppose, by way of contradiction, that $\psi \colon P_S \to P_{S'}$ is a monoid isomorphism. Because the only homomorphisms of Puiseux monoids are given by rational multiplications, there exists $q \in \qq_{> 0}$ such that $P_{S'} = q P_S$. In this case, all but finitely many primes in $\mathsf{D}_\pp(P_S)$ belong to $\mathsf{D}_\pp(P_{S'})$. Since $\mathsf{D}_\pp(P_S) \cap \mathsf{D}_\pp(P_{S'}) = \emptyset$ when $S \neq S'$, we get a contradiction.
\end{proof}
\medskip

\section{Molecules of Prime Reciprocal Monoids} \label{sec:mol prime reciprocals}

In this section, we focus our attention on the class consisting of all prime reciprocal Puiseux monoids.

\begin{proposition}[cf. Theorem~\ref{thm:more molecules than atoms in any NS}(1)] \label{prop:PPM has infinitely many molecules that are not atoms}
	There exist infinitely many non-finitely generated atomic Puiseux monoids $P$ such that $|\mathcal{M}(P) \! \setminus \! \mathcal{A}(P)| = \infty$.
\end{proposition}

\begin{proof}
	As in the proof of Theorem~\ref{thm:F(M) equals A(M)}, let $\mathcal{S} = \{S_n : n \in \nn\}$ be a collection of infinite and pairwise-disjoint subsets of $\pp \setminus \{2\}$. For every $n \in \nn$, let $P_n$ be a prime reciprocal Puiseux monoid over $S_n$. Fix $a \in \mathcal{A}(P_n)$, and take a factorization
	\[	
		z := \sum_{i=1}^k \alpha_i a_i  \in \mathsf{Z}(2a),
	\]
	for some $k \in \nn$, pairwise distinct atoms $a_1, \dots, a_k$, and $\alpha_1,\dots,\alpha_k \in \nn_0$. Since $\mathsf{d}(a) \neq 2$, after applying the $\mathsf{d}(a)$-adic valuation on both sides of the equality $2a = \sum_{i=1}^t \alpha_i a_i$, one obtains that $z = 2a$. So $2a \in \mathcal{M}(P_n) \setminus \mathcal{A}(P_n)$ and, as a result, $|\mathcal{M}(P_n) \setminus \mathcal{A}(P_n)| = \infty$. Now suppose, by way of contradiction, that $P_i \cong P_j$ for some $i,j \in \nn$ with $i \neq j$. Since the only isomorphisms of Puiseux monoids are given by rational multiplication, there exists $q \in \qq_{> 0}$ such that $P_j = q P_i$. However, this implies that only finitely many primes in $\mathsf{d}(P_i)$ are not contained in $\mathsf{d}(P_j)$, which contradicts that $S_i \cap S_j = \emptyset$. Hence no two monoids in $\{P_n : n \in \nn\}$ are isomorphic, and the proposition follows.
\end{proof}

%

Before characterizing the molecules of prime reciprocal monoids, let us introduce the concept of maximal multiplicity. Let $P$ be a Puiseux monoid. For $x \in P$ and $a \in \mathcal{A}(P)$ we define the \emph{maximal multiplicity} of $a$ in $x$ to be
\[
	\mathsf{m}(a,x) := \max\{n \in \nn_0 : na \mid_P x\}.
\]

\begin{proposition} \label{prop:sufficient condition for molecules in PPM}
	Let $P$ be a prime reciprocal monoid, and let $x \in P$. If $\mathsf{m}(a,x) < \mathsf{d}(a)$ for all $a \in \mathcal{A}(P)$, then $x \in \mathcal{M}(P)$.
\end{proposition}

\begin{proof}
	Suppose, by way of contradiction, that $x \notin \mathcal{M}(P)$. Then there exist $k \in \nn$, elements $\alpha_i, \beta_i \in \nn_0$ (for $i=1,\dots, k$), and pairwise distinct atoms $a_1, \dots, a_k$ such that
	\[
		z:= \sum_{i=1}^k \alpha_i a_i \ \text{ and } \ z' := \sum_{i=1}^k \beta_i a_i
	\]
	are two distinct factorizations in $\mathsf{Z}(x)$. As $z \neq z'$, there is an index $i \in \{1,\dots,k\}$ such that $\alpha_i \neq \beta_i$. Now we can apply the $\mathsf{d}(a_i)$-adic valuation to both sides of the equality
	\[
		\sum_{i=1}^k \alpha_i a_i = \sum_{i=1}^k \beta_i a_i
	\]
	to verify that $\mathsf{d}(a_i) \mid \beta_i - \alpha_i$. As $\alpha_i \neq \beta_i$, we obtain that
	\[
		\mathsf{m}(a_i, x) \ge \max\{\alpha_i, \beta_i\} \ge \mathsf{d}(a_i).
	\]
	However, this contradicts the fact that $\mathsf{m}(a,x) < \mathsf{d}(a)$ for all $a \in \mathcal{A}(P)$. As a consequence, $x \in \mathcal{M}(P)$.
\end{proof}

For $S \subseteq \pp$, we call the monoid $E_S := \langle 1/p : p \in S \rangle$ the \emph{elementary} prime reciprocal monoid over $S$; if $S = \pp$ we say that $E_S$ is \emph{the} elementary prime reciprocal monoid. It was proved in \cite[Section~5]{GG19} that every submonoid of the elementary prime reciprocal monoid is atomic. This gives a large class of non-finitely generated atomic Puiseux monoids, which contains each prime reciprocal monoid.

\begin{proposition} \label{prop:factorial equivalence for elementary primary monoids}
	Let $S$ be an infinite set of primes, and let $E_S$ be the elementary prime reciprocal monoid over $S$. For $x \in E_S$, the following conditions are equivalent:
	\begin{enumerate}
		\item $x \in \mathcal{M}(E_S)$;

		\item $1$ does not divide $x$ in $E_S$;

		\item$\mathsf{m}(a,x) < \mathsf{d}(a)$ for all $a \in \mathcal{A}(E_S)$;

		\item If $a_1, \dots, a_n \in \mathcal{A}(E_S)$ are distinct atoms and $\alpha_1, \dots, \alpha_n \in \nn_0$ satisfy that $\sum_{j=1}^n \alpha_i a_i \in \mathsf{Z}(x)$, then $\alpha_j < \mathsf{d}(a_j)$ for each $j = 1,\dots,n$.
	\end{enumerate}
\end{proposition}

\begin{proof}
	First, let us recall that since $E_S$ is atomic, $\mathcal{M}(E_S)$ is divisor-closed. On the other hand, note that for any two distinct atoms $a,a' \in \mathcal{A}(E_S)$, both factorizations $\mathsf{d}(a) \, a$ and $\mathsf{d}(a') \, a'$ are in $\mathsf{Z}(1)$. Therefore $1 \notin \mathcal{M}(E_S)$. Because the set of molecules of $E_S$ is divisor-closed, $1 \nmid_{E_S} m$ for any $m \in \mathcal{M}(E_S)$; in particular, $1 \nmid_{E_S} x$. Thus, (1) implies (2). If $\mathsf{m}(a,x) \ge \mathsf{d}(a)$ for $a \in \mathcal{A}(E_S)$, then
	\[
		x = \mathsf{m}(a,x) \, a + y = 1 + (\mathsf{m}(a,x) - \mathsf{d}(a)) \, a + y
	\]
	for some $y \in E_S$. As a result, $1 \mid_{E_S} x$, from which we can conclude that~(2) implies (3). It is obvious that (3) and (4) are equivalent conditions. Finally, the fact that (3) implies~(1) follows from Proposition~\ref{prop:sufficient condition for molecules in PPM}.
\end{proof}

\begin{corollary}
	Let $S$ be an infinite set of primes, and let $E_S$ be the elementary prime reciprocal monoid over $S$. Then $|\mathsf{Z}(x)| = \infty$ for all $x \notin \mathcal{M}(E_S)$.
\end{corollary}

In order to describe the set of molecules of an arbitrary prime reciprocal monoid, we need to cast its atoms into two categories.

\begin{definition}
	Let $P$ be a prime reciprocal monoid. We say that $a \in \mathcal{A}(P)$ is \emph{stable} if the set $\{a' \in \mathcal{A}(P) : \mathsf{n}(a') = \mathsf{n}(a)\}$ is infinite, otherwise we say that $a$ is \emph{unstable}. If every atom of $P$ is stable (resp., unstable), then we call $P$ \emph{stable} (resp., \emph{unstable}).
\end{definition}

For a prime reciprocal monoid $P$, we let $\mathcal{S}(P)$ denote the submonoid of $P$ generated by the set of stable atoms. Similarly, we let $\mathcal{U}(P)$ denote the submonoid of $P$ generated by the set of unstable atoms. Clearly, $P$ is stable (resp., unstable) if and only if $P = \mathcal{S}(P)$ (resp., $P = \mathcal{U}(P)$). In addition, $P = \mathcal{S}(P) + \mathcal{U}(P)$, and $\mathcal{S}(P) \cap \mathcal{U}(P)$ is trivial only when either $\mathcal{S}(P)$ or $\mathcal{U}(P)$ is trivial. Clearly, if $P$ is stable, then it cannot be finitely generated. Finally, we say that $u \in \mathcal{U}(P)$ is \emph{absolutely unstable} provided that $u$ is not divisible by any stable atom in $P$, and we let $\mathcal{U}^a(P)$ denote the set of all absolutely unstable elements of~$P$.

\begin{example}
	Let $(p_n)_{n \in \mathbb{N}}$ be the strictly increasing sequence with underlying set $\pp \setminus \{2\}$, and consider the prime reciprocal monoid $P$ defined as
	\[
		P := \bigg\langle \frac{3 + (-1)^n}{p_{2n-1}}, \frac{p_{2n} - 1}{p_{2n}} \ : \ n \in \nn \bigg\rangle.
	\]
	Set $a_n = \frac{3 + (-1)^n}{p_{2n-1}}$ and $b_n =  \frac{p_{2n} - 1}{p_{2n}}$. One can readily verify that $P$ is an atomic monoid with $\mathcal{A}(P) = \{a_n, b_n : n \in \nn\}$. As both sets
	\[
		\{n \in \nn : \mathsf{n}(a_n) = 2\} \quad \text{ and } \quad \{n \in \nn : \mathsf{n}(a_n) = 4\}
	\]
	have infinite cardinality, $a_n$ is a stable atom for every $n \in \nn$. In addition, since $(\mathsf{n}(b_n))_{n \in \mathbb{N}}$ is a strictly increasing sequence bounded below by $\mathsf{n}(b_1) = 4$ and $\mathsf{n}(a_n) \in \{2,4\}$, we have that $b_n$ is an unstable atom for every $n \in \nn_{\ge 2}$. Also, notice that $4/3 = 2a_1 \in \mathcal{S}(P)$, but $4/3 \notin \mathcal{U}(P)$ because $\mathsf{d}(4/3) = 3 \notin \mathsf{d}(\mathcal{U}(P))$. Furthermore, for every $n \in \nn$ the element $u_n := (p_{2n} - 1) b_n \in \mathcal{U}(P)$ is not in $\mathcal{S}(P)$ because $p_{2n} = \mathsf{d}(u_n) \notin \mathsf{d}(\mathcal{S}(P))$. However, $\mathcal{S}(P) \cap \mathcal{U}(P) \ne \emptyset$ since the element $4 = 6 a_1 = 5 b_1$ belongs to both $\mathcal{S}(P)$ and $\mathcal{U}(P)$. Finally, we claim that $2b_n$ is absolutely unstable for every $n \in \nn$. If this were not the case, then $2b_k \notin \mathcal{M}(P)$ for some $k \in \nn$. By Proposition~\ref{prop:sufficient condition for molecules in PPM} there exists $a \in \mathcal{A}(P)$ such that $\mathsf{m}(a, 2b_k) \ge \mathsf{d}(a)$. In this case, one would obtain that $2b_k \ge \mathsf{m}(a,2b_k)a \ge \mathsf{d}(a)a = \mathsf{n}(a) \ge 2$, contradicting that $b_n < 1$ for every $n \in \nn$. Thus, $2b_n \in \mathcal{U}^a(P)$ for every $n \in \nn$.
\end{example}

\begin{proposition} \label{prop: characterization of factorial elements in stable PPM}
	Let $P$ be a prime reciprocal monoid that is stable, and let $x \in P$. Then $x \in \mathcal{M}(P)$ if and only if~$\mathsf{n}(a)$ does not divide $x$ in $P$ for any $a \in \mathcal{A}(P)$.
\end{proposition}

\begin{proof}
	For the direct implication, assume that $x \in \mathcal{M}(P)$ and suppose, by way of contradiction, that $\mathsf{n}(a) \mid_P x$ for some $a \in \mathcal{A}(P)$. Since $a$ is a stable atom, there exist $p_1, p_2 \in \pp$ with $p_1 \neq p_2$ such that $\gcd(p_1 p_2, \mathsf{n}(a)) = 1$ and $\mathsf{n}(a)/p_1, \mathsf{n}(a)/p_2 \in \mathcal{A}(P)$. As $\mathsf{n}(a) \mid_P x$, we can take $a_1, \dots, a_k \in \mathcal{A}(P)$ such that $x = \mathsf{n}(a) + a_1 + \dots + a_k$. Therefore
	\[
		p_1 \frac{\mathsf{n}(a)}{p_1} + a_1 + \dots + a_k \ \text{ and } \ p_2 \frac{\mathsf{n}(a)}{p_2} + a_1 + \dots + a_k
	\]
	are two distinct factorizations in $\mathsf{Z}(x)$, contradicting that $x$ is a molecule.
	Conversely, suppose that $x$ is not a molecule. Consider two distinct factorizations $z := \sum_{i=1}^k \alpha_i a_i$ and $z' := \sum_{i=1}^k \beta_i a_i$ in $\mathsf{Z}(x)$, where $k \in \nn$, $\alpha_i, \beta_i \in \nn_0$, and $a_1, \dots, a_k \in \mathcal{A}(P)$ are pairwise distinct atoms. Pick an index $j \in \{1, \dots, k\}$ such that $\alpha_j \neq \beta_j$ and assume, without loss of generality, that $\alpha_j < \beta_j$. After applying the $\mathsf{d}(a_j)$-adic valuations on both sides of the equality
	\[
		\sum_{i=1}^k \alpha_i a_i = \sum_{i=1}^k \beta_i a_i
	\]
	one finds that the prime $\mathsf{d}(a_j)$ divides $\beta_j - \alpha_j$. Therefore $\beta_j > \mathsf{d}(a_j)$ and so
	\[
		x = \mathsf{n}(a_j) + (\beta_j - \mathsf{d}(a_j))a_j + \sum_{i \neq j} \alpha_i a_i.
	\vspace{-6pt}
	\]
	Hence $\mathsf{n}(a_j) \mid_P x$, which concludes the proof.
\end{proof}

Observe that the reverse implication of Proposition~\ref{prop: characterization of factorial elements in stable PPM} does not require $\mathcal{S}(P) = P$. However, the stability of $P$ is required for the direct implication to hold as the following example illustrates.

\begin{example}
	Let $(p_n)_{n \in \mathbb{N}}$ be the strictly increasing sequence with underlying set $\pp \setminus \{2\}$, and consider the unstable prime reciprocal monoid
	\[
		P := \bigg\langle \frac 12, \frac{p_n^2 -1}{p_n} \, : \, n \in \nn \bigg\rangle.
	\]
	Because the smallest two atoms of $P$ are $1/2$ and $8/3$, it immediately follows that $m := 2(1/2) + 8/3 \notin \langle 1/2 \rangle$ must be a molecule of $P$. In addition, notice that $1 = \mathsf{n}(1/2)$ divides $m$ in $P$.
\end{example}

We conclude this section characterizing the molecules of prime reciprocal monoids.

\begin{theorem} \label{thm:characterization of factorial elements in PPM}
	Let $P$ be a prime reciprocal monoid. Then $x \in P$ is a molecule if and only if $x = s + u$ for some $s \in \mathcal{S}(P) \cap \mathcal{M}(P)$ and $u \in \mathcal{U}^a(P) \cap \mathcal{M}(P)$.
\end{theorem}

\begin{proof}
	First, suppose that $x$ is a molecule. As $P = \mathcal{S}(P) + \mathcal{U}(P)$, there exist $s \in \mathcal{S}(P)$ and $u \in \mathcal{U}(P)$ such that $x = s + u$. The fact that $x \in \mathcal{M}(P)$ guarantees that $s,u \in \mathcal{M}(P)$. On the other hand, since $|\mathsf{Z}(u)| = 1$ and $u$ can be factored using only unstable atoms, $u$ cannot be divisible by any stable atom in $P$. Thus, $u \in \mathcal{U}^a(P)$, and the direct implication follows. 
	
	For the reverse implication, assume that $x = s + u$, where $s \in \mathcal{S}(P) \cap \mathcal{M}(P)$ and $u \in \mathcal{U}^a(P) \cap \mathcal{M}(P)$. We first check that $x$ can be uniquely expressed as a sum of two elements $s$ and $u$ contained in the sets $\mathcal{S}(P) \cap \mathcal{M}(P)$ and $\mathcal{U}^a(P) \cap \mathcal{M}(P)$, respectively. To do this, suppose that $x = s + u = s' + u'$, where $s' \in \mathcal{S}(P) \cap \mathcal{M}(P)$ and $u' \in \mathcal{U}^a(P) \cap \mathcal{M}(P)$. Take pairwise distinct stable atoms $a_1, \dots, a_k$ of $P$ for some $k \in \nn$ such that $z = \sum_{i=1}^k \alpha_i a_i \in \mathsf{Z}_P(s)$ and $z' = \sum_{i=1}^k \alpha'_i a_i \in \mathsf{Z}_P(s')$, where $\alpha_j, \alpha'_j \in \nn_0$ for $j = 1, \dots, k$. Because $u$ and $u'$ are absolutely unstable elements, they are not divisible in $P$ by any of the atoms $a_i$'s. Thus, $\mathsf{d}(a_j) \nmid \mathsf{d}(u)$ and $\mathsf{d}(a_j) \nmid \mathsf{d}(u')$ for any $j \in \{1,\dots, k\}$. Now for each $j = 1,\dots,k$ we can apply the $\mathsf{d}(a_j)$-adic valuation in both sides of the equality
	\[
		u + \sum_{i=1}^k \alpha_i a_i  = u' + \sum_{i=1}^k \alpha'_i a_i 
	\]
	to conclude that the prime $\mathsf{d}(a_j)$ must divide $\alpha_j - \alpha'_j$. Therefore either $z = z'$ or there exists $j \in \{1,\dots,k\}$ such that $|\alpha_j - \alpha'_j| > \mathsf{d}(a_j)$. Suppose that $|\alpha_j - \alpha'_j| > \mathsf{d}(a_j)$ for some $j$, and say $\alpha_j > \alpha'_j$. As $\alpha_j > \mathsf{d}(a_j)$, one can replace $\alpha_j a_j$ by $(\alpha_j - \mathsf{d}(a_j))a_j + \mathsf{n}(a_j)$ in $s = \phi(z) = \alpha_1 a_1 + \dots + \alpha_k a_k$ to find that $\mathsf{n}(a_j)$ divides $s$ in $\mathcal{S}(P)$, which contradicts Proposition~\ref{prop: characterization of factorial elements in stable PPM}. Then we have $z = z'$. Therefore $s' = s$ and $u' = u$.
	
	Finally, we argue that $x \in \mathcal{M}(P)$. Write $x = \sum_{i=1}^{\ell} \gamma_i a_i + \sum_{i=1}^\ell \beta_i b_i$ for $\ell \in \nn_{\ge k}$, pairwise distinct stable atoms $a_1, \dots, a_\ell$ (where $a_1, \dots, a_k$ are the atoms showing up in $z$), pairwise distinct unstable atoms $b_1, \dots, b_\ell$, and coefficients $\gamma_i, \beta_i \in \nn_0$ for every $i=1,\dots,\ell$. Set $z''' := \sum_{i=1}^{\ell} \gamma_i a_i$ and $w''' = \sum_{i=1}^\ell \beta_i b_i$. Note that, \emph{a priori}, $\phi(z''')$ and $\phi(w''')$ are not necessarily molecules.  As in the previous paragraph, we can apply $\mathsf{d}(a_j)$-adic valuation to both sides of the equality
	\[
		u + \sum_{i=1}^k \alpha_i a_i = \sum_{i=1}^{\ell} \gamma_i a_i + \sum_{i=1}^\ell \beta_i b_i
	\]
	to find that $z''' = z$. Hence $\phi(z''') = s$ and $\phi(w''') = u$ are both molecules. Therefore $z'''$ must be the unique factorization of $s$, while $w'''$ must be the unique factorization of $u$. As a result, $x \in \mathcal{M}(P)$.
\end{proof}
\medskip

%% file: tex/ch6.tex
\chapter{Puiseux algebras} \label{algebras}

\section{Monoid Algebras}
\label{sec:molecules of PA}

Let $M$ be a monoid and let $R$ be a commutative ring with identity. Then $R[X;M]$ denotes the ring of all functions $f \colon M \to R$ having finite \emph{support}, which means that $\supp(f) := \{s \in M : f(s) \neq 0 \}$ is finite. We represent an element $f \in R[X;M]$ by
\[
	f(X) = \sum_{i=1}^n f(s_i)X^{s_i},
\]
where $s_1, \dots, s_n$ are the elements in $\supp(f)$. The ring $R[X;M]$ is called the \emph{monoid ring} of $M$ \emph{over}~$R$, and the monoid $M$ is called the \emph{exponent monoid} of $R[X;M]$. For a field $F$, we will say that $F[X;M]$ is a \emph{monoid algebra}. As we are primarily interested in the molecules of monoid algebras of Puiseux monoids, we introduce the following definition.

\begin{definition}
	If $F$ is a field and $P$ is a Puiseux monoid, then we say that $F[X;P]$ is a \emph{Puiseux algebra}. If $N$ is a numerical monoid, then $F[X;N]$ is called a \emph{numerical monoid algebra}.
\end{definition}

Let $F[X;P]$ be a Puiseux algebra. We write any element $f \in F[X;P] \setminus \{0\}$ in \emph{canonical representation}, that is, $f(X) = \alpha_1 X^{q_1} + \dots + \alpha_k X^{q_k}$ with $\alpha_i \neq 0$ for every $i = 1, \dots, k$ and $q_1 > \dots > q_k$. It is clear that any element of $F[X;P] \setminus \{0\}$ has a unique canonical representation. In this case, $\deg(f) := q_1$ is called the \emph{degree} of $f$, and we see that the degree identity $\deg(fg) = \deg(f) + \deg(g)$ holds for all $f, g \in F[X;P] \setminus \{0\}$. As for polynomials, we say that $f$ is a \emph{monomial} if $k = 1$. It is not hard to verify that $F[X;P]$ is an integral domain with set of units $F^\times$, although this follows from~\cite[Theorem~8.1]{rG84} and~\cite[Theorem~11.1]{rG84}. Finally, note that, unless $P \cong (\nn_0,+)$, no monomial of $F[X;P]$ can be a prime element; this is a consequence of the trivial fact that non-cyclic Puiseux monoids do not contain prime elements. Puiseux algebras have been considered in~\cite{ACHZ07,CG19,fG18a,fG19e}.

For an integral domain $R$, we let $R_{\text{red}}$ denote the reduced monoid of the multiplicative monoid of $R$.

\section{Factorial Elements of Puiseux Algebras}

We proceed to study the factorial elements of a given Puiseux algebra.

\begin{definition}
	Let $R$ be an integral domain. We call a nonzero non-unit $r \in R$ a \emph{molecule} if $rR^\times$ is a molecule of $R_{\text{red}}$.
\end{definition}

Let $R$ be an integral domain. By simplicity, we let $\mathcal{A}(R)$, $\mathcal{M}(R)$, $\mathsf{Z}(R)$, and $\phi_R$ denote $\mathcal{A}(R_{\text{red}})$, $\mathcal{M}(R_{\text{red}})$, $\mathsf{Z}(R_{\text{red}})$, and $\phi_{R_{\text{\text{red}}}}$, respectively. In addition, for a nonzero non-unit $r \in R$, we let $\mathsf{Z}_R(r)$ and $\mathsf{L}_R(r)$ denote $\mathsf{Z}_{R_{\text{red}}}(rR^\times)$ and $\mathsf{L}_{R_{\text{red}}}(rR^\times)$, respectively.

\begin{proposition} \label{prop:monoid molecules}
	Let $F$ be a field, and let $P$ be a Puiseux monoid. For a nonzero $\alpha \in F$, a monomial $X^q \in \mathcal{M}(F[X;P])$ if and only if $q \in \mathcal{M}(P)$.
\end{proposition}

\begin{proof}
	Consider the canonical monoid monomorphism $\mu \colon P \to F[X;P] \setminus \{0\}$ given by $\mu(q) = X^q$. It follows from~\cite[Lemma~3.1]{CM11} that an element $a \in P$ is an atom if and only if the monomial $X^a$ is irreducible in $F[X;P]$ (or, equivalently, an atom in the reduced multiplicative monoid of $F[X;P]$). Therefore $\mu$ lifts canonically to the monomorphism $\bar{\mu} \colon \mathsf{Z}(P) \to \mathsf{Z}(F[X;P])$ determined by the assignments $a \mapsto X^a$ for each $a \in \mathcal{A}(P)$, preserving not only atoms but also factorizations of the same element. Put formally, this means that the diagram
	\[
		\begin{CD}
		\mathsf{Z}(P)      		 @>\bar{\mu}>>       \mathsf{Z}(F[X;P])        \\
		@V\phi_PVV                							 @V\phi_{F[X;P]}VV    	   \\
		P    							@>\mu>>          	  F[X;P]_{\text{red}}         \\     
		\end{CD}
	\]
	commutes, and the (fiber) restriction maps $\bar{\mu}_q \colon \mathsf{Z}_P(q) \to \mathsf{Z}_{F[X;P]}(X^q)$ of $\bar{\mu}$ are bijections for every $q \in P$. Hence $|\mathsf{Z}_P(q)| = 1$ if and only if $|\mathsf{Z}_{F[X;P]}(X^q)| = 1$ for all $q \in P^\bullet$, which concludes our proof.
\end{proof}

\begin{corollary} \label{cor:PA with infinitely more monomial molecules than atoms}
	For each field $F$, there exists an atomic Puiseux monoid $P$ whose Puiseux algebra satisfies that $|\mathcal{M}(F[X;P]) \setminus \mathcal{A}(F[X;P])| = \infty$.
\end{corollary}

\begin{proof}
	It is an immediate consequence of Proposition~\ref{prop:PPM has infinitely many molecules that are not atoms} and Proposition~\ref{prop:monoid molecules}.
\end{proof}

An element $x \in \gp(M)$ is called a \emph{root element} if it is contained in the root closure of $M$, i.e., $x \in \widetilde{M}$. Before providing a characterization for the irreducible elements of $F[X;P]$, let us argue the following two easy lemmas.

\begin{lemma} \label{lem:denominator sets of PM are closed under LCM}
	Let $P$ be a Puiseux monoid. Then $\mathsf{d}(P^\bullet)$ is closed under taking least common multiples.
\end{lemma}

\begin{proof}
	Take $d_1, d_2 \in \mathsf{d}(P^\bullet)$ and $q_1, q_2 \in P^\bullet$ with $\mathsf{d}(q_1) = d_1$ and $\mathsf{d}(q_2) = d_2$. Now set $d = \gcd(d_1, d_2)$ and $n = \gcd(\mathsf{n}(q_1), \mathsf{n}(q_2))$. It is clear that $n$ is the greatest common divisor of $(d_2/d) \mathsf{n}(q_1)$ and $(d_1/d) \mathsf{n}(q_2)$. So there exist $m \in \nn$ and $c_1, c_2 \in \nn_0$ such that
	\begin{equation} \label{eq:gcd}
		n \big(1 + m \, \lcm(d_1, d_2) \big) = c_1 \frac{d_2}{d} \mathsf{n}(q_1) + c_2 \frac{d_1}{d} \mathsf{n}(q_2).
	\end{equation}
	Using the fact that $d \, \lcm(d_1,d_2) = d_1 d_2$, one obtains that
	\[
		\frac{n \big(1 + m \, \lcm(d_1, d_2) \big)}{\lcm(d_1,d_2)} = c_1 q_1 + c_2 q_2 \in P
	\]
	after dividing both sides of the equality~(\ref{eq:gcd}) by $\lcm(d_1,d_2)$. In addition, note that $n (1 + m \, \lcm(d_1, d_2) )$ and $\lcm(d_1,d_2)$ are relatively prime. Hence $\lcm(d_1, d_2) \in \mathsf{d}(P^\bullet)$, from which the lemma follows.
\end{proof}

\begin{lemma} \label{lem:denominator division in a root-closed monoid}
	Let $P$ be a root-closed Puiseux monoid containing $1$. Then $1/d \in P$ for all $d \in \mathsf{d}(P^\bullet)$.
\end{lemma}

\begin{proof}
	Let $d \in \mathsf{d}(P^\bullet)$, and take $r \in P^\bullet$ such that $\mathsf{d}(r) = d$. As $\gcd(\mathsf{n}(r), \mathsf{d}(r)) = 1$, there exist $a,b \in \nn_0$ such that $a \mathsf{n}(r) - b \, \mathsf{d}(r) = 1$. Therefore
	\[
		\frac{1}{d} = \frac{a \mathsf{n}(r) - b \, \mathsf{d}(r)}{d} = a r - b \in \gp(P).
	\]
	This, along with the fact that $d (1/d) = 1 \in P$, ensures that $1/d$ is a root element of $P$. Since $P$ is root-closed, it must contain $1/d$, which concludes our argument.
\end{proof}

We are in a position now to characterize the irreducibles of $F[X;P]$.

\begin{proposition} \label{prop:irreducibles in a PA of a root closed PM}
	Let $F$ be a field, and let $P$ be a root-closed Puiseux monoid containing $1$. Then $f \in F[X;P] \setminus F$ is irreducible in $F[X;P]$ if and only if $f(X^m)$ is irreducible in $F[X]$ for every $m \in \mathsf{d}(P^\bullet)$ that is a common multiple of the elements of $\mathsf{d}(\supp(f))$.
\end{proposition}

\begin{proof}
	Suppose first that $f \in F[X;P] \setminus F$ is an irreducible element of $F[X;P]$, and let $m \in \mathsf{d}(P^\bullet)$ be a common multiple of the elements of $\mathsf{d}\big(\supp(f) \big)$. Then $f(X^m)$ is an element of $F[X]$. Take $g, h \in F[X]$ such that $f(X^m) = g(X) \, h(X)$. As $P$ is a root-closed and $m \in \mathsf{d}(P^\bullet)$, Lemma~\ref{lem:denominator division in a root-closed monoid} ensures that $g(X^{1/m}), h(X^{1/m}) \in F[X;P]$. Thus, $f(X) = g(X^{1/m}) h(X^{1/m})$ in $F[X;P]$. Since $f$ is irreducible in $F[X;P]$ either $g(X^{1/m}) \in F$ or $h(X^{1/m}) \in F$, which implies that either $g \in F$ or $h \in F$. Hence $f(X^m)$ is irreducible in $F[X]$.
	
	Conversely, suppose that $f \in F[X;P]$ satisfies that $f(X^m)$ is an irreducible polynomial in $F[X]$ for every $m \in \mathsf{d}(P^\bullet)$ that is a common multiple of the elements of the set $\mathsf{d}(\supp(f))$. To argue that $f$ is irreducible in $F[X;P]$ suppose that $f = g \, h$ for some $g, h \in F[X;P]$. Let $m_0$ be the least common multiple of the elements of $\mathsf{d}(\supp(g)) \cup \mathsf{d}(\supp(h))$. Lemma~\ref{lem:denominator sets of PM are closed under LCM} guarantees that $m_0 \in \mathsf{d}(P^\bullet)$. Moreover, $f = g \, h$ implies that $m_0$ is a common multiple of the elements of $\mathsf{d}(\supp(f))$. As a result, the equality $f(X^{m_0}) = g(X^{m_0})h(X^{m_0})$ holds in $F[X]$. Since $f(X^{m_0})$ is irreducible in $F[X]$, either $g(X^{m_0}) \in F$ or $h(X^{m_0}) \in F$ and, therefore, either $g \in F$ or $h \in F$. This implies that $f$ is irreducible in $F[X;P]$, as desired.
\end{proof}

We proceed to show the main result of this section.

\begin{theorem} \label{thm:U-UFD Puiseux algebras}
	Let $F$ be a field, and let $P$ be a root-closed Puiseux monoid. Then
	\[
		\mathcal{M}(F[X;P]) = \langle \mathcal{A}(F[X;P]) \rangle.
	\]
\end{theorem}

\begin{proof}
	As each molecule of $F[X;P]$ is a product of irreducible elements in $F[X;P]$, the inclusion $\mathcal{M}(F[X;P]) \subseteq \langle \mathcal{A}(F[X;P]) \rangle$ holds trivially. For the reverse inclusion, suppose that $f \in F[X;P] \setminus F$ can be written as a product of irreducible elements in $F[X;P]$. As a result, there exist $k,\ell \in \nn$ and irreducible elements $g_1, \dots, g_k$ and $h_1, \dots, h_\ell$ in $F[X;P]$ satisfying that
	\begin{equation} \label{eq:two products of irreducibles}
		g_1(X) \cdots g_k(X) = f(X) = h_1(X) \cdots h_\ell(X).
	\end{equation}
	Let $m$ be the least common multiple of all the elements of the set
	\[
		\bigg(\bigcup_{i=1}^k\mathsf{d}\big(\supp(g_i)\big) \bigg) \bigcup  \bigg(\bigcup_{j=1}^\ell\mathsf{d}\big(\supp(h_j)\big) \bigg).
	\]
	Note that $f(X^m)$, $g_i(X^m)$ and $h_j(X^m)$ are polynomials in $F[X]$ for $i = 1, \dots, k$ and $j = 1, \dots, \ell$. Lemma~\ref{lem:denominator sets of PM are closed under LCM} ensures that $m \in \mathsf{d}(P^\bullet)$. On the other hand, $m$ is a common multiple of all the elements of $\mathsf{d}(\supp(g_i))$ (or all the elements of $\mathsf{d}(\supp(h_i))$). Therefore Proposition~\ref{prop:irreducibles in a PA of a root closed PM} guarantees that the polynomials $g_i(X^m)$ and $h_j(X^m)$ are irreducible in $F[X]$ for $i=1,\dots,k$ and $j=1,\dots,\ell$. After substituting $X$ by $X^m$ in~(\ref{eq:two products of irreducibles}) and using the fact that $F[X]$ is a UFD, one finds that $\ell = k$ and $g_i(X^m) = h_{\sigma(i)}(X^m)$ for some permutation $\sigma \in S_k$ and every $i = 1, \dots, k$. This, in turns, implies that $g_i = h_{\sigma(i)}$ for $i = 1, \dots, k$. Hence $|\mathsf{Z}_{F[X;P]}(f)| = 1$, which means that $f$ is a molecule of $F[X;P]$.
\end{proof}

As we have seen before, Corollary~\ref{cor:PA with infinitely more monomial molecules than atoms} guarantees the existence of a Puiseux algebra $F[X;P]$ satisfying that $|\mathcal{M}(F[X;P]) \setminus \mathcal{A}(F[X;P])| = \infty$. Now we use Theorem~\ref{thm:U-UFD Puiseux algebras} to construct an infinite class of Puiseux algebras satisfying a slightly more refined condition.

\begin{proposition} \label{prop:U-UFD Puiseux algebras}
	For any field $F$, there exist infinitely many Puiseux monoids $P$ such that the algebra $F[X;P]$ contains infinitely many molecules that are neither atoms nor monomials.
\end{proposition}

\begin{proof}
	Let $(p_j)_{j \in \mathbb{N}}$ be the strictly increasing sequence with underlying set $\pp$. Then for each $j \in \nn$ consider the Puiseux monoid $P_j = \langle 1/p^n_j : n \in \nn \rangle$. Fix $j \in \nn$, and take $P := P_j$. The fact that $\gp(P) = P \cup -P$ immediately implies that $P$ is a root-closed Puiseux monoid containing~$1$. Consider the Puiseux algebra $\qq[X;P]$ and the element $X + p \in \qq[X;P]$, where $p \in \pp$. To argue that $X + p$ is an irreducible element in $\qq[X;P]$, write $X + p = g(X) \, h(X)$ for some $g, h \in \qq[X;P]$. Now taking $m$ to be the maximum power of $p_j$ in the set $\mathsf{d}( \supp(g) \cup \supp(h) )$, one obtains that $X^m + p = g(X^m) \, h(X^m)$ in $\qq[X]$. Since $\qq[X]$ is a UFD, it follows by Eisenstein's criterion that $X^m + p$ is irreducible as a polynomial over $\qq$. Hence either $g(X) \in \qq$ or $h(X) \in \qq$, which implies that $X + p$ is irreducible in $\qq[X;P]$. Now it follows by Theorem~\ref{thm:U-UFD Puiseux algebras} that $(X + p)^n$ is a molecule in $\qq[X;P]$ for every $n \in \nn$. Clearly, the elements $(X + p)^n$ are neither atoms nor monomials.
	
	Finally, we prove that the algebras we have defined in the previous paragraph are pairwise non-isomorphic. To do so suppose, by way of contradiction, that $\qq[X;P_j]$ and $\qq[X;P_k]$ are isomorphic algebras for distinct $j,k \in \nn$. Let $\psi \colon \qq[X;P_j] \to \qq[X;P_k]$ be an algebra isomorphism. Since $\psi$ fixes $\qq$, it follows that $\psi(X^q) \notin \qq$ for any $q \in P_j^\bullet$. This implies that $\deg(\psi(X)) \in P_k^\bullet$. As $\mathsf{d}(P_j^\bullet)$ is unbounded there exists $n \in \nn$ such that $p_j^n > \mathsf{n}( \deg( \psi(X) ) )$. Observe that
	\begin{align} \label{eq: degrees 1}
		\deg \big( \psi(X) \big) = \deg \big( \psi \big( X^{\frac{1}{p_j^n}} \big)^{p_j^n} \big) = p_j^n \deg \big( \psi \big( X^{\frac{1}{p_j^n}} \big) \big).
	\end{align}
	Because $\gcd(p_j,d) = 1$ for every $d \in \mathsf{d}(P_k^\bullet)$, from~(\ref{eq: degrees 1}) one obtains that $p_j^n$ divides $\mathsf{n}( \deg \psi (X) )$, which contradicts that $p_j^n > \mathsf{n}( \deg( \psi(X) ) )$. Hence the Puiseux algebras in $\{P_j : j \in \nn\}$ are pairwise non-isomorphic, which completes our proof.
\end{proof}
\medskip

%% file: bio.tex
\biography{%
	Marly received her Ph.D. in mathematics from the University of Florida in 2019. Her general area of research lies in Factorization Theory and Commutative Algebra. In particular, she has studied the atomic structure and factorization theory of Puiseux monoids. By the time she graduated, she had published eighth research papers and presented her research in more than ten conferences.
	
	Although she pursued a degree in pure mathematics, she is passionate about applied fields and technology. She interned for RStudio doing mathematical modeling and worked for more than two years as an Application Analyst for the College of Medicine at the University of Florida.
}